\documentclass[11pt]{amsart}
\usepackage[text={475pt,660pt},centering]{geometry}
\usepackage[english]{babel}
\usepackage[utf8]{inputenc}
\usepackage{amsmath,amssymb,amsthm,mathrsfs}
\usepackage{graphicx,caption,subcaption}
\usepackage{verbatim}
\usepackage{hyperref}
\hypersetup{
    colorlinks=true,
    linkcolor=blue,
    filecolor=magenta,      
    urlcolor=cyan,
    citecolor=cyan 
}
\usepackage{cleveref}
\numberwithin{equation}{section}
\usepackage{changepage}
\usepackage{empheq}
\usepackage{enumerate}

\binoppenalty=\maxdimen
\relpenalty=\maxdimen

\theoremstyle{plain}
\newtheorem{theorem}{Theorem}[section]
\newtheorem{corollary}{Corollary}[section]
\newtheorem{lemma}{Lemma}[section]
\newtheorem{proposition}{Proposition}[section]

\theoremstyle{definition}
\newtheorem{question}{Question}[section]
\newtheorem{definition}{Definition}[section]

\newcommand{\Prb}{\mathbb{P}}

\newcommand{\Expt}{\mathbb{E}}
\newcommand{\V}{\operatorname{Var}}
\newcommand{\Vc}{\mathbb{V}}
\newcommand{\GO}{\mathcal{O}}
\newcommand{\F}{\mathcal{F}}
\newcommand{\X}{\mathcal{X}}
\newcommand{\Y}{\mathcal{Y}}
\newcommand{\Z}{\mathcal{Z}}

\newcommand{\Ind}[1]{\mathbf{1}_{\left\{#1\right\}}}

\newcommand{\Rd}{\mathbb{R}^d}
\newcommand{\Zd}{\mathbb{Z}^d}
\newcommand{\cu}{\square}
\newcommand{\dist}{\operatorname{dist}}

\newcommand{\p}[1]{\partial_{x_{#1}}}
\newcommand{\supp}{\operatorname{supp}}
\newcommand{\bracket}[1]{\left\langle{#1}\right\rangle}
\newcommand{\norm}[1]{\left\Vert{#1}\right\Vert}
\newcommand{\edge}[1]{\{#1\}}

\newcommand{\D}{\mathcal{D}}
\newcommand{\Dr}[1]{\D_{\e_{#1}}}
\newcommand{\e}{e}

\newcommand{\ub}{\bar{u}}
\newcommand{\vb}{\bar{v}}

\newcommand{\aL}{\underline{L}}

\newcommand{\ec}{\lambda}
\newcommand{\phim}{\phi^{(\ec)}}
\newcommand{\Sm}{\S^{(\ec)}}
\newcommand{\ab}{\bar{\a}}
\newcommand{\Fi}{\mathbf{F}}
\newcommand{\g}{\mathbf{g}}
\newcommand{\Id}{\operatorname{Id}}

\newcommand{\clt}{\mathscr{C}}
\newcommand{\slt}{\clt_s}
\newcommand{\cltm}{\clt_*(\cu_m)}
\newcommand{\ac}{\a_{\clt}}
\newcommand{\acf}[1]{\a_{\clt,{#1}}}
\newcommand{\Tcu}{\mathcal{T}}
\newcommand{\Scu}{\mathcal{S}}
\newcommand{\Pcu}{\mathcal{P}}
\newcommand{\Gcu}{\mathcal{G}}
\newcommand{\M}{\mathcal{M}}
\newcommand{\A}{\mathcal{A}}
\newcommand{\coa}[1]{\left[#1\right]}
\newcommand{\Ed}{E_d}
\newcommand{\Eda}{E_d^{\a}}
\newcommand{\pc}{\mathfrak{p}_c}
\newcommand{\pp}{\mathfrak{p}}
\newcommand{\itr}{\mathrm{int}}
\newcommand{\rt}{\mathbf{r}}
\newcommand{\size}{\operatorname{size}}
\newcommand{\diam}{\operatorname{diam}}
\newcommand{\clp}{\operatorname{cl}_{\Pcu}}
\newcommand{\zb}{\operatorname{\bar{z}}}

\renewcommand{\leq}{\leqslant}
\renewcommand{\geq}{\geqslant}
\renewcommand{\epsilon}{\varepsilon}
\renewcommand{\S}{\mathbf{S}}

\renewcommand{\a}{\mathbf{a}}

\theoremstyle{remark}

\newtheorem*{remark}{Remark}

\author[Chenlin GU]{Chenlin GU}
\address[Chenlin GU]{DMA, Ecole Normale Supérieure, PSL University, Paris, France}
\email{chenlin.gu@ens.fr}
\title{An efficient algorithm for solving elliptic problems on percolation clusters}

\begin{document}

\begin{abstract}
 We present an efficient algorithm to solve elliptic Dirichlet problems defined on the cluster of $\Zd$ supercritical Bernoulli percolation, as a generalization of the iterative method proposed by S. Armstrong, A. Hannukainen, T. Kuusi and J.-C. Mourrat. We also explore the two-scale expansion on the infinite cluster of percolation, and use it to give a rigorous analysis  of the algorithm. 
\end{abstract}

\maketitle

\hspace{3cm}

\begin{figure}[h]
\centering
\includegraphics[scale=0.5]{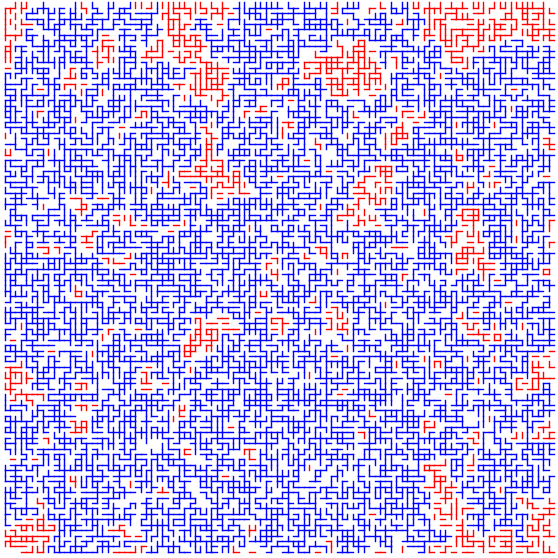}
\caption{A simulation of $2D$ Bernoulli bond percolation with $\pp = 0.51$ in a cube $\cu$ of size $100 \times 100$. The cluster in blue is the maximal cluster $\clt_*(\cu)$ while the clusters in red are the other small ones.}\label{Fig:cluster}
\end{figure}

\newpage

\tableofcontents

\section{Introduction}
\subsection{Motivation and main result}

The main goal of this paper is to study a fast algorithm for computing the solution of Dirichlet problems with random coefficients on percolation clusters. We consider a percolation model on the space $(\Zd, \Ed)$ for dimension $d \geq 2$, where $\Ed$ denotes the set of (unoriented) nearest-neighbor \emph{bonds} (or \emph{edges}), that is, two-element sets $\{x,y\}$ with $x,y \in \Zd$ satisfying $|x-y| = 1$. We also write $x \sim y$ whenever $\{x,y\} \in E_d$. More precisely, we give ourselves a constant $\Lambda > 1$ and a random field $\a : \Ed \rightarrow \{0\} \cup [\Lambda^{-1}, 1]$ such that the random variables $\{\a(e)\}_{e \in \Ed}$ are independent and identically distributed. We may refer to $\a(e)$ as the \emph{conductance} of the bond $e$, and say that $e$ is an \textit{open bond} if $\a(e) > 0$, and that $e$ is a \textit{closed bond} otherwise. We call the connected components on $(\Zd, \Ed)$ generated by the open bonds \textit{clusters}, and we are interested in the supercritical percolation case, that is, we assume that $\pp := \Prb[\a(e) > 0]$ is strictly larger than the critical percolation parameter, which we denote by $\pc(d)$. As a consequence, there exists a unique infinite percolation cluster $\clt_{\infty}$ \cite{kesten1982percolation}. We can also define its analogue in a finite cube $\cu_m := \left(-\frac{3^m}{2}, \frac{3^m}{2} \right)^d \cap \Zd$, which is a very connected maximal cluster denoted by $\cltm$, and our goal is to find an algorithm for solving Dirichlet problems on it. That is, given two functions $f,g : \cu_m \rightarrow \mathbb{R}$, we aim to define and study an efficient method for calculating the solution $u$ of
\begin{equation}
\label{eq:main}
\left\{
	\begin{array}{ll}
	-\nabla \cdot \a \nabla u = f  & \text{ in } \cltm, \\
	u = g & \text{ on } \cltm \cap \partial \cu_m,
	\end{array}
\right.
\end{equation}
where the divergence-form operator is defined as 
$$ 
-\nabla \cdot \a \nabla u(x) = \sum_{y \sim x}\a(x,y)\left(u(x) - u(y) \right).
$$
For convenience, we will only study this problem conditionally on the event that ``$\cu_m$ is a good cube''. This event has very high probability for large $m$, since there exists a positive constant $C(d,\pp)$ such that
$$
\Prb[ \cu_m \text{ is a good cube}] \geq 1 - C(d,\pp) \exp(-C(d,\pp)^{-1}3^m).
$$
The rigorous definitions of ``$\cu_m$ is a good cube" and of the maximal cluster $\cltm$ will be given in Definitions~\ref{def:goodcube} and~\ref{def:MaximalCluster} below. Informally, one can think of $\cltm$ as the largest cluster of $\clt_{\infty} \cap \cu_m$ (see \Cref{Fig:cluster} for an illustration.)

\Cref{eq:main} is very natural to describe many models in applied mathematics and other disciplines. For example, one can think of the electric potential in a porous medium: a domain is made of two types of composites, represented respectively by the open bonds and the closed bonds on the lattice graph $(\Zd, \Ed)$, and only the open bonds are available for the current to flow, while the closed bonds are insulating. 

The complex geometry of the percolation cluster causes significant perturbations to the electric potential, and this makes efficient numerical calculations challenging. Naive finite-difference schemes will become very costly as the size of the domain is increased, and the perforated geometry and low regularity of solutions does not allow for simple coarsening mechanisms. As is well-known, using the \textit{effective conductance} $\ab$, which is a constant matrix (in fact a scalar by the symmetries in our assumptions) whose definition will be recalled in eqs.\ \eqref{eq:defEnergy} and \eqref{eq:defEffectiveConductance}, one can replace the heterogeneous operator $-\nabla \cdot \a \nabla$ by the constant-coefficient operator $-\nabla \cdot \ab \nabla$, and thus obtain an approximation $\ub$ as the solution of a homogenized equation. This is a nice idea, but the gap between $\ub$ and $u$ always exists: on small scales, the homogenized solution $\ub$ will typically be very smooth, while $u$ has oscillations. Indeed, the homogenized solution $\ub$ can only approximate $u$ in $L^2$, but not in $H^1$. Moreover, the $L^2$ norm of $(u - \ub)$ depends on the size of $\cu_m$ and only goes to zero in the limit $m \to \infty$. In other words, $\ub$ converges to $u$ in $L^2$ only in the limit of ``infinite separation of scales".

The goal of the present work is to go beyond these limitations: we will devise an algorithm that produces a sequence of approximations which rapidly converge to $u$ in $H^1$, in a regime of large but finite separation of scales. The main idea is to look for a way to use the homogenized operator as a \emph{coarse operator} in a scheme analogous to a multigrid method. Let us introduce some more notations and state the main theorem. For any $V \subset \Zd$, the \textit{interior} of $V$ is defined as $\itr(V) := \{x\in V : y \sim x \Rightarrow y \in V\},$ and the \textit{boundary} is defined as $\partial(V) := V \backslash \itr(V).
$ The function space $C_0(V)$ is the set of the function with zero boundary condition. The $L^2$ integration of the gradient of $v$ on the percolation cluster is defined as 
$$
\norm{\nabla v \Ind{\a \neq 0}}_{L^2(V)} := \left(\frac{1}{2} \sum_{x,y \in V, x \sim y} (v(y) - v(x))^2 \Ind{\a(x,y) \neq 0} \right)^{\frac{1}{2}}.
$$

We denote the probability space by $\left(\{\a(e)\}_{e \in \Ed}, \F, \Prb\right)$, and for any $V \subset \Zd$ we denote by $\F(V) := \sigma\left( \{\a(e)\}_{e \cap V \neq \varnothing} \right)$. For a random variable $X$, we use two positive parameters $s,\theta$, and the notation $\GO$ to measure its size by 
$$X \leq \GO_s(\theta) \Longleftrightarrow \Expt\left[\exp((\theta^{-1}X)_+^s)\right] \leq 2,$$
where $(\theta^{-1} X)_+ := \max\{\theta^{-1} X, 0\}$. Roughly speaking, the statement $X \leq \GO_s(\theta)$ tells us that $X$ has a tail lighter than $\exp(-(\theta^{-1} x)^s)$. 
 We also define, for each $\ec > 0$, the mappings $\ec_{\clt,m} : \Zd \rightarrow \mathbb{R}$, and $\ell : \mathbb{R}^+ \rightarrow \mathbb{R}^+$ by
\begin{equation}\label{eq:defEcC}
\ec_{\clt,m}(x) := \left\{
	\begin{array}{ll}
	\ec & \text{if } x \in \cltm, \\	
	0 & \text{otherwise}.
	\end{array}
\right.
\qquad
\ell(\ec) := \left\{
	\begin{array}{ll}
	\log^{\frac{1}{2}} (1+\ec^{-1}) & \text{if } d = 2, \\	
	1 & \text{if }  d > 2.
	\end{array}
\right.
\end{equation} 

\begin{theorem}[Main theorem]\label{thm:main}
There exist two positive constants $s := s(d,\pp,\Lambda), C := C(d, \pp, \Lambda, s)$, and for every integer $m > 1$ and $\lambda \in \left(\frac{1}{3^m}, \frac{1}{2} \right)$, an $\F$-measurable random variable $\Z$ satisfying
$$\Z \leq \GO_s\left(C \ell(\lambda)^{\frac{1}{2}}\lambda^{\frac{1}{2}} m^{\frac{1}{s}+d}\right),$$  
such that the following holds.
Let $f,g:\cu_m  \rightarrow \mathbb{R}, u_0 \in g + C_0(\cu_m)$ and $u \in g + C_0(\cu_m)$ be the solution of \cref{eq:main}. On the event that $\cu_m$ is a good cube, for $u_1, \ub, u_2 \in C_0(\cu_m)$ solving (with null Dirichlet boundary condition)
\begin{equation}
\label{eq:iterativeCluster}
\left\{
	\begin{array}{lll}
	(\lambda^2 - \nabla \cdot \a \nabla)u_1 &= f + \nabla \cdot \a \nabla u_0  & \text{ in } \cltm\backslash \partial \cu_m ,\\
	-\nabla \cdot \ab \nabla \ub &= \ec_{\clt, m}^2 u_1  & \text{ in } \itr(\cu_m),\\
	(\lambda^2 - \nabla \cdot \a \nabla) u_2 &= (\lambda^2 - \nabla \cdot \ab \nabla) \ub & \text{ in } \cltm \backslash \partial \cu_m ,\\
	
	\end{array}
\right.
\end{equation}
and for $\hat{u}:= u_0 + u_1 + u_2$, we have the contraction estimate
\begin{equation}
\label{eq:contraction}
\norm{\nabla (\hat{u} - u) \Ind{\a \neq 0}}_{L^2(\cltm)} \leq  \Z \norm{ \nabla (u_0-u) \Ind{\a \neq 0}}_{L^2(\cltm)}.
\end{equation}
\end{theorem}

We explain a little more how this algorithm works. 
\begin{itemize}
    \item We start by an arbitrary guess $u_0 \in g + C_0(\cu_m)$ as an approximation of $u$, and repeat the \cref{eq:iterativeCluster} several rounds. At the end of every round, we use the $\hat{u}$ just obtained in place of $u_0$ in the new round of iteration. We hope that in every iteration, the error between our approximation and $u$ decreases by a multiplying factor $\Z < 1$, so we can get a desired precision once we repeat enough rounds of iteration.
    \item In fact, the random factor $\Z$ only depends on the conductance $\a$, the choice of our regularization $\ec$, and the size of the cube $\cu_m$, but it does not depend on the data $f,g,u_0$. We can choose $\ec$ such that $3^{-m} \ll \ec \ll m^{- 2\left(\frac{1}{s}+d\right)}$, then $\Z \leq \GO_s\left(C \ell(\lambda)^{\frac{1}{2}}\lambda^{\frac{1}{2}} m^{\frac{1}{s}+d}\right)$ tells us that $\Z$ has large probability to be smaller than $1$. One should think of $\lambda$ as being a small (but fixed) multiple of $m^{- 2\left(\frac{1}{s}+d\right)}$ (recall that $m$ is of the order of the logarithm of the diameter of the domain).
    \item The computational cost of each iteration is small. Indeed, we can first identify $\clt_*(\cu_m)$ once and for all by a \textit{``UnionFind"} algorithm \cite[Chapter 21]{cormen2009introduction}. The first and third steps in \cref{eq:iterative} are fast thanks to the regularization provided by the zero-order term $\ec^2$, while the second step can be done by a standard multigrid algorithm \cite{briggs2000multigrid} because it is a discrete Laplacian operator with constant coefficients. All these algorithms are not very expensive.
\end{itemize}

We can evaluate the computational complexity more precisely, assuming that the heterogeneous problems are solved by iterations of conjugate gradient descent (CGD). Denoting $r = 3^m$, we note that the spectral condition number of the operator $-\nabla \cdot \a \nabla$ is random but typically of the order of $r^2$, and thus $O(r)$ iterations of CGD are required if a direct approach is used. Meanwhile, for the operator $(\lambda^2-\nabla \cdot \a \nabla)$, the regularization helps reduce to $O\left((\log{r})^{\frac{2}{s}+2d}\right)$ iterations of CGD. Moreover, the complexity of standard multigrid and ``UnionFind" algorithms are of lower computational cost, so our iterative algorithm reduces the complexity from $O\left(r\right)$ iterations of CGD to $O\left((\log{r})^{\frac{2}{s}+2d}\right)$. Of course, in the actual numerical implementation, many improvements in intermediate steps are possible. In the rest of the paper, we focus on the convergence analysis of the method described in \Cref{thm:main}. 

We remark that \cref{eq:main} can be defined in a more general domain $\Zd \cap U_r$ where $U$ is a convex domain with $C^{1,1}$ boundary, $r>0$ is a length scale which we think of as being large, and $U_r:=\{rx | x \in U\}$. In this case $\clt_*(U_r)$ can be informally thought as the largest cluster in $U_r$. Our iterative algorithm \cref{eq:iterativeCluster} and its analysis can be adapted to this more general setting by following very similar arguments.

\subsection{Previous work}
The homogenization theory was first developed for elliptic or parabolic equations with periodic coefficients, and then generalised to the case of random stationary coefficients. There exist many classical references such as \cite{bensoussan2011asymptotic,kozlov1979averaging, tartar2009general,jikov2012homogenization,allaire1992homogenization}. Quantitative results in stochastic homogenization took a long time to emerge. The first partial results result were obtained by Yurinskii \cite{yurinskii1986averaging}. Recently, thanks to the work of Gloria, Neukamm and Otto \cite{gloria2011optimal,gloria2012optimal,gloria2014optimal,gloria2014regularity,gloria2015quantification}, and Armstrong, Kuusi, Mourrat and Smart \cite{armstrong2016lipschitz, armstrong2016mesoscopic, armstrong2016quantitative, armstrong2017additive}, we understand better the typical size of the fundamental quantities in the stochastic homogenization of uniformly elliptic equations, which provides us with the possibility to analyze the performance of numerical algorithms in this context.

The homogenization of environments that do not satisfy a uniform ellipticity condition also drew  attention. In \cite{zhikov2006homogenization}, Zhikov and Piatnitski establish many results qualitatively and explain how to formulate the effective equation on various types of degenerate stationary environments. In \cite{lamacz2015moment}, Lamacz, Neukamm and Otto obtain a bound of correctors on a simplified percolation model by imposing all the bonds in the first coordinate direction to be open. In \cite{benjamini2015disorder}, the Liouville regularity problem in a general context of random graphs is studied by Benjamini, Duminil-Copin, Kozma, and Yadin using the entropy method, and its complete description on infinite cluster of Bernoulli percolation is given by Armstrong and Dario in \cite{armstrong2018elliptic}. Dario also gives the moment estimate of the correctors of the same model in \cite{dario2018optimal}. 

Homogenization has a natural probabilistic interpertation in terms of \textit{random walks in random environment}, as a generalised central limit theorem. One fundamental work in this context is the paper \cite{kipnis1986central} by Kipnis and Varadhan, where the case of general reversible Markov chains is studied. The case of random walks on the supercritical percolation cluster attracted particular interest, and the \textit{quenched} central limit theorem was obtained by Berger, Biskup, Mathieu and Piatnitski in \cite{berger2007quenched,mathieu2007quenched}. We also refer to \cite{bisk, KLObook, kumagai} for overviews of this line of research.

Finally, concerning the construction of efficient numerical methods, our algorithm is inspired by the one introduced in \cite{armstrong2018iterative} by Armstrong, Hannukainen, Kuusi and Mourrat, which is designed to treat the same question in a uniform ellipticity context, and also \cite{gu2018uniform} where a uniform estimate is obtained. Besides the fact that the problem we consider here is not uniformly elliptic, we stress that a fundamental issue we need to address relates to the fact that the \emph{geometry} of the domain itself must be modified as we move from fine to coarse scales. Indeed, the fine scales must be resolved on the original, highly perforated domain, while the coarse scales are resolved in a homogeneous medium in which the wholes have been ``filled up''. As far as I know, this is the first work proposing a practical and rigorous method for the numerical approximation of elliptic problems posed in rapidly oscillating perforated domains. Notice that in our algorithm, we suppose that the effective conductance $\ab$ is known, because there are many excellent methods to do it quickly which can be naturally generalized in percolation setting, see for example \cite{gloria2012numerical, egloffe2014random, mourrat2016efficient,fischer2018choice, hannukainen2019computing}. Alternative numerical methods for computing the solution of elliptic problems in non-perforated domains have been studied extensively; we refer in particular to  \cite{bebendorf,BL,EGH,GGS,OZB,malpet,KY,roulette}, as well as to \cite{knap1,knap2,EL1,EL2} where the concept of homogenization is used explicitly.

\subsection{Ideas of proof and additional results}
\label{subsec:Idea}
In this part, we introduce some key concepts underlying the analysis of the algorithm and the proof of Theorem~\ref{thm:main}. We also present some useful results like estimates on the flux of the corrector, and a quantitative version of the two-scale expansion on the cluster of percolation, which are proved in this paper and are of independent interest. Some notations are explained quickly in the statement and their rigorous definitions will be given in \Cref{sec:Pre} or in the later part when they are used. 

The main strategy of the algorithm is very similar to an algorithm porposed in the previous work \cite{armstrong2018iterative, gu2018uniform} where we study the classical Dirichlet problem in $\Rd$ setting with symmetric $\mathbb{R}^{d \times d}$-valued coefficient matrix $\a$, which is random, stationary, of finite range correlation and satisfies the uniform ellipticity condition. We recall the idea in the previous work with a little abuse of notation that $\cu_m$ stands $\left(-\frac{3^m}{2}, \frac{3^m}{2}\right)^d$ in this paragraph: to solve a divergence-form equation $- \nabla \cdot \a \nabla u = f$ in $\cu_m$ with boundary condition $g$, we propose to compute $(u_1, \ub, u_2)$ with null Dirichlet boundary condition solving
\begin{equation}
\label{eq:iterative}
\left\{
	\begin{array}{lll}
	(\lambda^2 - \nabla \cdot \a \nabla)u_1 &= f + \nabla \cdot \a \nabla u_0  &  \text{ in } \cu_m ,\\
	-\nabla \cdot \ab \nabla \ub &= \ec^2 u_1  &  \text{ in } \cu_m,\\
	(\lambda^2 - \nabla \cdot \a \nabla) u_2 &= (\lambda^2 - \nabla \cdot \ab \nabla) \ub &  \text{ in } \cu_m.\\
	\end{array}
\right.
\end{equation}
In \cite{armstrong2018iterative, gu2018uniform} we proved that $\hat{u} := u_0 + u_1 + u_2$ satisfies
$$
\norm{\hat{u} - u}_{H^1(\cu_m)} \leq \Z \norm{u_0 - u}_{H^1(\cu_m)},
$$
with a random factor $\Z$ of size $\Z \leq \GO_s\left(C(\Lambda, s, d)\ell(\lambda)^{\frac{1}{2}}\lambda^{\frac{1}{2}} m^{\frac{1}{s}}\right)$ for any $s \in (0,2)$ and independent of $u,u_0,f,g$. The main ingredient in the proof is the two-scale expansion theorem: for $v, \vb$ with the same boundary condition and satisfying
\begin{equation}\label{eq:TwoScaleNaive}
(\mu^2 - \nabla \cdot \a \nabla) v = (\mu^2 - \nabla \cdot \ab \nabla) \vb  \qquad \text{ in } \cu_m,
\end{equation}
one can use $\vb + \sum_{k=1}^d  (\p{k} \vb) \phi_{\e_k}$ to approximate $v$ in $H^1$. Here $\{\e_k\}_{1 \leq k\leq d}$ stands for the canonical basis in $\Rd$, and $\phi_{\e_k}$ is \textit{the first order corrector} associated with the direction $\e_k$. In our algorithm \cref{eq:iterative}, combining the first equation, the second equation of \cref{eq:iterative} and $- \nabla \cdot \a \nabla u = f$, we can obtain that 
$$
-\nabla \cdot \ab \nabla \ub = -\nabla \cdot \a \nabla (u - u_0 - u_1)  \qquad \text{ in } \cu_m,
$$
which is an equation of type \cref{eq:TwoScaleNaive} with $\mu = 0$. Moreover, the third equation in \cref{eq:iterative} also follows the form of \cref{eq:TwoScaleNaive}, this time with $\mu = \lambda$. Thus, we have 
$$
(u - u_0 - u_1) \simeq w := \ub + \sum_{k=1}^d  (\p{k} \ub) \phi_{\e_k} \simeq u_2,
$$
up to a small error, so we can estimate $\vert\hat{u} - u\vert$ by studying
$$
\vert\hat{u} - u\vert = \vert u - (u_0 + u_1 + u_2)\vert \leq \vert (u - u_0 - u_1) - w\vert + \vert w - u_2 \vert.
$$ 
In \cite{armstrong2018iterative} the error in the two-scale expansion theorem is made quantitative, and in \cite{gu2018uniform} we refine this bound so that the contraction bound is uniform over the relevant data (most importantly: the bound is uniform over $u_0$, which guarantees that the algorithm can indeed be iterated).

It is easy to write down formal equations similar to those in \cref{eq:TwoScaleNaive} in our discrete context. However, the algorithm \cref{eq:iterativeCluster} explored in this work is not a simple adaption from $\Rd$ to the discrete setting of $(\Zd, \Ed)$. Indeed, the random geometry of the percolation cluster causes new difficulties, and finding a suitable way to handle the coarsening of this very irregular and non-uniform geometry is one of the main challenges addressed in this paper. In fact, a version of \cref{eq:TwoScaleNaive} on the cluster of percolation should be defined differently, because the cluster is determined by each realization of the random coefficient field $\a$, while the effective equation with $\ab$ is only useful if posed on a homogeneous geometry. As a consequence, a reasonable way to establish \cref{eq:TwoScaleNaive} is to define its left hand side on the cluster $\cltm$ and to define its right hand side on the homogenized geometry $\cu_m$. To overcome the obstacle when putting the two sides on two different domains, we use an operation of \textit{local mask} 
\begin{equation*}
\lambda_{\clt,m}(x) := \left\{
	\begin{array}{ll}
	\lambda & \text{if }x \in \cltm, \\	
	0 & \text{otherwise}.
	\end{array}
\right. \qquad
\acf{m}(x,y) := \left\{
	\begin{array}{ll}
	\a(x,y) & \text{if }x,y \in \cltm, \\	
	0 & \text{otherwise}.
	\end{array}
\right.
\end{equation*} 
Then a similar equation of \cref{eq:TwoScaleNaive} on $\cltm$ can be defined by  
\begin{equation}\label{eq:TwoScaleInformal}
(\ec_{\clt,m}^2 - \nabla \cdot \acf{m} \nabla) v = (\ec_{\clt,m}^2 - \nabla \cdot \ab \nabla) \vb  \qquad  \text{ in } \itr(\cu_m),
\end{equation}
and we hope to use a modified two-scale expansion 
\begin{equation}\label{eq:Expansion}
w := \vb + \sum_{k=1}^d (\Upsilon \Dr{k}\vb) \phim_{\e_k},
\end{equation}
to approximate $v$. Here $\Dr{k} \vb(x) := \vb(x + \e_k) - \vb(x)$ and $\Upsilon$ is a cut-off function supported in $\cu_m$, constant $1$ in the interior and decreases to $0$ linearly near the boundary defined as
\begin{equation}\label{eq:CutoffFunction}
\Upsilon := \Ind{\cu_m} \wedge \left(\frac{\dist(\cdot, \partial \cu_m) - \ell(\ec)}{\ell(\ec)} \right)_+,
\end{equation}
so the function $\Upsilon$ can help reduce the boundary layer effect of the two-scale expansion. The modified corrector $\{\phim_{\e_k}\}_{1 \leq k \leq d}$ is defined as  
\begin{equation}\label{eq:CorrectorModified}
\phim_{\e_k} := \phi_{\e_k} - \coa{\phi_{\e_k}}_{\Pcu}^{\eta} \star \Phi_{\ec^{-1}},
\end{equation}
where $\Phi_{\ec^{-1}}$ is a heat kernel of scale $\ec^{-1}$, i.e. $\Phi_{\ec^{-1}}(x) := \frac{1}{(4\pi \ec^{-2})^{d/2}}\exp\left(-\frac{x^2}{4\ec^{-2}}\right)$ and $\coa{\phi_{\e_k}}_{\Pcu}^{\eta}$ is a coarsened version of $\phi_{\e_k}$, whose proper definition will be given in \Cref{def:goodcube}. Although the corrector is only well-defined up to a constant, notice that \cref{eq:CorrectorModified} is well-defined. Notice also that by \eqref{eq:TwoScaleInformal}, the function $\vb$ is discrete-harmonic outside of $\cltm$.

In \Cref{sec:TwoScale}, we will prove the following quantitative two-scale expansion theorem as a main tool to prove the contraction estimate (Theorem~\ref{thm:main}). We remark here that we also add the condition $\cu_m \in \Pcu_*$, which is stronger than ``$\cu_m$ is good", and means that ``$\cltm$ is indeed a subset of $\clt_{\infty}$". 
\begin{theorem}[Two-scale expansion on percolation]\label{thm:TwoScale}
There exist two positive constants $s := s(d,\pp,\Lambda)$, $C := C(d, \pp, \Lambda, s)$, and for every integer $m > 1$ such that $\cu_m \in \Pcu_*$ and every $\lambda \in \left(\frac{1}{3^m}, \frac{1}{2} \right)$, there exists a random variable $\widetilde{\Z}$ controlled by
\begin{align*}
&\widetilde{\Z}  \leq \GO_s\left(C(d,\pp,\Lambda, s)\ell(\ec)m^{\frac{1}{s}+d}\right),
\end{align*}
such that the following is valid: for any $\mu \in [0, \ec]$ and any $v,\vb \in C_0(\cu_m)$ satisfying 
\begin{equation}\label{eq:TwoScale}
(\mu_{\clt,m}^2 - \nabla \cdot \acf{m} \nabla) v = (\mu_{\clt,m}^2 - \nabla \cdot \ab \nabla) \vb  \qquad  \text{ in } \itr(\cu_m),
\end{equation}
defining a two-scale expansion $w := \vb + \sum_{k=1}^d (\Upsilon \Dr{k}\vb) \phim_{\e_k}$, we have
\begin{equation}\label{eq:H1TwoScale}
\begin{split}
\norm{\nabla (w - v) \Ind{\a \neq 0}}_{L^2(\cltm)} &\leq \widetilde{\Z}  \left(\left(3^{-\frac{m}{2}} \ell^{-\frac{1}{2}}(\ec) + \mu \right) \norm{\nabla \vb}_{L^2(\cu_m)} + \norm{\Delta \vb}_{L^2(\itr(\cu_m))}\right.\\
& \qquad \qquad + \left. \norm{\nabla \vb}^{\frac{1}{2}}_{L^2(\cu_m)} \norm{\Delta \vb}^{\frac{1}{2}}_{L^2(\itr(\cu_m))}  \right).
\end{split}
\end{equation}
\end{theorem}

Another topic studied in detail in this paper is an object called \textit{centered flux} defined for each $p \in \Rd$ by
\begin{equation}\label{eq:centeredFlux}
 \g_{p} := \ac(\D\phi_{p} + p) - \ab p,    
\end{equation}
 where $\ac(x,y) := \a(x,y)\Ind{x,y \in \clt_{\infty}}$. Estimating the weak convergence of this quantity to zero is one of the fundamental ingredients required to prove the two-scale convergence \Cref{thm:TwoScale}. Its physical interpretation is clear: we define $l_{p}(x):= p \cdot x$ and recall that the harmonic function can be seen as an electric potential. Then $(l_{p} + \phi_{p})$ is the electric potential defined on $\clt_{\infty}$ with conductance~$\a$ associated to the direction $p$,  while $l_p$ is the one for the homogenized conductance~$\ab$. We know that $\phi_p$ as the difference between the electric potentials is small compared to $l_p$, and heuristically, it should also be the case for the electric current. By Ohm's law, the two electric currents are defined by $\ac \nabla (\phi_{p} + l_p)$ and $\ab \nabla l_p$, so we expect indeed that $\g_p$ will be small. This will however only be true in a weak sense, or equivalently, after a spatial convolution. In fact, we expect that $\g_p$ satisfies estimates that are very similar to those satisfied by $\nabla \phi_{p}$, and we will indeed prove an analogue of the result of \cite[Proposition 3.3]{dario2018optimal} which concerned the weak convergence of $\nabla \phi_p$ to $0$. Here we use the notation $[\cdot]$ to represent the constant extension on every cube of the form $z + (-\frac{1}{2}, \frac{1}{2})^{d}$ for some $z \in \Zd$.

\begin{proposition}[Spatial average]\label{prop:SpatialAverage2}
There exist two positive constants $s:=s(d,\pp, \Lambda), C:=C(d,\pp, \Lambda, s)$ such that for every $R\geq 1$ and every kernel $K_R : \Rd \rightarrow \mathbb{R}^+$ integrable and satisfying
\begin{equation}\label{eq:ConditionK}
\exists C_{K,R} < \infty, \quad \forall x \in \Rd, \quad K_R(x) \leq \frac{C_{K,R}}{R^d \left( \vert \frac{x}{R} \vert \vee 1 \right)^{\frac{d+1}{2}}},
\end{equation}
the quantity $\left( K_{R} \star \coa{\g_{p}}\right)(x)$ is well defined for every $x \in \Rd, p \in \Rd$ and it satisfies
\begin{equation}\label{eq:SpatialAverage2}
\left\vert K_{R} \star \coa{\g_{p}}\right\vert(x)  \leq \GO_s(C_{K,R} C |p| R^{-\frac{d}{2}}).
\end{equation}
\end{proposition}
Other results including $L^q$ or $L^{\infty}$ estimates on $\g_{p}$ can be deduced from \Cref{prop:SpatialAverage2} by the multi-scale Poincar\'e inequality, see \cite[Appendix A]{dario2018optimal}.

\subsection{Organization of the paper}
In \Cref{sec:Pre}, we define all the notations precisely and restate some important theorems in previous work. \Cref{sec:Flux} is devoted to the study of the centered flux $\g_{p}$ and to the proof of \Cref{prop:SpatialAverage2}. \Cref{sec:TwoScale} gives the proof of the two-scale expansion on the cluster of percolation (\Cref{thm:TwoScale}). In \Cref{sec:Analysis}, we use the two-scale expansion to analyze our algorithm. Finally, in \Cref{sec:Numeric}, we present numerical experiments confirming the usefulness of the algorithm.

\section{Preliminaires}\label{sec:Pre}
This part defines rigorously all the notations used throughout this article. We also record some important results developed in previous work. 

\subsection{Notations $\GO_s(1)$ and its operations}
We recall the definition of $\GO_s$
\begin{equation}\label{eq:ODef}
X \leq \GO_s(\theta) \iff \Expt\left[\exp((\theta^{-1}X)_{+}^s)\right] \leq 2,
\end{equation} 
where $(\theta^{-1}X)_+$ means $\max\{\theta^{-1}X, 0\}$. One can use the Markov inequality to obtain that 
$$
X \leq \GO_s(\theta) \Longrightarrow \forall x > 0, \Prb[X \geq \theta x] \leq 2 \exp(-x^s).
$$
For $\lambda \in \mathbb{R}^+, X \leq \GO_s(\theta) \Longrightarrow \lambda X \leq \GO_s(\lambda \theta)$. We list some results on the estimates of the random variables with respect of $\GO_s$ in \cite[Appendix A]{armstrong2018quantitative}. For the product of random variables, we have
\begin{equation}
\label{eq:OMultica}
|X| \leq \GO_{s_1}(\theta_1), |Y| \leq \GO_{s_2}(\theta_2) \Longrightarrow |XY| \leq \GO_{\frac{s_1 s_2}{s_1 + s_2}}(\theta_1 \theta_2).
\end{equation}
By choosing $Y=1$, one can always use the estimate above to get an estimate for smaller exponent i.e. for $0 < s' < s,$ there exists a constant $C_{s'} < \infty$ such that 
\begin{equation}
\label{eq:OSmaller}
X \leq \GO_s(\theta) \Longrightarrow X \leq \GO_{s'}(C_{s'}\theta).
\end{equation}
We have an estimate on the sum of a series of random variables: for a measure space $(E, \mathcal{S}, m)$ and $\{X(z)\}_{z \in E}$ a family of random variables, we have 
\begin{equation}
\label{eq:OSum}
\forall z \in E, X(z )\leq \GO_s(\theta(z)) \Longrightarrow \int_{E} X(z) m(dz) \leq \GO_s\left( C_s \int_{E} \theta(z) m(dz) \right),
\end{equation}
where $0 < C_s < \infty$ is a constant defined by 
\begin{equation}
C_s = \left\{
	\begin{array}{ll}
	\left(\frac{1}{s \log 2} \right)^{\frac{1}{s}}	 & s < 1,\\
	1 & s \geq 1.
	\end{array}
\right.
\end{equation}
Finally, we can also obtain the estimate of the maximum of a finite number of random variables, which is proved in \cite[Lemma 3.2]{gu2018uniform}: for all $N \geq 1$ and  family of random variables $\{X_i\}_{1\leq i \leq N}$ satisfying that $X_i \leq \GO_s(1)$, we have 
\begin{equation}
\label{eq:OMax}
\left( \max_{1 \leq i \leq N} X_i \right) \leq \GO_s\left( \left(\frac{\log(2N)}{\log(3/2)}\right)^{\frac{1}{s}}\right).
\end{equation}

\subsection{Discrete analysis}
This part is devoted to introducing notations and some functional inequalities on graphs or on lattices. We take two systems of derivative in our setting: $\nabla$ on graph and the finite difference $\D$ on $\Zd$. The notation $\nabla$ is more general, but it loses the sense of derivative with respect to a given direction, which is very natural in the system of $\D$. 
\subsubsection{Spaces and functions}
For every $ V \subset \Zd$, we can construct two types of geometry $(V, \Ed(V))$ and $(V, \Eda(V))$. The set of edges $\Ed(V)$ inherited from $(\Zd,\Ed)$ and $\Eda(V)$ inherited from the open bonds of the percolation are as 
$$
\Ed(V) := \left\{\{x,y\} | x,y \in V, x \sim y \right\}, \qquad \Eda(V) : =  \left\{\{x,y\} | x,y \in V, \a(x,y) \neq 0 \right\}.
$$
The \textit{interior} of $V$ with respect to $(V, \Ed(V))$ and $(V, \Eda(V))$ are defined
$$
\itr(V) := \{x\in V | y \sim x \Longrightarrow y \in V\}, \qquad \itr_{\a}(V) := \{x\in V | y \sim x, \a(x,y) \neq 0 \Longrightarrow y \in V\},
$$
and the \textit{boundaries} are defined as $\partial(V) := V \backslash \itr(V)$ and $\partial_{\a}(V) := V \backslash \itr_{\a}(V)$.
For any $x,y \in \Zd$, we say $x \xleftrightarrow{\a} y$ if there exists an open path connecting $x$ and $y$.

We denote by $\overrightarrow{\Ed}$ \textit{the oriented bonds} of $(\Zd, \Ed)$ i.e. $
\overrightarrow{\Ed} := \{(x,y) \in \Zd \times\ \Zd : |x-y| = 1\}$, and for any $E \subset \Ed$, we can associate it to a natural oriented bonds set $\overrightarrow{E}$. An \textit{(anti-symmetric) vector field} $\overrightarrow{F}$ on $\overrightarrow{\Ed}$ is a function $\overrightarrow{F} : \overrightarrow{\Ed} \rightarrow \mathbb{R}$ such that $\overrightarrow{F}(x,y) = - \overrightarrow{F}(y,x)$. Sometimes we also write $\overrightarrow{F}(e)$ for $e = \{x,y\} \in \Ed$ to give its value with an arbitrary orientation for $e$, in the case it is well defined (for example $\vert \overrightarrow{F} \vert(e)$). The \textit{discrete divergence} of $\overrightarrow{F}$ is defined as $\nabla \cdot \overrightarrow{F} : \Zd \rightarrow \mathbb{R}$
$$
\forall x \in \Zd, \quad \nabla \cdot \overrightarrow{F}(x) := \sum_{y \sim x} \overrightarrow{F}(x,y).
$$ 
For any $u : \Zd \rightarrow \mathbb{R}$, we define the discrete derivative $\nabla u : \overrightarrow{\Ed} \rightarrow \mathbb{R}$ as a vector field 
$$
\forall (x,y) \in \overrightarrow{\Ed}, \quad \nabla u(x,y):= u(y) - u(x),
$$ 
and $\a \nabla u, \nabla u \Ind{\a \neq 0}$ are vector fields $\overrightarrow{\Ed} \rightarrow \mathbb{R}$ defined by
\begin{align*}
\a \nabla u (x,y) := \a(x,y)\nabla u(x,y), \qquad \nabla u \Ind{\a \neq 0}(x,y) := \nabla u(x,y)\Ind{\a(x,y) \neq 0}.
\end{align*}
Then, the $\a$-Laplacian operator $-\nabla \cdot \a \nabla$ is well defined and we have 
$$
-\nabla \cdot \a \nabla u(x) := \sum_{y \sim x}\a(x,y)(u(x) - u(y)).
$$

\subsubsection{Finite difference derivative} 
We start by introducing the notation of translation: let $B$ be a Banach space, then for any $h \in \Zd$ and $u : \Zd \rightarrow B$ a $B$-valued function, we define $T_h$ as an operator
$$
\forall x \in \Zd, \quad (T_h u)(x) = u(x + h).
$$
We also define the operator $\D_h$ and its conjugate operator $\D^*_{h}$ for any $u : \Zd \rightarrow \mathbb{R}$,
$$
\D_{h}u := T_h u - u, \qquad \D^*_{h} u := T_{-h}u - u.
$$
It is easy to check $\D^*_{h} = - T_{-h}(\D_{h}u)$ and for two functions $f,g : \Zd \rightarrow \mathbb{R}$, we have
\begin{equation}\label{eq:DiffProduct}
\D_h (fg) = (\D_hf)g + (T_h f)(\D_h g).
\end{equation}
In this system, we also define \textit{vector field} $\widetilde{F} : \Zd \rightarrow \Rd, \widetilde{F}(x) = (\widetilde{F}_1(x), \widetilde{F}_2(x) \cdots \widetilde{F}_d(x))$ and this can be distinguished with the one defined on $\Ed$ by the context. We use $(\e_1, \e_2 \cdots \e_d)$ to represent the $d$ canonical directions in $\Zd$, and a discrete gradient $\D u : \Zd \rightarrow \mathbb{R}^d$ is a vector field
$$
\D u(x) := \left(\D_{\e_1}u(x), \D_{\e_2}u(x) \cdots \D_{\e_d}u(x)\right).
$$ 
Then the finite difference divergence operator is defined as the conjugate operator of $\D$
$$
\D^* \cdot \widetilde{F} := \sum_{j=1}^d \D^*_{\e_j}\widetilde{F}_j.
$$
As convention, we use the notation $\a \Dr{j}u$ and $\a \D u$ to represent 
\begin{align*}
\a \Dr{j}u(x) := \a(x, x +\e_j) \Dr{j}u(x), 
\qquad \a \D u := \left(\a \Dr{1}u, \cdots \a \Dr{d}u  \right),
\end{align*}
and $\Ind{\a \neq 0}  \Dr{j}u, \Ind{\a \neq 0}  \D u,$ to represent 
\begin{align*}
\Ind{\a \neq 0}  \Dr{j}u(x):= \Ind{\a(x, x+\e_j) \neq 0}  \Dr{j}u(x), 
\qquad \Ind{\a \neq 0} \D u := \left(\Ind{\a \neq 0}  \Dr{1}u, \cdots \Ind{\a \neq 0}  \Dr{d}u  \right).
\end{align*}
Thus $\a$-Laplacian operator $-\nabla \cdot \a \nabla$ can also be defined by finite difference $ \D^* \cdot \a\D$. We can prove by a simple calculation that
\begin{equation}\label{eq:DivergenceForm}
 - \nabla \cdot \a \nabla u = \D^* \cdot \a\D u.
\end{equation}

\subsubsection{Inner product and norm}
For $V \subset \Zd$ and $E \subset \Ed$, we define inner product $\bracket{\cdot,\cdot}_V$ for any function $u,v : V \rightarrow \mathbb{R}$ and $\bracket{\cdot,\cdot}_E$ for any vector field $\overrightarrow{F},\overrightarrow{G} : \overrightarrow{E} \rightarrow \mathbb{R}$ 
$$
\bracket{u,v}_V := \sum_{x \in V} u(x)v(x), \qquad \bracket{\overrightarrow{F},\overrightarrow{G}}_E := \sum_{\{x,y\}\in E} \overrightarrow{F}(x,y)\overrightarrow{G}(x,y),
$$
and this defines a norm $\norm{u}_{L^2(V)} := \sqrt{\bracket{u,u}_V}$ and $\norm{\overrightarrow{F}}_{L^2(E)} := \sqrt{\bracket{\overrightarrow{F},\overrightarrow{F}}_E}$. We also abuse a little the notation to define $\bracket{\cdot, \cdot}_V$ for vector field    
$$
\bracket{\overrightarrow{F},\overrightarrow{G}}_V := \bracket{\overrightarrow{F},\overrightarrow{G}}_{\Ed(V)} = \sum_{\{x,y\}\in \Ed(V)}\overrightarrow{F}(x,y)\overrightarrow{G}(x,y) = \frac{1}{2}\sum_{x,y\in V, y \sim x}\overrightarrow{F}(x,y)\overrightarrow{G}(x,y).
$$
We use the notation $\bracket{\cdot,\cdot}_{\Eda(V)}$ to represent the inner product of the vector field on $(V, \Eda(V))$.   
For two vector fields $\widetilde{F}, \widetilde{G} : V \rightarrow \Rd$, the inner product is defined as 
$$
\bracket{\widetilde{F}, \widetilde{G}}_V := \sum_{x \in V}\sum_{j=1}^d \widetilde{F}_j(x) \widetilde{G}_j(x),
$$
and similarly $\norm{\widetilde{F}}_{L^2(V)} = \sqrt{\bracket{\widetilde{F},\widetilde{F}}_E}$ also defines a norm. 

To define a general $L^p(V) (p \geq 1)$ norm for vector fields, we have to introduce its modules.  For any $\overrightarrow{F} : \overrightarrow{\Ed} \rightarrow \mathbb{R}$ or $\widetilde{F} : \Zd \rightarrow \Rd$, we write 
$$
\vert \overrightarrow{F} \vert(x) := \left(\frac{1}{2}\sum_{y \sim x} \overrightarrow{F}^2(x,y)\right)^{\frac{1}{2}}, \qquad \vert \widetilde{F} \vert(x) := \left(\sum_{j=1}^d \widetilde{F}^2_j(x)\right)^{\frac{1}{2}}.
$$
Then for $f$ (a function, an $\Rd$-valued vector field or a vector field on $\overrightarrow{\Ed}$)  
$$
\norm{f}^p_{L^p(V)} := \left(\sum_{x \in V} |f|^p(x)\right)^{\frac{1}{p}}, \qquad \norm{f}^p_{\aL^p(V)} := \left(\frac{1}{|V|}\sum_{x \in V} |f|^p(x)\right)^{\frac{1}{p}}.
$$

We recall $C_0(V)$ the space of functions supported on $V$ with null boundary condition. Then one can deduce integration by part formula: for any function $v \in C_0(V)$, $\overrightarrow{F} : \overrightarrow{\Ed}(V) \rightarrow \mathbb{R}$ and $\widetilde{F} : \Zd \rightarrow \Rd$, one can check
\begin{equation}\label{eq:IPP}
\bracket{v,-\nabla \cdot \overrightarrow{F}}_{\itr(V)} = \bracket{\nabla v, \overrightarrow{F}}_V, \qquad \bracket{v, \D^* \cdot \widetilde{F}}_{\itr(V)} = \bracket{\D v, \widetilde{F}}_V.
\end{equation}

\subsubsection{Some functional inequalities}

We list some discrete functional inequalities used throughout the article.
\begin{lemma}[Discrete functional inequality]
\begin{enumerate}[(i)]
\item (\textit{A naive estimate}) Given a $V \subset \Zd$ and for a function $v : V \rightarrow \mathbb{R}$, we have 
\begin{equation}\label{eq:NaiveInequality}
\bracket{\nabla v, \nabla v}_{V} \leq 2d \bracket{v,v}_V.
\end{equation}

\item (\textit{Poincar\'e's inequality}) For every $ v \in C_0(\cu_m)$, we have
\begin{equation}\label{eq:Poincare}
\norm{v}_{L^2(\cu_m)} \leq 3^m \norm{\nabla v}_{L^2(\cu_m)}.
\end{equation}

\item (\textit{$H^2$ interior regularity for discrete harmonic function}) Given two functions $v,f \in C_0(\cu_m)$ satisfying the discrete elliptic equation ($\Delta v = \nabla \cdot \nabla v$)
\begin{equation}\label{eq:Harmonic}
- \Delta v = f, \qquad \text{ in } \itr(\cu_m),
\end{equation}
then we have an interior estimate
\begin{equation}\label{eq:H2DisInter}
\norm{\D^* \D v}^2_{L^2(\itr(\cu_m))} := \sum_{i,j = 1}^d \norm{\Dr{i}^* \Dr{j} v}^2_{L^2(\itr(\cu_m))} \leq d \norm{f}^2_{L^2(\itr(\cu_m))}.
\end{equation}

\item (\textit{Trace inequality}) For every $u : \cu_m \rightarrow \mathbb{R}$ and $0 \leq K \leq \frac{3^m}{4}$, we have the following inequality 
\begin{equation}\label{eq:Trace}
\norm{u \Ind{\dist(\cdot, \partial \cu_m) \leq K}}^2_{L^2(\cu_m)} 
\leq C(d) (K+1) \left(3^{-m}\norm{u}^2_{L^2(\cu_m)} + \norm{u}_{L^2(\cu_m)} \norm{\nabla u}_{L^2(\cu_m)}\right).
\end{equation}
\end{enumerate}
\end{lemma}
The inequality \eqref{eq:NaiveInequality} is very elementary, and the proof of \cref{eq:Poincare} is similar to the standard case, so we skip their proofs.  The inequality \eqref{eq:H2DisInter} is also relatively standard, but involves a careful calculation. The argument for \cref{eq:Trace} is more combinatorial and non-trivial. We provide their proofs in \Cref{sec:Appendix}.

\subsection{Partition of good cubes}
One difficulty to treat the function defined on the percolation clusters comes from its random geometry. To overcome this problem, \cite{armstrong2018elliptic} introduces a Calder\'on-Zygmund type partition of good cubes, and we recall it here. 

We denote by $\Tcu$ the \textit{triadic cube} and $\cu_m(z)$ is defined by 
$$
\cu_m(z) := \Zd \bigcap \left(z + \left(-\frac{1}{2} 3^m, \frac{1}{2}3^m \right) \right), z \in 3^m \Zd, m \in \mathbb{N},
$$
where \textit{center} and \textit{size} of the cube above is respectively $z$ and $3^m$, and we use the notation $\size(\cdot)$ to refer to the size, i.e. $\size(\cu_m(z)) = 3^m$. In this paper, without further mention, we use the word ``cube" for short of ``triadic cube" and $\cu_m$ for short of $\cu_m(0)$. The collection of all the cubes of size $3^n$ is defined by $\Tcu_n$, i.e. $\Tcu_n := \{z + \cu_n : z \in 3^n \Zd\}$. Then we have naturally 
$ \Tcu = \bigcup_{n \in \mathbb{N}} \Tcu_n.$ Every cube of size $3^m$ can be divided into a partition of $3^{(m-n)}$ cubes in $\Tcu_n$, and two cubes in $\Tcu$ can be either disjoint or included one by the other. For each $\cu \in \Tcu$, the \textit{predecessor} of $\cu$ is the unique triadic cube $\widetilde{\cu} \in \Tcu$ satisfying 
$$
\cu \subset \widetilde{\cu}, \text{ and } \frac{\size(\widetilde{\cu})}{\size(\cu)} = 3,
$$
and reciprocally, we say $\cu$ is a \textit{successor} of $\widetilde{\cu}$.

The distance between two points $x,y \in \Rd$ is defined to be $\dist(x,y) = \max_{i\in \{1,2 \cdots d\}} |x_i - y_i|$ and the distance for $U,V \subset \Zd$ is $\dist(U,V) = \inf_{x \in U, y \in V}\dist(x,y)$. In particular, two $\cu, \cu'$ are neighbors if and only if $\dist(\cu, \cu') = 1$ and one is included in the other if and only if $\dist(\cu, \cu') = 0$. 

\subsubsection{General setting}
We state at first the general setting of partition of good cubes.

\begin{proposition}[Proposition 2.1 of \cite{armstrong2018elliptic}]\label{prop:GeneralPartion} 
Let $\Gcu \subset \Tcu$ a sub-collection of triadic cubes satisfying the following: for every $\cu = z + \cu_n \in \Tcu$,
$$
\{\cu \notin \Gcu \} \in \F(z + \cu_{n+1}),
$$ 
and there exist two positive constants $K,s$
$$
\sup_{z \in 3^n \Zd}\Prb[z + \cu_n \notin \Gcu] \leq K \exp(-K^{-1} 3^{ns}).
$$
Then, $\Prb$-almost surely there exists  $\Scu \subset \Tcu$ a partition of $\Zd$ with the following properties: 
\begin{enumerate}
\item{Cubes containing elements of $\Scu$ are good}: for every $\cu, \cu' \in \Tcu$,
$$
\cu \subset \cu', \cu \in \Scu \Longrightarrow \cu' \in \Gcu.
$$ 
\item{Neighbors of elements of $\Scu$ are comparable}: for every $\cu, \cu' \in \Scu$ such that $\dist(\cu, \cu') \leq 1$, we have 
$$
\frac{1}{3} \leq \frac{\size(\cu)}{\size(\cu')} \leq 3.
$$
\item{Estimate for the coarseness}: we use $\cu_{S}(x)$ to represent the unique element in $\Scu$ containing a point $x \in \Zd$, then there exists a positive constant $C:=C(s,K,d)$ such that, for every $x \in \Zd$,
$$
\size(\cu_{\Scu}(x)) \leq \GO_s(C).
$$
\end{enumerate} 
\end{proposition}

\subsubsection{Case of well-connected cubes}
The construction in \Cref{prop:GeneralPartion} works for all collection of good cubes $\Gcu$, here we give the concrete definition of good cubes we use in our context of percolation, as  appearing in the work \cite{penrose1996large}, \cite{pisztora1996surface} and \cite{antal1996chemical} of Antal and Pisztora. We remark that in \Cref{def:crossclt} and \Cref{def:goodcube} we use ``cube" exceptionally for a general lattice cube, and we will highlight explicitly ``triadic cube" when using it. The notation $\frac{3}{4} \overline{\cu}$ indicates that we take the convex hull of the lattice cube, and then change its size by multiplying by $\frac{3}{4}$ while keeping the center fixed.

\begin{definition}[Crossability and crossing cluster]\label{def:crossclt}
We say that a cube $\cu$ is \textit{crossable} with respect to the open edges defined by $\a$ if each of the $d$ pairs of opposite $(d-1)$-dimensional faces of $\cu$ can be joined by an open path in $\cu$. We say that a cluster $\clt \subset \cu$ is a \textit{crossing cluster} for $\cu$ if $\clt$ intersects each of the $(d-1)$-dimensional faces of $\cu$. 
\end{definition}

\begin{definition}[Well-connected cube and good cube, Theorem 3.2 of \cite{pisztora1996surface}]\label{def:goodcube}
We say that $\cu \in \Tcu$ is \textit{well-connected} if there exists a crossing cluster $\clt$ for $\cu$ such that :
\begin{enumerate}
\item \label{def:goodcube1} each cube $\cu'$ with $\frac{1}{10}\size(\cu) \leq \size(\cu') \leq \frac{1}{2}\size(\cu)$ and $\cu' \cap \frac{3}{4} \overline{\cu} \neq \emptyset$ is crossable.
\item \label{def:goodcube2} every path $\gamma \subset \cu'$ defined above with $\diam(\gamma) \geq \frac{1}{10} \size(\cu)$ is connected to $\clt$ within $\cu'$.
\end{enumerate}
We say that $\cu \in \Tcu$ is a \textit{good cube} if $\size(\cu) \geq 3$, $\cu$ is connected and all his $3^d$ successors are well-connected. Otherwise, we say that $\cu \in \Tcu$ is a \textit{bad cube}.
\end{definition}

The following estimates makes the the construction defined in \Cref{prop:GeneralPartion} work. 
\begin{lemma}[(2.24) of \cite{antal1996chemical}]\label{lem:AntalSize} For each $\pp \in (\pc, 1]$, there exists a positive constant $C:=C(d,\pp)$ such that for every $n \in \mathbb{N}$,
$$
\Prb[\cu_n \in \Gcu] \geq 1 - C \exp(-C^{-1}3^n).
$$
\end{lemma}

\begin{definition}[Partition of good cubes in percolation context]
We let $\Pcu \subset \Tcu$ be the partition $\Scu$ of $\Zd$ obtained by applying \Cref{prop:GeneralPartion} to the collection of good cubes defined in \Cref{def:goodcube}
$$
\Gcu := \left\{\cu \in \Tcu : \cu \text{ is good cube } \right\}.
$$
\end{definition}

A direct application of \Cref{lem:AntalSize} and \Cref{prop:GeneralPartion} gives us:
\begin{corollary}
There exists a positive constant $C(d, \pp)$, such that for every $z \in \Zd$, we have the two estimates  
\begin{equation}\label{eq:sizeP}
\size(\cu_{\Pcu}(z)) \leq \GO_1(C), \qquad \Ind{\size(\cu_{\Pcu}(z)) \geq n} \leq \GO_1(C 3^{-n}).
\end{equation} 
\end{corollary}

The maximal cluster is well defined on every good cube by \Cref{def:goodcube}.
\begin{definition}[Maximal cluster in good cubes]
For every good cube $\cu$, there exists a unique maximal crossing cluster in it, and we denote this cluster by $\clt_*(\cu)$. 
\end{definition}

Although $\clt_*(\cu)$ only uses local information, the next lemma shows that, for a $\cu \in \Pcu$ (stronger than $\cu$ is good), its maximal cluster $\clt_*(\cu)$ must belong to the infinite cluster $\clt_{\infty}$.
\begin{lemma}[Lemma 2.8 of \cite{armstrong2018elliptic}]\label{lem:FusionCluster}
Let $n,n' \in \mathbb{N}$ with $|n-n'|\leq 1$ and $z,z' \in 3^n \Zd$ such that 
$$
\dist(\cu_n(z), \cu_{n'}(z')) \leq 1.
$$
Suppose also that $\cu_n(z)$ and $\cu_{n'}(z')$ are all good cubes, then there exists a cluster $\clt$ such that 
$$
\clt_*(\cu_n(z)) \cup \clt_*(\cu_{n'}(z')) \subset \clt \subset \cu_n(z) \cup \cu_{n'}(z').
$$
\end{lemma}

This lemma helps us generalize the definition of maximal cluster in a general set $U \subset \Zd$, the idea is to define the union of the partition cubes that cover $U$, and then find the maximal cluster in it.

\begin{definition}[Maximal cluster in general set]\label{def:MaximalCluster}
For a general set $U \subset \Zd$, we define its \textit{closure with respect to} $\Pcu$ by 
\begin{equation}\label{eq:Closure}
\clp(U) := \bigcup_{z \in U}\cu_{\Pcu}(z),
\end{equation}
and $\clt_*(U)$ to be the cluster contained in $\clp(U)$ which contains all the clusters of $\clt_*(\cu_{\Pcu}(z)), z \in U$.
\end{definition}

One can check easily that \Cref{lem:FusionCluster} makes the definition $\clt_*(U)$ well-defined. However, we do not have necessarily $\clt_*(U) = \bigcup_{z \in U} \clt_*(\cu_{\Pcu}(z))$. We provide with a detailed discussion of this point in \Cref{sec:SmallCluster}.

Since the cubes in $\Tcu$ can be either included in one another or disjoint, if one cube $\cu \in \Tcu$ contains an element in $\Pcu$, then it can be decomposed as the disjoint union of elements in $\Pcu$ without enlarging the domain. Thus, we define:
\begin{definition}[Minimal scale for partition]\label{def:MinimalScale}
\begin{equation}\label{eq:defPstar}
\Pcu_* = \left\{ \cu \in \Tcu : \exists \cu' \subset \cu \text{ and } \cu' \in \Pcu \right\}.
\end{equation}
\end{definition}
The following observations are very useful and can be checked easily: for every $\cu \in \Tcu$ and $z$ as its center, we have 
\begin{align}\label{eq:PcuScale}
\cu \in \Pcu_* \Longrightarrow \clp(\cu) = \cu, \qquad
\Ind{\cu \notin \Pcu_*} \leq \Ind{\size(\cu_{\Pcu}(z)) > \size(\cu)} \leq \GO_1(C (\size(\cu))^{-1}). 
\end{align}

\subsubsection{Mask operation and coarsened function}
To overcome the problem of the passage between the two geometries $(\Zd, \Ed)$ and $(\clt_{\infty}, \Eda)$, one useful technique is the mask operation. 
\begin{definition}[Mask operation and local mask operation]\label{def:MaskOperation}
For $f : \Zd \rightarrow \mathbb{R}$ and $\a : \Ed \rightarrow \mathbb{R}$, we define a \textit{mask operation} $(\cdot)_{\clt}$ to make their support on $\clt_{\infty}$ and $\Eda(\clt_{\infty})$ respectively
\begin{equation}\label{eq:MaskOperation}
f_{\clt}(x) := \left\{
	\begin{array}{ll}
	f(x) & x \in \clt_{\infty}, \\	
	0 & \text{otherwise}.
	\end{array}
\right. \qquad
\ac(x,y) := \left\{
	\begin{array}{ll}
	\a(x,y) & x,y \in \clt_{\infty}, \\	
	0 & \text{otherwise}.
	\end{array}
\right.
\end{equation}
Moreover, we also define a \textit{local mask operation} for $\cu_m \in \Gcu$ as  
\begin{equation}\label{eq:LocalMaskOperation}
f_{\clt,m}(x) := \left\{
	\begin{array}{ll}
	f(x) & x \in \clt_{*}(\cu_m), \\	
	0 & \text{otherwise}.
	\end{array}
\right. \qquad
\acf{m}(x,y) := \left\{
	\begin{array}{ll}
	\a(x,y) & x,y \in \clt_{*}(\cu_m), \\	
	0 & \text{otherwise}.
	\end{array}
\right.
\end{equation}
Then we call $f_{\clt} ( f_{\clt, m} ), \ac ( \acf{m} )$ \textit{(local) masked function and (local) masked conductance}.
\end{definition}
Reciprocally, for a function only defined on the clusters, sometimes we have to extend them to the whole space. We can apply the technique of coarsening the function defined on the percolation cluster. 
\begin{definition}[Coarsened function]\label{def:CoarsenedFunction}
Given $\cu \in \Pcu$, we let $\zb(\cu)$ represent the vertex in $\clt_{*}(\cu)$ which is closest to its center. Given a function $u : \clt_{\infty} \rightarrow \mathbb{R}$, we define the coarsened function with respect to $\Pcu$ to be $\coa{u}_{\Pcu} : \Zd \rightarrow \mathbb{R}$ that 
$$
\coa{u}_{\Pcu}(x) := u(\zb(\cu_{\Pcu}(x))).
$$
We also use the notation $[\cdot]$ to mean doing constant extension on every cube i.e. given $ v : \Zd \rightarrow \mathbb{R}$, we define $[v] :\Rd \rightarrow \mathbb{R}$ such that for every $z \in \Zd$ and every $x \in z + \left[- \frac{1}{2},  \frac{1}{2}\right)^d$, $
[v](x) := v(z)$.
\end{definition}

The advantage of the coarsened function is that it allows to extend the support of function from $\clt_{\infty}$ to the whole space, and constant in every cube by paying a small cost of errors.
\begin{proposition}[Lemmas 3.2 and 3.3 of \cite{armstrong2018elliptic}]\label{prop:ErrorCoarse}
For every $1 \leq s < \infty$, there exists a positive constant $C(s, d, \pp)$, such that for every $\cu \in \Pcu_{*}$, $u : \clt_{\infty} \rightarrow \mathbb{R}$, we have 
\begin{align}
\sum_{x \in \clt_*(\cu)} \vert u(x) - \coa{u}_{\Pcu}(x) \vert^s \leq C^s \sum_{\{y,z\} \in \Eda(\clt_*(\cu))} \size(\cu_{\Pcu}(y))^{sd} \vert \nabla u\vert^s(y,z),  \label{eq:ErrorCoarse1}
\end{align}
\begin{multline}
\sum_{\{x,y\} \in \Ed(\clp(\cu))} \vert \nabla \coa{u}_{\Pcu}(x,y) \vert^s \\
 \leq C^s \sum_{\{x,y\} \in \Eda(\clt_*(\cu))} \left(\size(\cu_{\Pcu}(x))^{sd-1} + \size(\cu_{\Pcu}(y))^{sd-1} \right) \vert \nabla u\vert^s(x,y). \label{eq:ErrorCoarse2}
\end{multline}
\end{proposition}
\begin{remark}\label{rmk:CoarsenedFunctionLocal} The main idea of coarsened function is to give function a constant value in every cube, but the value does not have to be of the one closest to the center. Following the same idea of proof of \cite[Lemmas 3.2 and 3.3]{armstrong2018elliptic}, one can prove that for $\cu \in \Pcu_*, u \in C_0(\cu)$
\begin{equation}\label{eq:CoarsenedFunctionLocal}
\coa{u}_{\Pcu,\cu}(x) = \left\{
	\begin{array}{ll}
	\coa{u}_{\Pcu}(x) & \dist(\cu_{\Pcu}(x), \partial \clp(\cu)) \geq 1,\\
	0 	& \dist(\cu_{\Pcu}(x), \partial \clp(\cu)) = 0,
	\end{array}
\right.
\end{equation}
we have the same inequality as \cref{eq:ErrorCoarse1} and \cref{eq:ErrorCoarse2} by putting $\coa{u}_{\Pcu,\cu}$ in the place of $\coa{u}_{\Pcu}$.
\end{remark}

\subsection{Harmonic function on the infinite cluster}
We define $\A(U)$, the set of $\a$-harmonic functions on $U \subset \Zd$, by 
$$
\A(U) := \{v:\clt_{\infty} \rightarrow \mathbb{R} | -\nabla \cdot \ac \nabla v = 0, \forall x \in \itr_{\a}(U) \},
$$
and $\A(\clt_{\infty})$ the set $\a$-harmonic functions on $\clt_{\infty}$. The $\a$-harmonic function $\A_k(\clt_{\infty})$ is the subspace of $\a$-harmonic functions which grows more slowly than a polynomial of degree $k+1$:
$$
\A_k(\clt_{\infty}) := \left\{ u \in \A(\clt_{\infty}) \left| \limsup_{R \rightarrow \infty} R^{-(k+1)}\norm{u}_{\aL^2(\clt_{\infty} \cap B_R)} = 0 \right. \right\}.
$$
Similarly, we can define the spaces $\bar{\A}, \bar{\A}_k$ for harmonic functions on $\Rd$. It is well-known that the space $\bar{\A}_k$ is a finite-dimensional vector space of polynomials. A recent remarkable result about $\a$-harmonic functions on the infinite cluster of percolation conjectured in \cite{benjamini2015disorder} and proved in \cite{armstrong2018elliptic} is that the space $\A_k(\clt_{\infty})$ also has this property, and in fact has the same dimension as $\bar{\A}_k$. Here we only recall the structure of $\A_1(\clt_{\infty})$: for every $\a$-harmonic functions $u \in \A_{1}(\clt_{\infty})$, there exists $c \in \mathbb{R}, p \in \mathbb{R}$ such that 
$$
\forall x \in \clt_{\infty}, \quad u(x) = c + p \cdot x + \phi_p (x),
$$
where the functions $\{\phi_{p}\}_{p \in \Rd}$ are called \textit{the first order correctors}.  The first order correctors have sublinear growth: there exists a positive exponent $\delta(d, \pp, \Lambda) < 1$ and a minimal scale $\mathcal{M} \leq \GO_s(C(d,\pp, \Lambda))$ such that, for every $r \geq \mathcal{M}$ and $p \in \Rd$,
\begin{equation}\label{eq:Correctors}
\norm{\phi_p}_{\aL^2(\clt_{\infty} \cap B_r)} \leq C |p| r^{1-\delta}.
\end{equation}
This property plays an important role in many proofs, and \cite{dario2018optimal} gives a more precise description of these correctors. We recall that $\Phi_R(x) := \frac{1}{(4\pi R^2)^{d/2}}\exp\left(-\frac{x^2}{4R^2}\right)$, and $\coa{\phi_p}_{\Pcu}^{\eta}:=\coa{\phi_p}_{\Pcu} \star \eta$ where $\eta \in C_0^{\infty}(B_1)$ is positive, and $\eta \equiv 1$ in $B_{\frac{1}{2}}$.
\begin{proposition}[Local estimate and spatial average estimate, Proposition 3.3 of \cite{dario2018optimal}]
There exist two positive constants $s:=s(d,\pp, \Lambda)$, $C:=C(d,\pp, \Lambda)$ such that for each $R\geq 1$ and each $p \in \Rd$,
\begin{align}
\forall x \in \Zd, \quad \vert \nabla \phi_{p} \Ind{\a \neq 0} \vert (x) &\leq \GO_s(C|p|), \label{eq:LocalAverage1} \\
\forall x \in \Rd, \quad \left\vert \nabla \left( \Phi_{R} \star \coa{\phi_p}_{\Pcu}^{\eta}\right)(x) \right\vert &\leq \GO_s(C\vert p\vert R^{-\frac{d}{2}}). \label{eq:SpatialAverage1}
\end{align}
\end{proposition}

\begin{proposition}[Theorem 1 and 2 of \cite{dario2018optimal}, $L^q$ estimates on $\clt_{\infty}$]\label{prop:LpCluster}
There exist three positive constants $s:=s(d,\pp, \Lambda), k:= k(d,\pp,\Lambda)$ and $C:=C(d,\pp, \Lambda)$ such that for each $q \in [1, \infty),$ $R \geq 1$ and $p \in \Rd$,
\begin{equation}
\left(R^{-d} \int_{\clt_{\infty} \cap B_R} \vert \phi_p - (\phi_p)_{\clt_{\infty} \cap B_R} \vert^q \right)^{\frac{1}{q}} \leq \left\{
	\begin{array}{ll}
	\GO_s(C\vert p\vert q^k \log^{\frac{1}{2}}(R) )  & d=2, \\	
	\GO_s(C\vert p\vert q^k )  & d=3,
	\end{array}
\right.
\end{equation}
and for every $x,y \in \Zd$ and $p \in \Rd$,
\begin{equation}
\vert \phi_{p}(x) - \phi_p(y) \vert \Ind{x,y \in \clt_{\infty}} \leq \left\{
	\begin{array}{ll}
	\GO_s(C\vert p\vert \log^{\frac{1}{2}}\vert x - y\vert )  & d=2, \\	
	\GO_s(C\vert p\vert )  & d=3.
	\end{array}
\right.
\end{equation}
\end{proposition}

\section{Centered flux on the cluster}\label{sec:Flux}
In this part, we will study an object $\g_{p} : \Zd \rightarrow \Rd$ called \textit{centered flux} defined by
$$
\g_{p} := \ac(\D\phi_{p} + p) - \ab p,
$$
where $\ac$ is the masked conductance defined in \cref{eq:MaskOperation} and it is $\a$ restricted on the infinite cluster $\clt_{\infty}$. Because $\g_{p}$ satisfies
 $\D^* \cdot \g_{p} = 0$ on $ \Zd$, following the spirit of Helmholtz-Hodge decomposition, in the later part of this section we will also study another object $\S_p : \Zd \rightarrow \mathbb{R}^{d \times d} $ called \textit{flux corrector} such that $\g_p = \D^* \cdot \S_p$ on $\Zd$, in the sense $\g_{p,i} = \sum_{j=1}^d \Dr{j}^* \S_{p, ij}$.

The quantities $\g_p$ and $\S_p$ are fundamental to the quantitative analysis of the two-scale expansion, see for instance \cite{gloria2011optimal} and 
\cite[Chapter 6]{armstrong2018quantitative}. Roughly speaking,  
$\D \phi_{p}, \g_{p}$ and $\D \S_p$ should satisfy similar estimates. The goal of this section is to study various quantities like spatial averages and $L^p$ and $L^{\infty}$ estimates on $\g_{p}$ and $\D \S_p$, as a counterpart of the work \cite{dario2018optimal} concerning $\D \phi_p$.

We can prove at first a very simple result.
\begin{proposition}[Local average]\label{prop:LocalAverage2}
There exit two positive constants $s:=s(d,\pp, \Lambda)$ and $C(d,\pp, \Lambda)$ such that 
\begin{equation}\label{eq:LocalAverage2}
\forall x \in \Zd, \quad \vert \g_{p} \vert(x) \leq \GO_s(C |p|).
\end{equation}
\end{proposition}
\begin{proof}
We have, by \cref{eq:LocalAverage1}
\begin{align*}
\vert \g_{p} \vert = \vert \ac(\D\phi_{p} + p) - \ab p \vert  \leq  \vert \ac\D\phi_{p} \vert + \vert \ac p \vert  + \vert \ab p \vert  \leq \vert \nabla \phi_{p} \Ind{\a \neq 0} \vert   + 2|p| \leq \GO_s(C|p|).
\end{align*}
\end{proof}

\subsection{Spatial average of centered flux}
In this part, we focus on the spatial average quantity $ K_{R} \star \coa{\g_{p}}$ and prove \Cref{prop:SpatialAverage2}. The spirit of the proof can go back to the spectral gap method (or Efron-Stein type inequality) in the work of Naddaf and Spencer \cite{naddaf1998estimates}, which is also employed in the work of Gloria and Otto \cite{gloria2011optimal,gloria2012optimal,gloria2014optimal}. \Cref{prop:SpatialAverage2} is more technical in two aspects: 
\begin{itemize}
    \item In the percolation context, the perturbation of the geometry of clusters has to be taken into consideration when applying the spectral gap method. 
    \item The result stated with $\GO_s$ notation requires a stronger concentration analysis.
\end{itemize}
Our proof follows generally the main idea of that of \cref{eq:SpatialAverage1} appearing in \cite[Proposition 3.3]{dario2018optimal}, and the main tool used in this proof is a variant of the Efron-Stein type inequality, combined with the Green function and Meyers' inequality on $\clt_{\infty}$. 
\begin{proof}[Proof of \Cref{prop:SpatialAverage2}]
Without loss of generality we suppose that $|p| = 1$, and the proof is decomposed into 4 steps.

\textit{Step 1: Spectral gap inequality and double environment.} We introduce the Efron-Stein type inequality used for the proof, which is proved first in \cite[Proposition 2.2]{armstrong2017optimal} and also used in \cite[Proposition 2.18]{dario2018optimal}. (We remark that there is a typo in the exponent in \cite[Proposition 2.2]{armstrong2017optimal}, which should be $\frac{2-\beta}{2}$; see also  \cite[Appendix A]{armstrong2017optimal} where the exponent is correct.) 

\begin{proposition}[Exponential Efron-Stein inequality, Proposition 2.2 of \cite{armstrong2017optimal}] Fix $\beta \in (0,2)$ and let $X$ be a random variable defined in the random space $(\Omega, \F, \Prb)$ generated by $\{\a(e)\}_{e \in \Ed}$, and we define 
\begin{align}
\F(\Ed \backslash \{e\}) &:= \sigma\left(\{\a(e')\}_{e' \in \Ed \backslash e} \right), \\
X_e := \Expt\left[X|\F(\Ed \backslash \{e\})\right], &\qquad \Vc[X] := \sum_{e \in \Ed} (X - X_e)^2.
\end{align}
Then, there exists a positive constant $C:=C(d, \beta)$ such that 
\begin{equation}\label{eq:EfronStein}
\Expt\left[\exp\left(|X - \Expt[X]|^{\beta}\right)\right] \leq C \Expt \left[\exp\left( (C \Vc[X])^{\frac{\beta}{2 - \beta}}\right) \right]^{\frac{2-\beta}{2}}.
\end{equation}
\end{proposition}
In the proof of \Cref{prop:LocalAverage2}, we apply this inequality by posing $X := \left( K_{R} \star \coa{\g_{p}}\right)(x)$ and we claim that it suffices to verify two conditions 
\begin{align}
\Expt[X] &\leq  C_1R^{-\frac{d}{2}}, \label{eq:EfronSteinCondition1}\\
\Vc[X] &\leq \GO_{s'}\left(C_2 R^{-d} \right) \label{eq:EfronSteinCondition2}.
\end{align}
It is also very natural, because the two conditions say that the average and fluctuation of $X$ are of the order of $R^{-\frac{d}{2}}$. We choose a $s$ such that $\frac{s}{2-s} = s'$ where $s'$ is the exponent in \cref{eq:EfronSteinCondition2} and $C_3 := (C_1 \vee C_2)C(d,\beta)$ where $C(d,\beta)$ is the constant in \cref{eq:EfronStein} and $C_1,C_2$ the one in \cref{eq:EfronSteinCondition1}, \cref{eq:EfronSteinCondition2}, then
\begin{align*}
\Expt\left[\exp \left(\left(\frac{X}{C_3R^{-\frac{d}{2}}} \right)^s \right)\right]  & \leq \Expt\left[\exp \left( \left(\frac{X - \Expt[X]}{C_3R^{-\frac{d}{2}}} + \frac{\Expt[X]}{C_3R^{-\frac{d}{2}}}  \right)^{s} \right) \right]\\
& \leq C \underbrace{\Expt\left[\exp \left( \left(\frac{ \vert X - \Expt[X] \vert}{C_3R^{-\frac{d}{2}}} \right)^{s} \right) \right]}_{\text{Using }\cref{eq:EfronStein}}\\
& \leq C \Expt\left[\exp\left( \left(\frac{\Vc[X]}{C_2R^{-d}} \right)^{\frac{s}{2-s}}\right)\right]^{\frac{2-s}{2}}\\
& \leq 2C. 
\end{align*}
Finally, we increase $C_3$ with respect to $s$ so that we get $X \leq \GO_s(CR^{-\frac{d}{2}})$. 

We focus on the two conditions \cref{eq:EfronSteinCondition1}, \cref{eq:EfronSteinCondition2}. In fact, we can check the condition \cref{eq:EfronSteinCondition1} by proving $\Expt[\ac(\D \phi_p + p)] = \ab p$, which is a well-known result in classic homogenization. In percolation context, it is also true by a careful check of the several equivalent definitions of $\ab$. We put its proof in \Cref{thm:CharaConductance}. 

To prove the condition \cref{eq:EfronSteinCondition2}, we use a useful technique in Efron-Stein type inequality of ``doubling" the probability space: we sample a copy of random conductance $\{\widetilde{\a}(e')\}_{e' \in \Ed}$ with the same law but independent to $\{\a(e')\}_{e' \in \Ed}$, and the two probability spaces generated by the two copies are denoted respectively by $\left(\Omega_{\a}, \F_{\a}, \Prb_{\a} \right), \left(\Omega_{\tilde{\a}}, \F_{\widetilde{\a}}, \Prb_{\widetilde{\a}} \right)$. Then we put the two copies of random conductance together and make a larger probability space $\left(\Omega', \F', \Prb' \right)=\left(\Omega_{\a} \times \Omega_{\widetilde{\a}}, \F_{\a} \otimes \F_{\widetilde{\a}}, \Prb_{\a} \otimes \Prb_{\widetilde{\a}} \right)$, and we also use the notation $\GO'_s$ to represent the same definition \cref{eq:ODef} in the larger probability space $\left(\Omega', \F', \Prb' \right)$. We also introduce the another random environment $\{\a^{e}(e')\}_{e' \in \Ed}$, obtained by replacing one conductance $\a(e)$ by $\widetilde{\a}(e)$, i.e.
\begin{equation}
\a^{e}(e') = \left\{ 
\begin{array}{ll}
\a(e') & e' \neq e,\\
\widetilde{\a}(e') & e' = e.
\end{array}
\right.
\end{equation}
We use $X^{e}, \clt_{\infty}^e, \phi_{p}^e$ to represent respectively the random variable, the infinite cluster and the corrector in the environment $\{\a^{e}(e')\}_{e' \in Ed}$. The definition of $\Vc[X]$ says that the variance comes from the fluctuation caused by the perturbation of every conductance, which suggests the following lemma:
\begin{lemma}\label{lem:DoubleTrickInequality}
We have the following estimate
\begin{equation}
\sum_{e \in \Ed}(X - X^e)^2 \leq \GO'_{s}(CR^{-d}) \Longrightarrow \Vc[X] \leq \GO_s(C_sR^{-d}).
\end{equation}
\end{lemma}
\begin{proof}
We use the double environment trick to see that 
$$
X_e = \Expt\left[X|\F(\Ed \backslash \{e\})\right] = \Expt_{\widetilde{\a}}[X^{e}],
$$
and Jensen's inequality to reformulate at first the inequality
\begin{align*}
\Expt\left[\exp\left( \left(\frac{\Vc[X]}{CR^{-d}} \right)^s \right)\right] &= \int_{\Omega} \exp\left( \left(\frac{\sum_{e \in \Ed}(X - X_e)^2}{CR^{-d}} \right)^s \right) \, d\Prb_{\a}(\omega)\\
&= \int_{\Omega} \exp \left( \left(\frac{\sum_{e \in \Ed}\left( \int_{\Omega_{\widetilde{\a}}} (X - X^e) \, d\Prb_{\widetilde{\a}}(\omega) \right)^2}{CR^{-d}} \right)^s \right) \, d\Prb_{\a}(\omega)\\
&\leq \int_{\Omega} \exp \left( \left(\int_{\Omega_{\widetilde{\a}}} \frac{ \sum_{e \in \Ed} (X - X^e)^2  }{CR^{-d}} \, d\Prb_{\widetilde{\a}}(\omega) \right)^s \right) \, d\Prb_{\a}(\omega)\\
\end{align*}
In the next step, we want to add a constant $t_s$ to make $\exp\left((\cdot + t_s)^s\right)$ convex, and then exchange the expectation and $\exp\left((\cdot +  t_s)^s\right)$ by Jensen's inequality. We can choose $t_s = 0$ for $s \geq 1$, and $t_s = \left(\frac{1-s}{s}\right)^{\frac{1}{s}}$ for $0 < s < 1$. (The spirit is the same as \cref{eq:OSum} and see \cite[Lemma A.4]{armstrong2018quantitative} for details of this proof.)
\begin{align*}
\Expt\left[\exp\left( \left(\frac{\Vc[X]}{CR^{-d}} \right)^s \right) \right] &\leq \int_{\Omega} \int_{\Omega_{\widetilde{\a}}} \exp \left( \left(\frac{ \sum_{e \in \Ed} (X - X^e)^2  }{CR^{-d}} + t_s \right)^s \right) \, d\Prb_{\widetilde{\a}}(\omega) \, d\Prb_{\a}(\omega)\\
&\leq \tilde{C}\int_{\Omega} \int_{\Omega_{\widetilde{\a}}} \exp\left( \left(\frac{ \sum_{e \in \Ed} (X - X^e)^2  }{CR^{-d}}\right)^s \right) \, d\Prb_{\widetilde{\a}}(\omega) \, d\Prb_{\a}(\omega)\\
& \leq 2\tilde{C}.
\end{align*}
In the last step we use the condition $\sum_{e \in \Ed}(X - X^e)^2 \leq \GO'_{s}(CR^{-d})$ and we reduce the constant $C$ to get the desired result. 
\end{proof}
By \Cref{lem:DoubleTrickInequality}, to prove \cref{eq:EfronSteinCondition2} it suffices to focus on the quantity $\sum_{e \in \Ed}(X - X^e)^2$ and in our context and to prove
\begin{equation}\label{eq:QuantDoubleEnvironment1}
\sum_{e \in \Ed}\left\vert K_{R} \star \left(\coa{\g_{p}} - \coa{\g^e_{p}} \right) \right\vert^2 (x) \leq \GO'_s(C R^{-d}).
\end{equation}
We will distinguish several cases and attack them one by one. 

\textit{Step 2: Case $\clt_{\infty} \neq \clt^e_{\infty}$, proof of $\g_{p} = \g_{p}^e$.} We have to consider the perturbation of the geometry between $\clt_{\infty}$ and $\clt^e_{\infty}$.

\begin{lemma}[Pivot edge]\label{lem:PivotEdge}
In the case $\clt_{\infty} \neq \clt^e_{\infty}$, we have necessarily the relation of inclusion, without loss of generality we suppose that $\clt^e_{\infty} \subset \clt_{\infty}$ , then we have:
\begin{enumerate}
\item \label{lem:PivotEdge1} The part $\clt_{\infty} \backslash \clt^e_{\infty}$ is connected to $\clt^e_{\infty} $ by $e$ (called \textit{pivot edge}), and $\vert \clt_{\infty} \backslash \clt^e_{\infty} \vert < \infty$.
\item \label{lem:PivotEdge2} We denote by $e:=\{e_*, e^*\}, e_* \in \clt^e_{\infty} \cap \clt_{\infty}$ and $e^* \in \clt_{\infty} \backslash \clt^e_{\infty}$, then the function $(\phi_{p}+l_{p})$ is constant on $\clt_{\infty} \backslash \clt^e_{\infty}$ and equals to $(\phi_{p} + l_{p})(e_*)$. 
\item \label{lem:PivotEdge3} The function $\phi^e_{p}$ has a representation that $\phi^e_{p} = \phi_{p}\Ind{\clt^e_{\infty}}$ up to a constant and satisfies $\ac\nabla (\phi_{p}+l_{p}) = \a^e_{\clt} \nabla (\phi^e_{p} + l_{p})$ on $\Ed$.
\end{enumerate} 
\end{lemma}
\begin{proof}
\begin{enumerate}
\item It comes from the fact that $\a$ and $\a^e$ are different only by one edge, thus $\clt^e_{\infty} \subsetneq \clt_{\infty}$ means that the edge $e$ is non-null in $\clt_{\infty}$ but zero in $\clt^e_{\infty} $ and makes one part disconnected from $\clt_{\infty}$. It is well-known that in the supercritical percolation, almost surely there exists one unique infinite cluster, thus we have $\vert \clt_{\infty} \backslash \clt^e_{\infty} \vert < \infty$.

\item We study the harmonic function $-\nabla \cdot \ac \nabla (\phi_{p} + l_{p}) = 0$ on the part $\clt_{\infty} \backslash \clt^e_{\infty}$. This is a non-degenerate linear system  with $\vert \clt_{\infty} \backslash \clt^e_{\infty} \vert$ equations and $\vert \clt_{\infty} \backslash \clt^e_{\infty} \vert+1$ variables, thus the solution is of $1$ dimension and we know this constant is $(\phi_{p}+l_{p})(e_*)$.

\begin{figure}[h]
\centering
\includegraphics[scale=0.5]{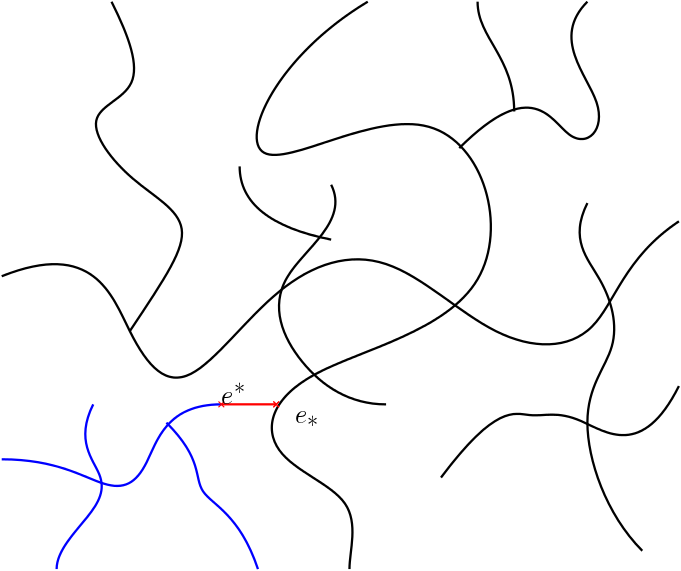}
\caption{In the image the segment in red is the edge $e = \{e_*, e^*\}$ and the part in blue is the cluster $\clt_{\infty} \backslash \clt^e_{\infty}$, where $\a$-harmonic function $(\phi_p + l_p)$ is constant of value $(\phi_p + l_p)(e_*)$.}
\end{figure}

\item We prove that at first that $\ac\nabla (\phi_{p}+l_{p}) = \a^e_{\clt} \nabla (\phi_{p} \Ind{\clt^e_{\infty}} + l_{p})$ on every $\Ed$. 
\begin{itemize}
\item For the edge $e'$ such that $\ac(e') = 0$, as $\a^e_{\clt}(e') \leq \ac(e') $, the two functions $\ac\nabla (\phi_{p}+l_{p})(e')$ and $ \a^e_{\clt} \nabla (\phi_{p} \Ind{\clt^e_{\infty}} + l_{p})(e')$ are null. 
\item For the only pivot edge $e$ that $\ac(e) > 0, \a^e_{\clt}(e) = 0$, thanks to the second term of \Cref{lem:PivotEdge}, we have $\nabla (\phi_{p}+l_{p})(e) = 0$. Thus, the equation also establishes. 
\item For the edge that $\ac(e') > 0, \a^e_{\clt}(e') > 0$, we know that this implies that the two endpoints are on $\clt^e_{\infty}$ so that we have 
$$
\nabla (\phi_{p}+l_{p})(e') = \nabla (\phi_{p} \Ind{\clt^e_{\infty}} + l_{p})(e'),
$$
and $\ac(e') = \a^e_{\clt}(e')$, so the equation is also established.
\end{itemize}
$\ac\nabla (\phi_{p}+l_{p}) = \a^e_{\clt} \nabla (\phi_{p} \Ind{\clt^e_{\infty}} + l_{p})$ implies directly that $- \nabla \cdot \a^e_{\clt} \nabla (\phi_{p} \Ind{\clt^e_{\infty}} + l_{p}) = 0$ on $\Zd$, therefore, by the Liouville regularity, we obtain that $\phi^e_{p} = \phi_{p}\Ind{\clt^e_{\infty}}$ on $\clt^e_{\infty}$ up to a constant. 
\end{enumerate}
\end{proof}

A direct corollary of the third part of \Cref{lem:PivotEdge} is that $\g_{p} = \g_{p}^e$ on $\Ed$ when $\clt_{\infty} \neq \clt^e_{\infty}$, thus $K_R \star (\coa{\g_{p}} - \coa{\g_{p}^e}) = 0$. So, it suffices to consider $\sum_{e \in \Ed}\left\vert K_{R} \star \left(\coa{\g_{p}} - \coa{\g^e_{\e}} \right) \right\vert^2 (x)$ under the condition $\clt_{\infty} = \clt^e_{\infty}$. Then, we can reformulate the quantity in \cref{eq:QuantDoubleEnvironment1} as following:

\begin{equation*}
\begin{split}
K_R \star (\coa{\g_{p}} - \coa{\g_{p}^e})(x) &= K_R \star \left(\coa{\ac \D(\phi_{p} + l_{p})} - \coa{\a^e_{\clt} \D(\phi^e_{p} + l_{p})} \right)(x) \Ind{\clt_{\infty} = \clt^e_{\infty}}\\
&= \underbrace{K_R \star \coa{(\ac - \a^e_{\clt}) \D(\phi^e_{p} + l_{p})}(x) \Ind{\clt_{\infty} = \clt^e_{\infty}}}_{:= A_e(x)} + \underbrace{K_R \star \coa{\ac \D(\phi_{p} - \phi^e_{p})}(x) \Ind{\clt_{\infty} = \clt^e_{\infty}}}_{:= B_e(x)}.
\end{split}
\end{equation*}
In order to prove \cref{eq:QuantDoubleEnvironment1}, we study $A_e(x)$ and $B_e(x)$ separately.

\textit{Step 3: Case $\clt_{\infty} = \clt^e_{\infty}$, proof of $\sum_{e \in \Ed} |A_e(x)|^2 \leq \GO'_s(CR^{-d}) $.}
\Cref{lem:PivotEdge} helps us simplify the discussion on the case $\clt_{\infty} = \clt^e_{\infty}$ and following lemma carries the convolution to the cluster $\clt_{\infty}$.
\begin{lemma}\label{lem:Convolution}
For a kernel $K_R$ as in \Cref{prop:SpatialAverage2} and every $x \in \Rd$, there exists a function $\Gamma_{K,R}^x : \Zd \rightarrow \mathbb{R}^+$ such that for every function $\xi \in L^2(\clt_{\infty})$, we have 
\begin{equation}\label{eq:LemConvolution1}
(K_R \star \coa{\xi})(x) = \bracket{\Gamma_{K,R}^x , \xi}_{\clt_{\infty}},
\end{equation}
and we have the estimate 
\begin{equation}\label{eq:LemConvolution2}
\Gamma_{K,R}^x(z) \leq  \frac{2^d C_{K,R}}{R^d \left( \vert \frac{x-z}{R} \vert \vee 1 \right)^{\frac{d+1}{2}}}.
\end{equation}
\end{lemma}
\begin{proof}
We can do the calculation directly  
\begin{align*}
(K_R \star \coa{\xi})(x) &= \int_{\Rd} \coa{\xi}(y) K_R(x - y) \, dy \\
&= \int_{\Rd} \left(\sum_{z \in \clt_{\infty}}\Ind{y \in z + \cu}  \coa{\xi}(z)\right) K_R(x - y) \, dy \\
&= \sum_{z \in \clt_{\infty}} \left(\int_{\Rd} \Ind{y \in z + \cu} K_R(x - y) \, dy\right) \xi(z). 
\end{align*}
Thus we can define
\begin{equation}\label{eq:defGamma}
\Gamma_{K,R}^x(z) := \int_{y \in z + \cu} K_R(x - y) \, dy.
\end{equation} 
The estimate \cref{eq:LemConvolution2} comes directly from this expression and $K_{R} \leq \frac{C_{K,R}}{R^d \left( \vert \frac{x}{R} \vert \vee 1 \right)^{\frac{d+1}{2}}}$.
\end{proof}

We want to apply directly \Cref{lem:Convolution} to every random environment $\a^e$ to $A_e(x)$. We see that it suffices to study the case $e \in \Ed(\clt_{\infty})$, otherwise the condition $\clt_{\infty} = \clt^e_{\infty}$ will not be satisfied or $(\ac - \a^e_{\clt})\D(\phi^e_{p} + l_{p})$ will be $0$. Thus, $\supp\left( (\ac - \a^e_{\clt})\Dr{i}(\phi^e_{p} + l_{p})\right) \subset \clt_{\infty}$ and \Cref{lem:Convolution} works.
\begin{equation}\label{eq:AEstimate1}
\begin{split}
\sum_{e \in \Ed} \left| A_e \right|^2(x) &= \sum_{i=1}^d \sum_{e \in \Ed(\clt_{\infty})} \left \vert K_R \star \coa{(\ac - \a^e_{\clt}) \Dr{i}(\phi^e_{p} + l_{p})} \right \vert^2 (x) \Ind{\clt_{\infty} = \clt^e_{\infty}} \\
&= \sum_{i=1}^d\sum_{e \in \Ed(\clt_{\infty})} \left\vert \bracket{\Gamma_{K,R}^x ,\coa{(\ac - \a^e_{\clt})\D(\phi^e_{p} + l_{p})}_i }_{\clt_{\infty}}\right\vert^2 \Ind{\clt_{\infty} = \clt^e_{\infty}}\\
&\leq \sum_{\i=1}^d \sum_{z \in \clt_{\infty}} \vert \Gamma_{K,R}^x (\ac - \a^{\{z, z+\e_i\}}_{\clt})\Dr{i}(\phi^{\{z,z+\e_i\}}_{\e} + l_{p})\vert(z)^2 \Ind{\clt_{\infty} = \clt^e_{\infty}}.
\end{split}
\end{equation}
The third line comes from the fact that under the condition $\clt_{\infty} = \clt^e_{\infty}$, only one conductance $\a(e)$ and $\a^e(e)$ is different.

Using \cref{eq:LemConvolution2}, we know the part $\vert\Gamma_{K,R}^x\vert^2$ is integrable. It remains to estimate the vector field $\Theta$
\begin{equation}\label{eq:Theta}
\forall e = \{e_*, e^*\} \in \Ed, \qquad \Theta(e_*, e^*) := (\ac - \a^e_{\clt})\nabla(\phi^e_{p} + l_{p})(e_*, e^*),
\end{equation}
and prove a stochastic integrability for \cref{eq:AEstimate1}. Since the quantity $\Theta(e)$ plays an important role in our analysis and we will use it several times, we prove the following lemma:
\begin{lemma}\label{lem:Theta}
There exist two positive constants $s:=s(d,\pp,\Lambda)$ and $C:=C(d,\pp,\Lambda)$ such that 
\begin{equation}
\forall e \in \Ed, \qquad \vert \Theta(e) \vert \Ind{\clt_{\infty} = \clt^e_{\infty}} \leq \GO'_s(C).
\end{equation}
\end{lemma}
\begin{proof}
This estimate is very easy when $\ac^e(e) > 0$, where we use \Cref{prop:LocalAverage2} directly that
$$
\vert(\ac - \ac^e)\D(\phi^e_{p} + l_{p})\vert \Ind{\clt_{\infty} = \clt^e_{\infty}} \leq (1+\Lambda) \vert \ac^e\D(\phi^e_{p} + l_{p})\vert \Ind{\clt_{\infty} = \clt^e_{\infty}} \leq \GO'_s(C).
$$
The less immediate part comes from the case $\ac^e(e) = 0$ while $\ac(e) > 0$, where $\ac^e\D(\phi^e_{p} + l_{p}) = 0$ and cannot be used to dominate $\vert \Theta \vert(e)$. We treat this case differently: we denote by $e = \left\{e_*, e^* \right\}$, $\clt_{\infty} = \clt^e_{\infty}$ implies the existence of another open path $\gamma$ in $\clt^e_{\infty}$ connecting $e_*$ and $e^*$. This path can be chosen in $\clt_*(\cu_{\Pcu^e}(e_*))\cup\clt_*(\cu_{\Pcu^e}(e^*))$ applying \Cref{lem:FusionCluster} to the partition cube $\Pcu^e$.
\begin{equation}\label{eq:PathConnecting}
\begin{split}
\vert \Theta(e_*, e^*) \vert \Ind{\clt_{\infty} = \clt^e_{\infty}} &= \vert (\ac^e - \ac) \nabla(\phi^e_{p} + l_{p}) \vert (e_*, e^*) \Ind{\clt_{\infty} = \clt^e_{\infty}} \\
&\leq 2  \sum_{e' \in \gamma \subset \clt_*(\cu_{\Pcu^e}(e_*))\cup\clt_*(\cu_{\Pcu^e}(e^*))} \vert \nabla(\phi^e_{p} + l_{p})\vert (e')  \Ind{\clt_{\infty} = \clt^e_{\infty}} \\
&\leq 2 \vert \gamma \vert^{\frac{1}{2}} \left( \sum_{e' \in \gamma \subset \clt_*(\cu_{\Pcu^e}(e_*))\cup\clt_*(\cu_{\Pcu^e}(e^*))} \vert \nabla(\phi^e_{p} + l_{p})\vert (e') \Ind{\a^{e}(e') \neq 0} \right)^{\frac{1}{2}}\\
&\leq \vert \size(\cu_{\Pcu^e}(e_*))^d + \size(\cu_{\Pcu^e}(e^*))^d \vert^{\frac{1}{2}} \\
& \qquad \times\left( \sum_{e' \in \gamma \subset \clt_*(\cu_{\Pcu^e}(e_*))\cup\clt_*(\cu_{\Pcu^e}(e^*))} \vert \nabla(\phi^e_{p} + l_{p})\vert (e') \Ind{\a^{e}(e') \neq 0} \right)^{\frac{1}{2}} \\
&\leq \GO'_s(C).
\end{split}
\end{equation} 
\begin{figure}[h]
\centering
\includegraphics[scale=0.5]{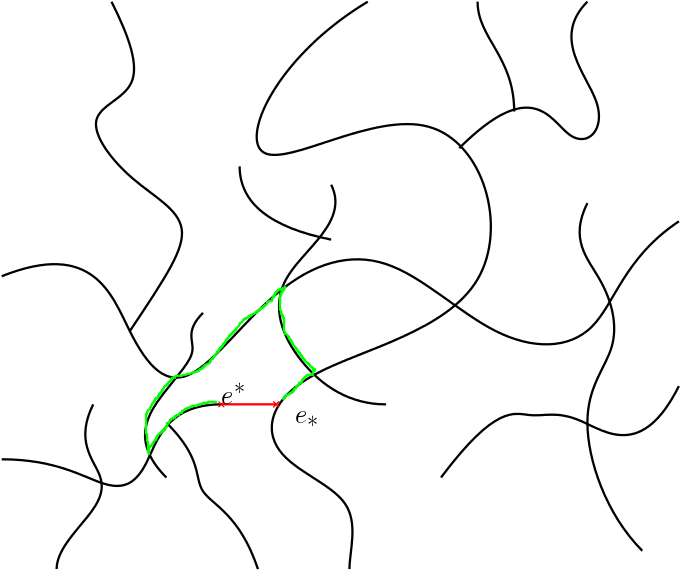}
\caption{This image shows the case $\clt_{\infty} = \clt^e_{\infty}$ and $\ac^e(e) = 0$ while $\ac(e) > 0$, so $e = \{e_*, e^*\}$ (segment in red) is an open bond in $\clt_{\infty}$ but not in $\clt^e_{\infty}$. The condition $\clt_{\infty} = \clt^e_{\infty}$ assures another open path (segment in green) in $\clt^e_{\infty}$ connecting $e_*$ and $e^*$.}
\end{figure}
We combine the two cases and prove \Cref{lem:Theta}.
\end{proof}

We conclude from \cref{eq:AEstimate1}, \cref{eq:LemConvolution2} and \Cref{lem:Theta} that 
\begin{align*}
\sum_{e \in \Ed} \left| A_e \right|^2(x) &\leq \sum_{\i=1}^d \sum_{z \in \clt_{\infty}} \vert \Gamma_{K,R}^x \vert^2(z) \vert \Theta\vert^2(z,z+\e_i) \Ind{\clt_{\infty} = \clt^e_{\infty}} \\
& \leq  \sum_{i=1}^d \sum_{z \in \Zd} \frac{4^dC^2_{K,R}}{R^{2d}(\left\vert \frac{x-z}{R}   \right\vert \vee 1)^{d+1}} \vert \Theta\vert^2(z,z+\e_i) \Ind{\clt_{\infty} = \clt^e_{\infty}} \\
& \leq \GO'_s\left( \frac{C_{K,R}^2 C(d,\pp,\Lambda)}{R^d}\right).  
\end{align*}

\textit{Step 4: Case $\clt_{\infty} = \clt^e_{\infty}$, proof of $\sum_{e \in \Ed} \vert B_e \vert^2 (x) \leq \GO'_s(CR^{-d})$.} This step is similar to that for $A_e$ but more technical. We define a space 
$$
\dot{H}^1(\clt_{\infty}) := \{v : \clt_{\infty} \rightarrow \mathbb{R}, \bracket{\nabla v, \nabla v }_{\Eda(\clt_{\infty})} < \infty\},
$$ 
and use the Green function on $(\clt_{\infty}, \Eda(\clt_{\infty}))$ \cite[Proposition 2.15]{dario2018optimal}:
\begin{proposition}[Green function on $\clt_{\infty}$]\label{prop:PoissonEquation}
Let $\a \in \Omega$ be an environment with an infinite cluster $\clt_{\infty}$ and $x,y \in \clt_{\infty}$, then there exists a constant $C:=C(d,\Lambda) < \infty$ and a Green function $G^{x,y} \in \dot{H}^1(\clt_{\infty}) $ such that 
$$
- \nabla \cdot \ac \nabla G^{x,y} = \delta_{y} - \delta_{x} \text{ on } \clt_{\infty}, 
$$
in the sense $\forall v \in \dot{H}^1(\clt_{\infty})$, we have 
$$
\bracket{\nabla G^{x,y}, \ac \nabla v}_{\Eda(\clt_{\infty})} = v(y) - v(x).
$$
In the case that $e = (x, y) \in \overrightarrow{\Eda}(\clt_{\infty})$, we note $G^{x,y} := G^e$. The Green function $G^{x,y}$ has the following properties
\begin{itemize}
\item \textit{Symmetry}: For every $x,y,x',y' \in \clt_{\infty}$, we have $G^{x,y}(y') - G^{x,y}(x') = G^{x',y'}(y) - G^{x',y'}(x)$.

\item \textit{Representation}: For every $v \in \dot{H}^1(\clt_{\infty})$, every vector field $\xi : \overrightarrow{\Eda}(\clt_{\infty}) \rightarrow \mathbb{R}$, and $u_{\xi} \in \dot{H}^1(\clt_{\infty}) $ such that
$$
\bracket{\nabla u_{\xi}, \ac \nabla v}_{\Eda(\clt_{\infty})} = \bracket{\xi, \nabla v}_{\Eda(\clt_{\infty})},
$$
we have the representation 
\begin{equation}\label{eq:PoissonRepresentation}
\nabla u_{\xi} = \sum_{e \in \Eda(\clt_{\infty})} \xi(e) \nabla G^e.
\end{equation}
In this formula, we give an arbitrary orientation for $e \in \Eda(\clt_{\infty})$, and the equation is well-defined.
\end{itemize}  
\end{proposition}
\begin{proof}
The proof of the existence and uniqueness of the function $G^{x,y}$ comes from the Lax-Milgram theorem on the space $\dot{H}^1(\clt_{\infty})$ where the conductance satisfies the quenched uniform ellipticity condition. The symmetry comes from testing the equation $- \nabla \cdot \ac \nabla G^{x,y} = \delta_{y} - \delta_{x}$ by $G^{x',y'}$ and testing the equation $- \nabla \cdot \ac \nabla G^{x',y'} = \delta_{y'} - \delta_{x'}$ by $G^{x,y}$ that 
$$
G^{x,y}(y') - G^{x,y}(x') = \bracket{\nabla G^{x,y}, \ac \nabla G^{x',y'}}_{\Eda(\clt_{\infty})} = G^{x',y'}(y) - G^{x',y'}(x).
$$
The final representation formula can be checked easily by the linear combination of the Green function.
\end{proof}

The proof of $\sum_{e \in \Ed} \vert B_e \vert^2 (x) \leq \GO'_s(CR^{-d})$ can be divided in $4$ steps. \newline
\textit{Step 4.1: Identification of $\D(\phi_{p} - \phi^e_{p})$ using Green function.}
We identify at first $\D(\phi_{p} - \phi^e_{p})$ by using the Green function $G^e$ introduced in \Cref{prop:PoissonEquation} and then estimate its size by this representation. We want to carry all the analysis on the geometry $(\clt_{\infty}, \Ed(\clt_{\infty}))$ and to claim the following lemma:

\begin{lemma}\label{lem:GreenRepresentaion}
We denote $e := \{e_*, e^*\} \in \Ed(\clt_{\infty})$, under the condition $\clt_{\infty} = \clt^e_{\infty}$,  then we have the following representation for $(\phi^e_{p} - \phi_{p})$, using \Cref{prop:PoissonEquation} and the definition $\Theta$ in \cref{eq:Theta},
\begin{equation}\label{eq:GreenRepresentation}
\nabla (\phi^e_{p} - \phi_{p})(\cdot) = \Theta(e_*, e^*)\nabla G^{e_*,e^*}(\cdot) \text{ on } \overrightarrow{\Eda}(\clt_{\infty}).
\end{equation}
\end{lemma}
\begin{proof}
Using the $\a$-harmonic equation and $\a^{e}$-harmonic equation for their correctors, we have at first
\begin{equation*}
- \nabla \cdot \ac \nabla (\phi_{p} + l_{p}) = - \nabla \cdot \a^e_{\clt} \nabla (\phi^e_{p} + l_{p}) \text{ on } \Zd,
\end{equation*}
then we obtain that 
\begin{equation}\label{eq:GreenRepresentation2}
-\nabla \cdot \ac \nabla(\phi_{p} - \phi^e_{p}) = - \nabla \cdot (\ac^e - \ac) \nabla(\phi^e_{p} + l_{p}) \text{ on } \Zd.
\end{equation}
Using the definition $\Theta(e_*, e^*) = (\ac - \ac^e)\nabla(\phi^e_{p} + l_{p})(e_*, e^*)$ and under the condition $\clt_{\infty} = \clt^e_{\infty}$, the right hand side of \cref{eq:GreenRepresentation2} equals to $\Theta(e)(\delta_{e^*} - \delta_{e_*})$. Moreover, since $e_*, e^* \in \clt_{\infty}$, \cref{eq:GreenRepresentation2} can be seen restricted on the cluster $(\clt_{\infty}, \Eda(\clt_{\infty}))$. Thus we solve the in $\dot{H}^1(\clt_{\infty})$ the equation 
$$
-\nabla \cdot \ac \nabla \tilde{w}^e = - \nabla \cdot (\ac^e - \ac) \nabla(\phi^e_{p} + l_{p}) \text{ on } \clt_{\infty},
$$
and by \Cref{prop:PoissonEquation}, it has a unique solution up to a constant that $\tilde{w}^e = \Theta(e_*, e^*)G^{e_*,e^*} $. 

Now we have $(\phi_{p} - \phi^e_{p})$ and $\tilde{w}^e$ solving the same equation, but we do not yet know that $(\phi_{p} - \phi^e_{p})$ belongs to $\dot{H}^1(\clt_{\infty})$. We hope to identify that $\phi_{p} - \phi^e_{p} = \tilde{w}^e$ and the argument is to use the Liouville regularity theorem: notice that $(\phi_{p} - \phi^e_{p} - \tilde{w}^e)$ is an $\a$-harmonic function on $\clt_{\infty}$ and $\bracket{\nabla \tilde{w}^e, \nabla \tilde{w}^e}_{\Eda} < \infty$ implies that $(\phi_{p} - \phi^e_{p} - \tilde{w}^e) \in \A_1$. We claim that it is in fact in $\A_0$ and prove by 
contradiction: suppose that $(\phi_{p} - \phi^e_{p} - \tilde{w}^e) \in \A_1 \backslash \A_0$, then by the Liouville regularity there exists $h \neq 0$ such that  
$$
\phi_{p} - \phi^e_{p} - \tilde{w}^e  = l_{h} + \phi_h.
$$
However, this implies that $\tilde{w}^e  = \phi_{p} - \phi^e_{p} - \phi_{h} - l_h$, so $\tilde{w}^e$ has an asymptotic linear increment at infinity, which contradicts the fact that $\tilde{w}^e \in \dot{H}^1(\clt_{\infty})$. In conclusion, we have $\phi_{p} - \phi^e_{p} - \Theta(e_*, e^*) G^{e_*,e^*} = c$ and we get \cref{eq:GreenRepresentation}. 
\end{proof}

\textit{Step 4.2: Carry the analysis on $(\clt_{\infty}, \Eda(\clt_{\infty}))$.}
Observing that we do the sum of $\ac \D( \phi_{p} - \phi^e_{p} )$, it suffices to do the sum over $\Eda(\clt_{\infty})$ and with the help of \Cref{lem:GreenRepresentaion} , we have the formula 
\begin{equation}\label{eq:BSimplification1}
\begin{split}
\sum_{e \in \Ed}  \vert B_e \vert^2(x) &= \sum_{i=1}^d \sum_{e \in \Eda(\clt_{\infty})} \vert K_R \star \coa{\ac \Dr{i}( \phi_{p} - \phi^e_{p} )}\vert^2(x)\Ind{\clt_{\infty} = \clt^e_{\infty}} \\
&= \sum_{i=1}^d \sum_{e \in \Eda(\clt_{\infty})}\left\vert K_R \star \coa{\ac \Dr{i}G^{e}} \right\vert^2(x) \Theta^2(e) \Ind{\clt_{\infty} = \clt^e_{\infty}} \\
&= \sum_{i=1}^d \sum_{e \in \Eda(\clt_{\infty})}\left\vert \bracket{\Gamma^x_{K,R}, \ac \Dr{i}G^{e}}_{\clt_{\infty}} \right\vert^2(x)\Theta^2(e) \Ind{\clt_{\infty} = \clt^e_{\infty}}.
\end{split}
\end{equation}

We analyze $\bracket{\Gamma^x_{K,R}, \ac \Dr{i}G^{e}}_{\clt_{\infty}} \Ind{\clt_{\infty} = \clt^e_{\infty}} $ by defining the notation that
\begin{equation*}
\Ind{\Ed^i}(e') = \left\{ 
\begin{array}{ll}
1 & \text{ if } \exists z \in \Zd \text{ such that } e' = \{z, z+ \e_i\} \\
0 & \text{ Otherwise }
\end{array},
\right.
\end{equation*}
and a vector field $\widetilde{\Gamma}^x_{K,R,i} : \Eda(\clt_{\infty}) \rightarrow \mathbb{R}$ that
\begin{equation}\label{eq:defGamma2}
\widetilde{\Gamma}^x_{K,R,i} := \Gamma^x_{K,R}(e'_{*,i}) \ac(e')\Ind{\Ed^i}(e'), \qquad  e'_{*,i}\in {\clt_{\infty}}  \text{ such that } e' = \left\{e'_*, e'_* + \e_i\right\}.
\end{equation}
 
Then, we can send $\bracket{\Gamma^x_{K,R}, \ac \Dr{i}G^{e}}_{\clt_{\infty}} \Ind{\clt_{\infty} = \clt^e_{\infty}} $ to the inner product of vector field on $\Eda(\clt_{\infty})$
\begin{equation}\label{eq:BSimplification2}
\begin{split}
\bracket{\Gamma^x_{K,R}, \ac \Dr{i}G^{e}}_{\clt_{\infty}} \Ind{\clt_{\infty} = \clt^e_{\infty}} = \bracket{\widetilde{\Gamma}^x_{K,R,i} , \nabla G^{e}}_{\Eda(\clt_{\infty})}\Ind{\clt_{\infty} = \clt^e_{\infty}}.
\end{split}
\end{equation}

\textit{Step 4.3 : Apply once again the representation with the Green function.}
Since $\widetilde{\Gamma}^x_{K,R,i} $ is defined on $\Eda(\clt_{\infty})$, we can apply \Cref{prop:PoissonEquation} to define $w_{\widetilde{\Gamma}^x_{K,R,i}} \in \dot{H}_1(\clt_{\infty})$ the solution of the equation 
\begin{equation}\label{eq:defWideGamma}
- \nabla \cdot \ac \nabla w_{\widetilde{\Gamma}^x_{K,R,i}} = - \nabla \cdot \widetilde{\Gamma}^x_{K,R,i}, \text{ on } \clt_{\infty},
\end{equation}
and it has a representation $\nabla w_{\widetilde{\Gamma}^x_{K,R,i}}(e) = \sum_{e' \in \Eda(\clt_{\infty})} \widetilde{\Gamma}^x_{K,R,i}(e') \nabla G^{e'}(e)$. We use the symmetry $\nabla G^{e'}(e) = \nabla G^{e}(e')$
\begin{equation}\label{eq:BSimplification3}
\nabla w_{\widetilde{\Gamma}^x_{K,R,i}}(e) = \sum_{e' \in \Eda(\clt_{\infty})} \widetilde{\Gamma}^x_{K,R,i}(e') \nabla G^{e}(e') = \bracket{\widetilde{\Gamma}^x_{K,R,i}, \nabla G^e}_{\Eda(\clt_{\infty})}.
\end{equation}
We combine \cref{eq:BSimplification1} \cref{eq:BSimplification2} and \cref{eq:BSimplification3} together and obtain that
\begin{equation}\label{eq:BSimplification4}
\sum_{e \in \Ed}  \vert B_e \vert^2(x)  \leq \sum_{i=1}^d \sum_{e \in \Eda(\clt_{\infty})} \vert \nabla w_{\widetilde{\Gamma}^x_{K,R,i}}\vert^2(e) \Theta^2(e) \Ind{\clt_{\infty} = \clt^e_{\infty}}.
\end{equation}

\textit{Step 4.4: Meyers' inequality and minimal scale.}
From the \cref{eq:defWideGamma}, we obtain a $\dot{H}^1(\clt_{\infty})$ estimate using \cref{eq:LemConvolution2}
\begin{align*}
&\bracket{\nabla w_{\widetilde{\Gamma}^x_{K,R,i}}, \ac\nabla w_{\widetilde{\Gamma}^x_{K,R,i}}}^{\frac{1}{2}}_{\Eda(\clt_{\infty})} = \bracket{\nabla w_{\widetilde{\Gamma}^x_{K,R,i}}, \widetilde{\Gamma}^x_{K,R,i}}^{\frac{1}{2}}_{\Eda(\clt_{\infty})} \\
\Longrightarrow &\norm{ \nabla w_{\widetilde{\Gamma}^x_{K,R,i}}}_{L^2(\Eda(\clt_{\infty}))} \leq  \Lambda \norm{ \widetilde{\Gamma}^x_{K,R,i}}_{L^2(\Eda(\clt_{\infty}))}  \leq  \Lambda \norm{ \Gamma^x_{K,R}}_{L^2(\Eda(\clt_{\infty}))} \leq C_{K,R}^2 R^{-\frac{d}{2}}.
\end{align*} 
Combining \cref{eq:BSimplification4} and the estimate on $\Theta(e)$, one may want to argue that 
$$
\sum_{e \in \Ed}  \vert B_e \vert^2(x)  \leq \sum_{i=1}^d \sum_{e \in \Eda(\clt_{\infty})} \GO'_s \left(\vert \nabla w_{\widetilde{\Gamma}^x_{K,R,i}}\vert^2(e)\right) \leq \GO'_s( C^2_{K,R} R^{-d}). 
$$
However, this argument is not correct since \cref{eq:OSum} does not work for our case where $w_{\widetilde{\Gamma}^x_{K,R,i}}$ is stochastic. A rigorous proof needs an argument as in \cite[Lemma 3.6]{dario2018optimal} using the minimal scale: We construct a collection of good cubes $\Gcu'$ such that not only Meyers' inequality \cite[Proposition 2.14]{dario2018optimal} is established, but also there exists $\epsilon(d, \pp, \Lambda) > 0$ and $C(d, \pp, \Lambda) < \infty$ for all $\cu \in \Gcu'$ 
\begin{equation}\label{eq:MinimalScale}
\frac{1}{|\cu|} \left(\int_{\Eda(\clt_{\infty} \cap \cu)} \Theta^{\frac{2(2+\epsilon)}{\epsilon}}(e)\right)^{\frac{\epsilon}{2 + \epsilon}} < C(d,\pp, \Lambda).
\end{equation}
Then we do the Calder\'on-Zygmund decomposition \Cref{prop:GeneralPartion} for $\Gcu'$ to obtain a partition of cubes $\mathcal{U}$, and apply Meyers' inequality for  \cref{eq:BSimplification4}

\begin{align*}
& \sum_{e \in \Ed}  \vert B_e \vert^2(x) 
\leq \sum_{i=1}^d \sum_{\cu \in \mathcal{U}} \sum_{e \in \Eda(\clt_{\infty} \cap \cu)} \vert \nabla w_{\widetilde{\Gamma}^x_{K,R,i}}\vert^2(e) \Theta^2(e) \Ind{\clt_{\infty} = \clt^e_{\infty}} 
\\
& \quad \leq \sum_{i=1}^d \sum_{\cu \in \mathcal{U}} |\cu|  \underbrace{\left(\frac{1}{|\cu|} \sum_{e \in \Eda(\clt_{\infty} \cap \cu)} \vert \nabla w_{\widetilde{\Gamma}^x_{K,R,i}}\vert^{2+\epsilon}(e) \right)^{\frac{2}{2+\epsilon}}}_{\text{Applying Meyers' inequality}}  \underbrace{\left(\frac{1}{|\cu|} \sum_{e \in \Eda(\clt_{\infty} \cap \cu)} \Theta^{\frac{2(2+\epsilon)}{\epsilon}}(e)  \right)^{\frac{\epsilon}{2+\epsilon}}}_{\leq C \text{ after } \cref{eq:MinimalScale}} \Ind{\clt_{\infty} = \clt^e_{\infty}}
\\
& 
\quad \leq C  \sum_{i=1}^d \sum_{\cu \in \mathcal{U}} |\cu| \left( \frac{1}{|\frac{4}{3}\cu|} \sum_{e \in \Eda(\clt_{\infty} \cap \frac{4}{3}\cu)} \vert \nabla w_{\widetilde{\Gamma}^x_{K,R,i}}\vert^{2}(e) + \left(\frac{1}{|\frac{4}{3}\cu|} \sum_{e \in \Eda(\clt_{\infty} \cap  \frac{4}{3}\cu)} \vert \widetilde{\Gamma}^x_{K,R,i}\vert^{2+\epsilon}(e) \right)^{\frac{2}{2+\epsilon}} \right).
\end{align*} 

The first term can be controlled by the $\dot{H}^1$ estimate for $w_{\widetilde{\Gamma}^x_{K,R,i}}$ that $$
\norm{ \nabla w_{\widetilde{\Gamma}^x_{K,R,i}}}^2_{L^2(\Eda(\clt_{\infty}))} \leq  \Lambda^2 \norm{ \Gamma^x_{K,R}}^2_{L^2(\Eda(\clt_{\infty}))} \leq C_K^2 R^{-d}.
$$
While for the second term, we can now apply \cref{eq:OSum} as $ \widetilde{\Gamma}^x_{K,R,i}$ is deterministic
\begin{multline*}
\sum_{i=1}^d \sum_{\cu \in \mathcal{U}} |\cu|  \left(\frac{1}{|\frac{4}{3}\cu|} \sum_{e \in \Eda(\clt_{\infty} \cap  \frac{4}{3}\cu)} \vert \widetilde{\Gamma}^x_{K,R,i}\vert^{2+\epsilon}(e) \right)^{\frac{2}{2+\epsilon}} 
\\
\leq \sum_{i=1}^d \sum_{\cu \in \mathcal{U}} |\cu|^{\frac{\epsilon}{2+\epsilon}}  \sum_{e \in \Eda(\clt_{\infty} \cap  \frac{4}{3}\cu)} \vert \widetilde{\Gamma}^x_{K,R,i}\vert^{2}(e)  \leq \GO'_s(C_{K,R}^2R^{-d}).
\end{multline*}
This concludes the proof.
\end{proof}

\subsection{Construction of flux conrrectors}\label{subsec:FluxCorrector}
In this part, we prove a Helmholtz-Hodge type decomposition for $\g_p$, which is another quantity $\S_p$ used in the further quantification of algorithm.  We recall that we use $\g_{p,i}$ to represent the $i$-th component of the vector field $\g_{p} : \Zd \rightarrow \Rd$ and the standard heat kernel is defined as $\Phi_R(x) := \frac{1}{(4\pi R^2)^{d/2}}\exp\left(-\frac{x^2}{4R^2}\right)$.

\begin{proposition}\label{prop: Helmholtz-Hodge }
For each $p \in \Rd$, almost surely there exists a vector field $\S_{p} : \Zd \rightarrow \mathbb{R}^{d \times d}$ called \textit{flux corrector} of $\g_{p}$, which takes values in the set of anti-symmetric matrices (that is, $\S_{p,ij} = - \S_{p,ji}$) and satisfying the following equations:
\begin{equation}\label{eq:defS}
\left\{
\begin{array}{ll}
\D^* \cdot \S_{p} &= \g_{p},\\
-\Delta \S_{p,ij} &= \Dr{j} \g_{p,i} - \Dr{i} \g_{p,j}, 
\end{array}
\right.
\end{equation}
where the first equation means that for every $i \in \{1,2 \cdots d\}$, $\sum_{j = 1}^{d} \Dr{j}^* \S_{p,ij} = \g_{p,i}.$ 

The quantity satisfies similar estimation as \cref{eq:LocalAverage1} and \cref{eq:SpatialAverage1}: there exist two positive constants $s:=s(d,\pp, \Lambda),$ $C:=C(d,\pp, \Lambda,s)$ such that 
\begin{align}\label{eq:LocalAverage3}
\forall 1 \leq i,j \leq d, \quad  \forall x \in \Zd, \quad \vert \Dr{k} \S_{p,ij} \vert(x) \leq \GO_s(C|p|),
\end{align}
and for the heat kernel $\Phi_{R}$, we have 
\begin{equation}\label{eq:SpatialAverage3}
\left\vert \Phi_{R} \star \coa{ \Dr{k}\S_{p,ij}}\right\vert (x)  \leq \GO_s( C |p| R^{-\frac{d}{2}}).
\end{equation}
\end{proposition}

\subsubsection{Heuristic analysis}
The following discussion gives a little heuristic analysis before a rigorous proof. In fact, if we define a field $H_{p}:\Zd \rightarrow \Rd$ such that 
\begin{equation}\label{eq:InformalS}
- \Delta H_{p,i} = \g_{p,i},
\end{equation}
where $- \Delta := - \nabla \cdot \nabla = \D^* \cdot \D$ is the discrete Laplace and then we define $\S_p$ such that 
\begin{equation}\label{eq:InformalDefS}
\S_{p,ij} = \Dr{j}H_{p,i} - \Dr{i}H_{p,j}.
\end{equation}
We see that this definition gives us a solution of \cref{eq:defS} since 
\begin{align*}
-\Delta \S_{p,ij} &= -\Delta \left( \Dr{j}H_{p,i} - \Dr{i}H_{p,j} \right) = \Dr{j}(-\Delta H_{p,i}) - \Dr{i}(-\Delta H_{p,j}) = \Dr{j} \g_{p,i} - \Dr{i} \g_{p,j}. \\
(\D^* \cdot \S_{p})_i &= \sum_{j=1}^d \Dr{j}^* \left( \Dr{j}H_{p,i} - \Dr{i}H_{p,j}\right) = -\Delta H_{p,i} - \Dr{i}\sum_{j=1}^d \Dr{j}^* H_{p,j} = \g_{p,i}.
\end{align*}
Here we use one property that $H_{p,i} = (-\Delta)^{-1}\g_{p,i}$ so that $H_{p}$ is also divergence free. This idea works on periodic homogenization problem \cite[Lemma 3.1]{kenig2012convergence}, but in our context, one key problem is to well define \cref{eq:InformalS}. In the present work, we apply an elementary probabilistic approach: Let $(S_k)_{k \geq 0}$ defines a lazy discrete time simple random walk on $\Zd$ with probability $\frac{1}{2}$ to stay unmoved and $\frac{1}{4d}$ to move towards one of the nearest neighbors on $\Zd$, and we use $(P_t)_{t\in \mathbb{N}}$ to define its semigroup, with the notation 
\begin{equation}\label{eq:defSemigroup}
\begin{split}
P_t(x,y) &:= P_t(y-x)= \Prb[S_t = y - x] \\
(P_t f)(x) &:= \sum_{y \in \Zd} P_t(x,y)f(y) = (\coa{P_t} \star \coa{f})(x), \qquad \forall f \in L^1(\Zd), 
\end{split}
\end{equation}
where $[\cdot]$ denotes a constant extension on every $z + \left(-\frac{1}{2}, \frac{1}{2}\right)^{d}$. 
Using the representation of the solution of harmonic function by a simple random walk   
$$
H_{p,i}(x) = \frac{1}{4d}\sum_{t = 0}^{\infty} (P_t \g_{p,i})(x),
$$ 
and we deduce from the definition of $\S_{\e}$ in \cref{eq:InformalDefS} 
$$
\S_{p,ij}(x) = \frac{1}{4d}\sum_{t = 0}^{\infty}\Dr{j}(P_t \g_{p,i})(x) - \Dr{i}(P_t \g_{p,j})(x).
$$
If we believe that $P_t$ is close to the heat kernel that $P_t(x,y) \simeq \frac{1}{(\pi t)^{d/2}}\exp\left(- \frac{|y-x|^2}{t}\right)$, and that the operator $\D$ helps to gain another factor of $t^{-\frac{1}{2}}$, then \Cref{prop:SpatialAverage2} would give us that $|\Dr{j}(P_t \g_{p,i})|(x) \leq \GO_s(t^{-\frac{1}{2}-\frac{d}{4}})$. We expect that this upper bound is sharp in general, and the fact that $\sum_{t = 1}^{\infty}t^{-\frac{1}{2}-\frac{d}{4}} = \infty$ prevents us from being able to define $\S_{p,ij}$ directly in dimension $d = 2$. Nevertheless we can make sense of 
\begin{equation}\label{eq:defDS}
(\Dr{k}\S_{p,ij})(x) = \frac{1}{4d}\sum_{t = 0}^{\infty}\Dr{k}\Dr{j}(P_t \g_{p,i})(x) - \Dr{k}\Dr{i}(P_t \g_{p,j})(x),
\end{equation}
because differentiating $P_t$ a second time will allow us to gain an extra factor of $t^{-\frac 1 2}$, and thus give us that that $|\Dr{k}\Dr{j}P_t \g_{p,i}| \leq \GO_s(C t^{-1-\frac{d}{4}})$. Then we can apply \cref{eq:OSum} to say that $\Dr{k}\S_{p,ij}$ is well-defined and prove other properties.

\subsubsection{Rigorous construction of $\D \S$}
\begin{proof}[Proof of \Cref{prop: Helmholtz-Hodge }]
We will give a rigorous proof that \cref{eq:defDS} gives a well-defined anti-symmetric valued vector field $\S_{e}$. The proof can be divided into three steps.

\textit{Step 1: Stochastic integrability of $\Dr{k}\Dr{j}(P_t \g_{p,i})$.}
In the first step, we prove that \cref{eq:defDS} makes sense, that is the part $\Dr{k}\Dr{j}(P_t \g_{p,i})(x) - \Dr{k}\Dr{i}(P_t \g_{p,j})(x)$ is summable. In the heuristic analysis, we compare $P_t$ with the heat kernel, which can be reformulated carefully by the local central limit theorem.
\begin{lemma}[Page 61, Exercise 2.10 of \cite{lawler2010random}]\label{lem:LCLT}
We denote by $\bar{P}_t(x,y) = \frac{1}{(\pi t)^{d/2}}\exp\left(- \frac{|y-x|^2}{t}\right)$, then there exists a positive constant $C(d)$, such that for all $t > 0$,
\begin{equation}\label{eq:LCLT}
\sup_{x \in \Zd} \left|\Dr{k}\Dr{j}P_t - \Dr{k}\Dr{j}\bar{P}_t \right|(x) \leq C(d)t^{-\frac{d+3}{2}}.
\end{equation}
\end{lemma}
\begin{proof}
The proof follows the idea in \cite[Theorem 2.3.5]{lawler2010random} and also relies on \cite[Lemmas 2.3.3 and 2.3.4]{lawler2010random} where we have 
$$
P_t(x) = \bar{P}_t(x) + V_t(x,r) + \frac{1}{(2\pi)^d t^{\frac{d}{2}}} \int_{|\theta| \leq r} e^{-\frac{ ix \cdot \theta}{\sqrt{t}}} e^{-\frac{|\theta|^2}{4}} F_t(\theta) \, d\theta,
$$
and there exits $\zeta > 0$ such that for every $0 < \theta, r \leq t^{\frac{1}{8}}$,  $V_t(x,r),F_t(\theta)$ satisfy 
$$
\left| F_t\right|(\theta) \leq \frac{|\theta|^4}{t}, \qquad \left| V_t(x,r) \right| \leq c(d)t^{-\frac{d}{2}} e^{- \zeta r^2}.
$$
We apply $\Dr{k}\Dr{j}$ with respect to $x$ and obtain that 
$$
\left|\Dr{k}\Dr{j}P_t -  \Dr{k}\Dr{j}\bar{P}_t\right|(x) = \left| \Dr{k}\Dr{j}V_t(x,r) + \frac{1}{(2\pi)^d t^{\frac{d}{2}}} \int_{|\theta| \leq r} \Dr{k}\Dr{j}e^{-\frac{ ix \cdot \theta}{\sqrt{t}}} e^{-\frac{|\theta|^2}{4}} F_t(\theta) \, d\theta \right|.
$$
We take $r = t^{\frac{1}{8}}$, then the term $V_t(x,r)$ has an error of exponential type 
$$
 \left| V_t(x,t^{\frac{1}{8}}) \right| \leq c(d)t^{-\frac{d}{2}} e^{- \zeta t^{\frac{1}{4}}} \leq c'(d) t^{- \frac{d+3}{2}}.
$$
So we focus on another part, by a simple finite difference calculus we have$\Dr{k}\Dr{j}e^{-\frac{ ix \cdot \theta}{\sqrt{t}}} \leq  \frac{2\theta}{\sqrt{t}}$. Moreover, we apply $\left| F_t\right|(\theta) \leq \frac{|\theta|^4}{t}$ and have 
\begin{align*}
\left| \frac{1}{(2\pi)^d t^{\frac{d}{2}}} \int_{|\theta| \leq r} \Dr{k}\Dr{j}e^{-\frac{ ix \cdot \theta}{\sqrt{t}}} e^{-\frac{|\theta|^2}{4}} F_t(\theta) \, d\theta \right| \leq \left| \frac{1}{(2\pi)^d t^{\frac{d+3}{2}}} \int_{|\theta| \leq r} |\theta|^5 e^{-\frac{|\theta|^2}{4}} \, d\theta \right| \leq C t^{-\frac{d+3}{2}}.
\end{align*}
This concludes the proof.
\end{proof}

We prove that \cref{eq:defDS} is well defined by showing that  
\begin{equation}\label{eq:DSIntergrable}
\Prb\text{-a.s  } \quad \forall 1\leq i,j,k \leq d,  \forall x \in \Zd, \quad \sum_{t = 0}^{\infty}\left | \Dr{k}\Dr{j}(P_t \g_{p,i})\right |(x) < \infty.
\end{equation}
We break this term into two
\begin{equation}\label{eq:DSIntergrableBreak}
\sum_{t = 0}^{\infty} \left| \Dr{k}\Dr{j}(P_t \g_{p,i}) \right|(x) = \underbrace{\sum_{t = 0}^{\infty} \left|\Dr{k}\Dr{j}(\bar{P}_t \g_{p,i})\right|(x)}_{\text{\cref{eq:DSIntergrableBreak}-a}} + \underbrace{\sum_{t = 0}^{\infty} \left|\Dr{k}\Dr{j}(P_t - \bar{P}_t) \g_{p,i}\right|(x)}_{\text{\cref{eq:DSIntergrableBreak}-b}}.
\end{equation}
$\bar{P}_t$ is better than $P_t$ since it is a standard heat kernel and we can do explicit calculation. We observe that 
$$
\forall t \geq 1, \forall y \in \Zd, \left|\Dr{k}\Dr{j}\bar{P}_t\right|(y) \leq \frac{C(d)}{t}\bar{P}_{2t}(y)=\frac{C(d)}{t(2\pi t)^{d/2}}\exp\left(-\frac{|y|^2}{2t}\right) \leq \frac{C(d)}{t^{\frac{d+2}{2}} \left( \vert \frac{y}{\sqrt{t}} \vert \vee 1 \right)^{\frac{d+1}{2}}}.
$$ 
Then $K_{\sqrt{t}} := \Dr{k}\Dr{j}\bar{P}_t$ is a kernel described in \Cref{prop:SpatialAverage2}  with the constant $C_{K,\sqrt{t}} := \frac{C(d)}{t}$, so we have
$$
\left| \Dr{k}\Dr{j}(\bar{P} \g_{p,i}) \right|(x) = \left| K_{\sqrt{t}} \star \coa{\g_{p,i}} \right|(x) \leq  \GO_s(Ct^{-1-\frac{d}{4}}).
$$ 
We put these estimates with \Cref{prop:LocalAverage2} in the \cref{eq:DSIntergrableBreak}-a and get  
\begin{align*}
\text{\cref{eq:DSIntergrableBreak}-a} \leq |\g_{p,i}|(x) + \sum_{t=1}^{\infty}\left| \Dr{k}\Dr{j}(P_t \g_{p,i}) \right|(x) \leq \GO_s(C) + \sum_{t=1}^{\infty}\GO_s(C t^{-1-\frac{d}{4}}) \leq \GO_s(C).
\end{align*}

On the other hand, to handle \cref{eq:DSIntergrableBreak}-b, we choose $\epsilon > 0$ and study at first 
\begin{align*}
\left|\Dr{k}\Dr{j}(P_t - \bar{P}_t) \g_{p,i}\right|(x) & \leq \left \vert \int_{|y| \leq t^{\frac{1}{2}+ \epsilon}} \coa{\Dr{k}\Dr{j}(P_t - \bar{P}_t)}(y) \coa{\g_{p,i}}(x-y) \, dy  \right \vert \\
& \qquad + \left \vert 4\int_{|y| \geq t^{\frac{1}{2} + \epsilon}} \left(\coa{P_t} + \coa{\bar{P}_t}\right)(y) \coa{\g_{p,i}}(x-y) \, dy  \right \vert \\
& \stackrel{\Cref{lem:LCLT}}{\leq} \int_{|y| \leq t^{\frac{1}{2}+ \epsilon}} \GO_s\left( \frac{C}{t^{\frac{d+3}{2}}}\right) \, dy +  4\int_{|y| \geq t^{\frac{1}{2} + \epsilon}} \left(\coa{P_t} + \coa{\bar{P}_t}\right)(y) \GO_s(C) \, dy \\
& \stackrel{\cref{eq:OSum}}{\leq} \GO_s\left(C t^{-(\frac{3}{2}-d\epsilon)} \right) + \GO_s\left(\int_{|y| \geq t^{\frac{1}{2} + \epsilon}} \left(\coa{P_t} + \coa{\bar{P}_t}\right)(y) \, dy\right).
\end{align*}
We divide the estimation into two terms since \Cref{lem:LCLT} is uniform but not optimal for the tail probability, which is of type sub-Gaussian so that the mass outside $|t|^{\frac{1}{2}+\epsilon}$ should be very small. By direct calculation, we have that $\int_{|y| \geq t^{\frac{1}{2} + \epsilon}} \coa{\bar{P}_t}(y) \, dy \leq C(d)e^{-t^{2\epsilon}}$ and by Hoeffding's inequality for the lazy simple random walk $(S_t)_{t \geq 0}$
\begin{align*}
\int_{|y| \geq t^{\frac{1}{2} + \epsilon}} \coa{P_t}(y) \, dy = \Prb\left[|S_t| \geq t^{\frac{1}{2} + \epsilon} \right] \leq 2\exp\left(- \frac{2 t^{1+2\epsilon}}{\V[S_t]} \right) \leq 2e^{-4t^{2\epsilon}}.
\end{align*}
Combining these tail event estimates and by choosing $\epsilon = \frac{1}{4d}$, we obtain that 
$$
\left|\Dr{k}\Dr{j}(P_t - \bar{P}_t) \g_{p,i}\right|(x) \leq \GO_s(C t^{-\frac{5}{4}}),
$$ 
and this concludes that $\text{\cref{eq:DSIntergrableBreak}-b} \leq \GO_s(C)$, so \cref{eq:DSIntergrable} holds, \cref{eq:LocalAverage3} holds and that $\Dr{k}\S_{p,ij}$ is well defined.
\begin{remark}\label{rmk:DSIntergrable}
In the proof, we also obtained one quantitative estimate of the following type: There exist two constants $s:=s(d,\pp,\Lambda)$, $ C:=C(d,\pp,\Lambda, s)$ such that for every random field $X:\Zd \rightarrow \mathbb{R}$ satisfying for every $ z \in \Zd, |X(z)| \leq \GO_s(\theta)$, we have
\begin{equation}\label{eq:DDgQuant}
\forall 1 \leq i,j,k \leq d, x \in \Zd, \left| \Dr{k}\Dr{j}(P_t - \bar{P}_t) X \right|(x) \leq \GO_s(C \theta t^{-\frac{5}{4}}).
\end{equation}
By a similar approach with the classical local central limit theorem \cite[Theorem 2.3.5]{lawler2010random}, we can also prove that 
\begin{equation}\label{eq:PgQuant}
\forall x \in \Zd, |P_t \g_{p,i}|(x) \leq \GO_s(C t^{-\frac{1}{4}}).
\end{equation}
\end{remark}

\textit{Step 2: Verification of \cref{eq:defS}.}
The verification of \cref{eq:defS} is direct thanks to \cref{eq:DSIntergrable}. We will also use the semigroup property that 
\begin{equation}\label{eq:Generator}
P_t(x) - P_{t-1}(x) = \frac{1}{4d}\Delta P_{t-1}(x).
\end{equation}
\begin{align*}
\left(\D^* \cdot \S_{p}\right)_i(x) &=  \frac{1}{4d} \sum_{t=0}^{\infty} \sum_{j=1}^{d} \Dr{j}^* \Dr{j} P_t \g_{p,i}(x) - \Dr{j}^* \Dr{i} P_t \g_{p,j}(x)\\
&= \frac{1}{4d} \sum_{t=0}^{\infty} (\underbrace{-\Delta P_t}_{\cref{eq:Generator}}) \g_{p,i}(x) -\Dr{i}P_t\underbrace{(\D^* \cdot \g_{p})}_{=0}(x)\\
&  = \sum_{t=0}^{\infty} (P_t - P_{t+1}) \g_{p,i}(x)\\
& = \g_{p,i}(x).
\end{align*}
In the last step, we use implicitly that $\lim_{t\rightarrow \infty}P_t \g_{p,i}(x) = 0 $ almost surely. This is true by Borel-Cantelli lemma and the estimation $|P_t \g_{p,i}|(x) \leq \GO_s(Ct^{-\frac{1}{4}})$ (see \cref{eq:PgQuant}):
$$
\sum_{t=1}^{\infty}\Prb[|P_t \g_{p,i}|(x) \geq \epsilon] \leq \sum_{t=1}^{\infty} \exp\left(- \left(C\epsilon t^{\frac{1}{4}}\right)^s\right) < \infty.
$$
The second part of \cref{eq:defS}, concerning $-\Delta \S_{p,ij}$, is easy to verify by a similar calculation
\begin{align*}
-\Delta \S_{p,ij}(x) &= \frac{1}{4d}\sum_{t = 0}^{\infty}\Dr{k}^*\Dr{k}\Dr{j}(P_t \g_{p,i})(x) - \Dr{k}^*\Dr{k}\Dr{i}(P_t \g_{p,j})(x)\\
&= \frac{1}{4d}\sum_{t = 0}^{\infty} (-\Delta P_t \Dr{j}\g_{p,i})(x) - (-\Delta
 P_t \Dr{i}\g_{p,j})(x)\\
&= (\Dr{j}\g_{p,i})(x) - (\Dr{i}\g_{p,j})(x).
\end{align*} 
Finally, by the definition, we can define $\S_{p,ij}$ just by integration of $\D \S_{p,ij}$ along a path. This construction does not depend on the choice of path since $\D \S_{p,ij} $ is a potential field.

\textit{Step 3: Estimation of $\left\vert \Phi_{R} \star \coa{\Dr{k}\S_{p,ij}}\right\vert (x)$.} 
This is a result of the convolution. Thanks to the \cref{eq:DSIntergrable}, we can apply Fubini lemma to  $\left\vert \Phi_{R} \star \coa{\Dr{k}\S_{p,ij}}\right\vert (x)$ that 
\begin{align*}
\left\vert \Phi_{R} \star \coa{\Dr{k}\S_{p,ij}}\right\vert (x) &= \left\vert \frac{1}{4d}\sum_{t = 0}^{\infty} \Phi_{R} \star \coa{\Dr{k}\Dr{j}P_t} \star \coa{\g_{p,i}} - \Phi_{R} \star \coa{\Dr{k}\Dr{i} P_t} \star \coa{\g_{p,j}}\right\vert (x)\\
&=  \left\vert \frac{1}{4d}\sum_{t = 0}^{\infty} \coa{\Dr{k}\Dr{j}P_t} \star \left(\Phi_{R} \star  \coa{\g_{p,i}}\right) - \coa{\Dr{k}\Dr{i} P_t} \star \left(\Phi_{R} \star  \coa{\g_{p,j}}\right)\right\vert (x)
\end{align*}
The main idea is that $\Phi_{R} \star  \coa{\g_{p,j}} \leq \GO_s\left(CR^{-\frac{d}{2}}\right)$ by \Cref{prop:SpatialAverage2} , then we repeat the main argument of stochastic integrability of $\Dr{k}\S_{p,ij}$ to get a better estimate. We focus on just one term:
\begin{equation}\label{eq:DSConvolEstimation}
\begin{split}
\left\vert \frac{1}{4d}\sum_{t = 0}^{\infty} \coa{\Dr{k}\Dr{j}P_t} \star \left(\Phi_{R} \star  \coa{\g_{p,i}}\right) \right\vert  & \leq \underbrace{\left\vert \frac{1}{4d}\sum_{t = 0}^{\infty} \coa{\Dr{k}\Dr{j}P_t - \Dr{k}\Dr{j}\bar{P}_t} \star \left(\Phi_{R} \star  \coa{\g_{p,i}}\right) \right\vert }_{\text{\cref{eq:DSConvolEstimation}-a}} \\
&\quad + \underbrace{\left\vert \frac{1}{4d}\sum_{t = 0}^{\infty} \left(\coa{\Dr{k}\Dr{j}\bar{P}_t} - \Dr{k}\Dr{j}\Phi_{\sqrt{\frac{t}{2}}}\right) \star \left(\Phi_{R} \star  \coa{\g_{p,i}}\right) \right\vert }_{{\text{\cref{eq:DSConvolEstimation}-b}}} \\
&\quad + \underbrace{\left\vert \frac{1}{4d}\sum_{t = 0}^{\infty} \Dr{k}\Dr{j}\Phi_{\sqrt{\frac{t}{2}}} \star \left(\Phi_{R} \star  \coa{\g_{p,i}}\right) \right\vert }_{{\text{\cref{eq:DSConvolEstimation}-c}}}.
\end{split}
\end{equation}
We treat the three terms one by one. For \cref{eq:DSConvolEstimation}-a, we apply \cref{eq:DDgQuant} with $X := \Phi_{R} \star  \coa{\g_{p,i}}$ and we use also \cref{eq:OSum}
\begin{equation}\label{eq:DSConvolEstimation-a}
\begin{split}
\text{\cref{eq:DSConvolEstimation}-a} \leq &  \frac{1}{d} \left \vert \Phi_{R} \star  \coa{\g_{p,i}} \right\vert (x) +  \frac{1}{4d}\sum_{t = 1}^{\infty} \left\vert \coa{\Dr{k}\Dr{j}(P_t - \bar{P}_t)} \star \left(\Phi_{R} \star  \coa{\g_{p,i}}\right) \right\vert (x) \\
\leq & \GO_s(CR^{-\frac{d}{2}}) + \sum_{t = 1}^{\infty}\GO_s(Ct^{-\frac{5}{4}}R^{-\frac{d}{2}})\\
\leq & \GO_s(CR^{-\frac{d}{2}}).
\end{split}
\end{equation}

For the term \text{\cref{eq:DSConvolEstimation}-b}, we observe that for every $y \in \Rd,$
$$
\left \vert \coa{\Dr{k}\Dr{j}\bar{P}_t}(y) - \Dr{k}\Dr{j}\Phi_{\sqrt{\frac{t}{2}}}(y) \right \vert = \left \vert \coa{\Dr{k}\Dr{j} \Phi_{\sqrt{\frac{t}{2}}}}(y) - \Dr{k}\Dr{j}\Phi_{\sqrt{\frac{t}{2}}}(y)\right \vert \leq \frac{C(d)}{t^{\frac{3}{2}}} \Phi_{\sqrt{t}}(y).
$$
We apply this estimate and use \cref{eq:OSum} to obtain that 
\begin{equation}\label{eq:DSConvolEstimation-b}
\begin{split}
& \left\vert \frac{1}{4d}\sum_{t = 0}^{\infty} \left(\coa{\Dr{k}\Dr{j}\bar{P}_t} - \Dr{k}\Dr{j}\Phi_{\sqrt{\frac{t}{2}}}\right) \star \left(\Phi_{R} \star  \coa{\g_{p,i}}\right) \right\vert (x)\\ 
\leq & \frac{1}{4d}\sum_{t = 0}^{\infty} \left \vert \coa{\Dr{k}\Dr{j}\bar{P}_t} - \Dr{k}\Dr{j}\Phi_{\sqrt{\frac{t}{2}}}\right\vert \star \left\vert  \Phi_{R} \star  \coa{\g_{p,i}} \right\vert (x) \\
\leq &  \frac{1}{d} \left \vert \Phi_{R} \star  \coa{\g_{p,i}} \right\vert (x) + \frac{1}{4d}\sum_{t = 1}^{\infty} \frac{C(d)}{t^{\frac{3}{2}}} \Phi_{\sqrt{t}} \star \left\vert  \Phi_{R} \star  \coa{\g_{p,i}} \right\vert (x) \\
\leq & \GO_s(CR^{-\frac{d}{2}}) + \sum_{t=1}^{\infty} \GO_s(Ct^{-\frac{3}{2}}R^{-\frac{d}{2}}) \\
\leq & \GO_s(CR^{-\frac{d}{2}}).
\end{split}
\end{equation}

For the last term $\left\vert \frac{1}{4d}\sum_{t = 0}^{\infty} \Dr{k}\Dr{j}\Phi_{\sqrt{\frac{t}{2}}} \star \left(\Phi_{R} \star  \coa{\g_{p,i}}\right) \right\vert (x)$, we use the property of semigroup, the linearity of the finite difference operator and we apply \Cref{prop:SpatialAverage2} to the kernel $\Dr{k}\Dr{j}\Phi_{\sqrt{\frac{t}{2}}}$
\begin{equation}\label{eq:DSConvolEstimation-c}
\begin{split}
\left\vert \frac{1}{4d}\sum_{t = 0}^{\infty} \Dr{k}\Dr{j}\Phi_{\sqrt{\frac{t}{2}}} \star \left(\Phi_{R} \star  \coa{\g_{p,i}}\right) \right\vert (x) &= \left\vert \frac{1}{4d}\sum_{t = 0}^{\infty} \Dr{k}\Dr{j} \left(\Phi_{\sqrt{\frac{t}{2} + R^2}} \star  \coa{\g_{p,i}}\right) \right\vert (x) \\
& = \left\vert \frac{1}{4d}\sum_{t = 0}^{\infty} \left(\Dr{k}\Dr{j}\Phi_{\sqrt{\frac{t}{2} + R^2}} \right) \star  \coa{\g_{p,i}} \right\vert (x) \\
& \leq \sum_{t = 0}^{\infty} \GO_s\left( C \left(\frac{t}{2} + R^2\right)^{-1-\frac{d}{2}} \right) \\
& \leq \GO_s(C R^{-\frac{d}{2}}).
\end{split}
\end{equation} 

This concludes the proof as we put the three estimates \cref{eq:DSConvolEstimation-a},\cref{eq:DSConvolEstimation-b} and \cref{eq:DSConvolEstimation-c} in \cref{eq:DSConvolEstimation} and \cref{eq:SpatialAverage3}.
\end{proof}

\subsection{$L^q, L^{\infty}$ estimate of $\g_p, \D \S_p$}
\begin{proposition}\label{prop:LqFlux}
Let $F : \Zd \rightarrow \mathbb{R}$ stand for a random field of the form $\g_{p_i}$  or $\Dr{k}\S_{p,ij}$, for any $p \in \Rd, 1\leq i,j,k \leq d$. This field satisfies: 

There exist three positive constants $s:=s(d,\pp, \Lambda), k:= k(d,\pp,\Lambda)$ and $C:=C(d,\pp, \Lambda,s)$ such that for each $q \in [1, \infty)$,
\begin{equation}
\left(R^{-d} \int_{\clt_{\infty} \cap B_R} \vert F - (F)_{\clt_{\infty} \cap B_R} \vert^q \right)^{\frac{1}{q}} \leq \left\{
	\begin{array}{ll}
	\GO_s(C\vert p\vert q^k \log^{\frac{1}{2}}(R) )  & d=2, \\	
	\GO_s(C\vert p\vert q^k )  & d=3,
	\end{array}
\right.
\end{equation}
and for each $x,y \in \Zd$,
\begin{equation}
\vert F(x) - F(y) \vert \Ind{x,y \in \clt_{\infty}} \leq \left\{
	\begin{array}{ll}
	\GO_s(C\vert p\vert \log^{\frac{1}{2}}\vert x - y\vert )  & d=2, \\	
	\GO_s(C\vert p\vert )  & d=3.
	\end{array}
\right.
\end{equation}
\end{proposition}
\begin{proof}
Similar to \cite[Theorems 1 and 2]{dario2018optimal}, these estimates are the results of local estimate and spatial average estimates proved in \cref{eq:LocalAverage2}, \cref{eq:SpatialAverage2}, \cref{eq:LocalAverage3} and \cref{eq:SpatialAverage3} by applying a heat kernel type multi-scale Poincar\'e's inequality. We refer to \cite[Sections 4 and 5]{dario2018optimal}.
\end{proof}

\section{Two-scale expansion on the cluster}
	\label{sec:TwoScale}
In this part, we prove \Cref{thm:TwoScale} which is the heart of all the analysis of our algorithm as stated in \Cref{subsec:Idea}. Here we prove a more detailed version of the theorem.
\begin{proposition}[Two-scale expansion on percolation]\label{prop:TwoScaleDetail}
Under the same context of \Cref{thm:TwoScale}, there exist three random variables $\X, \Y_1, \Y_2$ satisfying
\begin{align*}
\X \leq \GO_1(C(d, \pp, \Lambda)m), \qquad \Y_1  \leq \GO_s\left(C(d,\pp,\Lambda, s)\ell(\ec)m^{\frac{1}{s}}\right), \qquad \Y_2  \leq \GO_s\left(C(d,\pp,\Lambda, s)\ec^{\frac{d}{2}}m^{\frac{1}{s}}\right),
\end{align*}
and we have the estimate 
\begin{equation*}
\begin{split}
\norm{\nabla (w - v) \Ind{\a \neq 0}}_{L^2(\cltm)}
&\leq C(d,\Lambda) \left(\norm{\D \vb}_{L^2(\cu_m)} \left(3^{-\frac{m}{2}} \ell^{\frac{1}{2}}(\ec)\X^{d} + 3^{-\frac{m}{2}} \ell^{-\frac{1}{2}}(\ec)\Y_1\X^{d} + \mu \Y_1 + \Y_2\X^d  \right)\right.\\
& \qquad + \norm{\D \vb}^{\frac{1}{2}}_{L^2(\cu_m)} \norm{\D^* \D \vb}^{\frac{1}{2}}_{L^2(\itr(\cu_m))}\left(\ell^{\frac{1}{2}}(\ec)\X^{d} + \ell^{-\frac{1}{2}}(\ec)\Y_1\X^{d} \right) \\
& \qquad + \left. \norm{\D^* \D \vb}_{L^2(\itr(\cu_m))} \Y_1 \X^d \right).
\end{split}
\end{equation*}
\end{proposition}

\subsection{Main part of the proof}
The main idea of the proof is to use the quantities $\{\phi_{\e_k}\}_{k=1, \ldots, d}$ and $\{\S_{\e_k, ij}\}_{i,j,k=1, \ldots, d}$ analyzed in previous work and in \Cref{sec:Flux}, under the condition $\cu_m \in \Pcu_*$. We do some simple manipulations at first. Throughout the proof, we use the notation $h := v - w$.
\begin{proof}
\textit{Step 1: Setting up.} We define a modified coarsened function $\tilde{h}$
\begin{equation}
\label{eq:TestFunction}
\tilde{h}(x) = \left\{
	\begin{array}{ll}
	h(x) 	& x \in \clt_*(\cu_m), \\
	\coa{h}_{\Pcu}(x) & x \in \cu_m \backslash \clt_*(\cu_m) , \dist(\cu_{\Pcu}(x), \partial \cu_m) \geq 1,\\
	0 	& x \in  \cu_m \backslash \clt_*(\cu_m) , \dist(\cu_{\Pcu}(x), \partial \cu_m) = 0.
	\end{array}
\right.
\end{equation}
We put it as a test function in \cref{eq:TwoScale} 
\begin{align*}
\bracket{\tilde{h}, (\mu_{\clt, m}^2 - \nabla \cdot \acf{m} \nabla) v}_{\itr(\cu_m)} = \bracket{\tilde{h}, (\mu_{\clt,m}^2 - \nabla \cdot \ab \nabla) \vb }_{\itr(\cu_m)}. 
\end{align*}
Since $\tilde{h} \in C_0(\cu_m)$, we can apply the formula \cref{eq:IPP} and get 
\begin{align}
\bracket{\mu_{\clt,m}\tilde{h}, \mu_{\clt,m} v}_{\cu_m} + \bracket{\nabla \tilde{h},  \acf{m} \nabla v}_{\cu_m} &= \bracket{\mu_{\clt,m}\tilde{h}, \mu_{\clt,m} \vb}_{\cu_m} + \bracket{\nabla \tilde{h},  \ab \nabla \vb }_{\cu_m}.  \label{eq:TestIn} 
\end{align}
We subtract a term of $w$ on the two sides to get 
\begin{multline*}
\bracket{\mu_{\clt,m}\tilde{h}, \mu_{\clt,m}(v-w)}_{\cu_m} +  \bracket{\nabla \tilde{h},  \acf{m} \nabla (v-w)}_{\cu_m} 
\\
 = \bracket{\mu_{\clt,m}\tilde{h}, \mu_{\clt,m} (\vb - w)}_{\cu_m} + \bracket{\nabla \tilde{h},  \ab \nabla \vb - \a_{\clt, m} \nabla w }_{\cu_m}. 
\end{multline*}
We put $v-w = h$ into the identity and obtain that 
\begin{multline}\label{eq:TestInOut}
\bracket{\mu_{\clt,m}\tilde{h}, \mu_{\clt,m}h}_{\cu_m} +  \bracket{\nabla \tilde{h},  \acf{m} \nabla h}_{\cu_m} 
\\
 = \bracket{\mu_{\clt,m}\tilde{h}, \mu_{\clt,m} (\vb - w)}_{\cu_m} + \bracket{\nabla \tilde{h},  \ab \nabla \vb - \acf{m}\nabla w }_{\cu_m}.
\end{multline}

\textit{Step 2: Restriction tricks.}
There are three observations:
\begin{itemize}
\item{Observation 1.} The effect of $\mu_{\clt,m}$ restricts the inner product to $\cltm$, and on $\cltm$ we have $\tilde{h} = h$ by \cref{eq:TestFunction}. Thus we have
$$
\bracket{\mu_{\clt, m}\tilde{h}, \mu_{\clt, m}h}_{\cu_m} = \mu^2\bracket{h, h}_{\cltm}, \qquad  \bracket{\mu_{\clt, m}\tilde{h}, \mu_{\clt, m} (\vb - w)}_{\cu_m} =  \mu^2\bracket{h,  \vb - w}_{\cltm}.
$$
\item{Observation 2.} The definition of $\acf{m}$ also restricts the inner product on $\Ed(\cltm)$ and we have 
$$
\bracket{\nabla \tilde{h},  \acf{m} \nabla h}_{\cu_m} = \bracket{\nabla h,  \a \nabla h}_{\cltm},$$ 
as $\acf{m} = 0$ outside $\Ed(\cltm)$ by \cref{eq:LocalMaskOperation}.
 
\item{Observation 3.} This step is the key where we gain much in the estimate and where we use the condition $\cu_m \in \Pcu_*$. We apply the formula \cref{eq:DivergenceForm} to $\bracket{\nabla \tilde{h},  \ab \nabla \vb - \acf{m}\nabla w }_{\cu_m}$ to obtain that
\begin{align*}
\bracket{\nabla \tilde{h},  \ab \nabla \vb - \acf{m}\nabla w }_{\cu_m} &= \bracket{\tilde{h},  -\nabla \cdot (\ab \nabla \vb - \acf{m}\nabla w) }_{\itr(\cu_m)}\\
&=\bracket{\tilde{h}, \D^* \cdot (\ab \D \vb - \acf{m}\D w)}_{\itr(\cu_m)} \\
& =\bracket{\tilde{h}, \D^* \cdot (\ab \D \vb - \ac \D w)}_{\itr(\cu_m)} \\
& \qquad + \bracket{\tilde{h}, \D^* \cdot (\ac - \acf{m}) \D w}_{\itr(\cu_m)}.
\end{align*}
We use the condition $\cu_m \in \Pcu_*$, which implies that $\cltm \subset \clt_{\infty}$ and 
$$
\supp \left(\D^* \cdot (\ac - \acf{m}) \D w \right) \subset (\clt_{\infty} \cap \cu_m) \backslash \cltm.
$$ 
In \Cref{def:SmallCluster} and \Cref{lem:SmallCluster}, we prove that $ (\clt_{\infty} \cap \cu_m) \backslash \cltm $ is the union of small clusters contained in the partition cubes $\cu_{\Pcu}$ with distance $1$ to $\partial \cu_m$, where $\tilde{h}$ equals $0$. Therefore, we obtain that 
 \begin{align*}
\bracket{\tilde{h}, \D^* \cdot (\ab \D \vb - \acf{m}\D w)}_{\itr(\cu_m)} = & \bracket{\tilde{h}, \D^* \cdot (\ab \D \vb - \ac \D w)}_{\itr(\cu_m)}.
\end{align*}
Using an identity 
$$
\D^* \cdot \left(\ac \D w - \ab \D \vb \right) = \D^* \cdot \Fi,
$$
which will be proved later in \Cref{lem:defF} and $\Fi$ is a vector field $\Fi : \Zd \rightarrow \Rd$, we conclude
\begin{align*}
\bracket{\nabla \tilde{h},  \ab \nabla \vb - \acf{m}\nabla w }_{\cu_m} \Ind{\cu_m \in \Pcu_*} &= \bracket{\tilde{h}, - \D^* \cdot \Fi}_{\itr(\cu_m)} \Ind{\cu_m \in \Pcu_*}\\
&= -\bracket{\D \tilde{h},   \Fi }_{\cu_m} \Ind{\cu_m \in \Pcu_*}
\end{align*}
\end{itemize}
Combining all these observations, we transform \cref{eq:TestInOut} to 
\begin{equation}\label{eq:TestInOut2}
\begin{split}
\left(\mu^2\bracket{h, h}_{\cltm} + \bracket{\nabla h,  \a \nabla h}_{\cltm}\right)  & \Ind{\cu_m \in \Pcu_*}\\
& = \left( \mu^2\bracket{h,  \vb - w}_{\cltm} - \bracket{\D \tilde{h},   \Fi }_{\cu_m}\right)  \Ind{\cu_m \in \Pcu_*}. 
\end{split}
\end{equation}
Using H\"older's inequality and Young's inequality, we obtain that 
\begin{equation} \label{eq:TestInOut3}
 \bracket{\nabla h,  \a \nabla h}_{\cltm} \Ind{\cu_m \in \Pcu_*} \leq \left(\frac{\mu^2}{4}\norm{\vb - w}^2_{L^2(\cltm)} +  \norm{\D \tilde{h}}_{L^2(\cu_m)} \norm{ \Fi }_{L^2(\cu_m)}\right)\Ind{\cu_m \in \Pcu_*}. 
\end{equation}

\begin{figure}[h]
\includegraphics[scale=0.6]{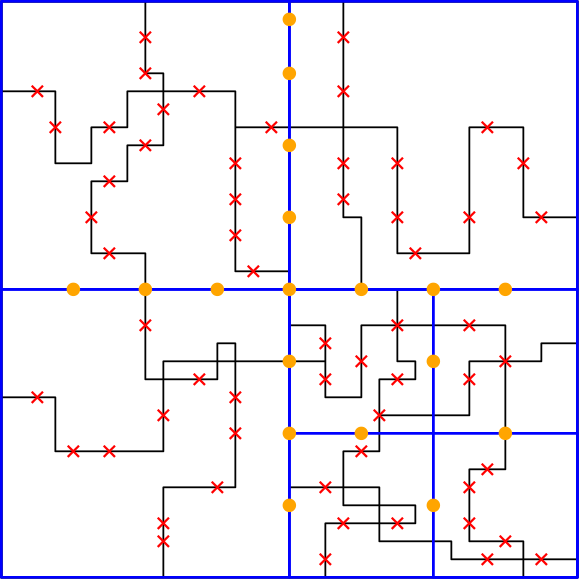}
\centering
\caption{The figure shows the sources contributing to $\norm{\nabla \tilde{h}}_{L^2(\cu_m)}$. The black segments represent the cluster $\cltm$ while the blue segments represent the partition of good cubes. Using the coarsened function, we see that the quantity can be controlled by the sum of three terms: the difference between $\nabla \tilde{h}$ and $\coa{h}_{\Pcu,\cu_m}$ near the cluster $\clt_*(\cu_m)$, marked with red cross in the image; the gradient $\coa{h}_{\Pcu,\cu_m}$ at the interface of different partition cubes $\cu_{\Pcu}$, marked with orange disk.}
\end{figure}

\textit{Step 3: Study of $\norm{\D \tilde{h}}_{L^2(\cu_m)}$.} 
The next step is to estimate the size of $\norm{\D \tilde{h}}_{L^2(\cu_m)}$. Since $\tilde{h} \in C_0(\cu_m)$, we have that $\norm{\D \tilde{h}}_{L^2(\cu_m)} = \norm{\nabla \tilde{h}}_{L^2(\cu_m)}$. We use the function $\coa{h}_{\Pcu,\cu_m}$ defined in \cref{eq:CoarsenedFunctionLocal}
\begin{equation}
\coa{h}_{\Pcu,\cu_m}(x) = \left\{
	\begin{array}{ll}
	\coa{h}_{\Pcu}(x) & \dist(\cu_{\Pcu}(x), \partial \cu_m) \geq 1,\\
	0 	& \dist(\cu_{\Pcu}(x), \partial \cu_m) = 0.
	\end{array}
\right.
\end{equation}
as a function to do comparison and apply \cref{eq:NaiveInequality} 
\begin{align*}
\norm{\nabla \tilde{h}}_{L^2(\cu_m)} &= \norm{\nabla (\tilde{h} - \coa{h}_{\Pcu,\cu_m})}_{L^2(\cu_m)} + \norm{\nabla \coa{h}_{\Pcu,\cu_m}}_{L^2(\cu_m)} \\
& \leq 2d \norm{\tilde{h} - \coa{h}_{\Pcu,\cu_m}}_{L^2(\cu_m)} + \norm{\nabla \coa{h}_{\Pcu,\cu_m}}_{L^2(\cu_m)}\\
& \leq 2d \norm{h - \coa{h}_{\Pcu,\cu_m}}_{L^2(\itr(\cu_m)) \cap \cltm} + \norm{\nabla \coa{h}_{\Pcu,\cu_m}}_{L^2(\cu_m)}
\end{align*}
The last step is correct since $\tilde{h}$ and $\coa{h}_{\Pcu,\cu_m}$ coincide at the boundary and also on the part $\itr(\cu_m)\backslash \cltm$. We define 
\begin{equation}\label{eq:defX}
\X := \max_{x \in \cu_m} \size(\cu_{\Pcu}(x)),
\end{equation}
which can be estimated using \cref{eq:OMax} and \cref{eq:MinimalScale} as
\begin{equation}\label{eq:QuantX}
\X \leq \GO_1(C(d, \pp, \Lambda)m),
\end{equation}
and apply \Cref{prop:ErrorCoarse}
\begin{align*}
\norm{\nabla \tilde{h}}_{L^2(\cu_m)} &\leq 2d \norm{h - \coa{h}_{\Pcu,\cu_m}}_{L^2(\itr(\cu_m)) \cap \cltm} + \norm{\nabla \coa{h}_{\Pcu,\cu_m}}_{L^2(\cu_m)}\\
& \leq  C \left( \sum_{\{x,y\} \in \Ed(\cltm)} \size(\cu_{\Pcu}(x))^{2d} \vert \nabla h\vert^2(x,y)\Ind{\a(x,y) \neq 0}\right)^{\frac{1}{2}}\\
& \qquad +  C\left( \sum_{\{x,y\} \in \Ed(\cltm} \left(\size(\cu_{\Pcu}(x))^{2d-1} + \size(\cu_{\Pcu}(y))^{2d-1} \right) \vert \nabla h\vert^2(x,y)\Ind{\a(x,y) \neq 0}\right)^{\frac{1}{2}}\\
& \leq C \X^{d} \left( \sum_{\{x,y\} \in \Ed(\cltm)}  \vert \nabla h\vert^2(x,y)\Ind{\a(x,y) \neq 0} \right)\\
& = C \X^d \norm{\nabla h \Ind{\a \neq 0}}_{L^2(\cltm)}.
\end{align*}
We put it back to \cref{eq:TestInOut3}   
\begin{align*}
\bracket{\nabla h,  \a \nabla h}_{\cltm} \Ind{\cu_m \in \Pcu_*} 
&\leq \left(\frac{\mu^2}{4}\norm{\vb - w}^2_{L^2(\cltm)} +  C \X^d \norm{\nabla h \Ind{\a \neq 0}}_{L^2(\cltm)} \norm{ \Fi }_{L^2(\cu_m)}\right) \Ind{\cu_m \in \Pcu_*}
\end{align*}
and use Young's inequality, finally we get 
\begin{equation}\label{eq:TestInOut4}
\norm{\nabla h \Ind{\a \neq 0}}_{L^2(\cltm)} \Ind{\cu_m \in \Pcu_*}  \leq C(d)\Lambda\left(\mu\norm{\vb - w}_{L^2(\cltm)} +   \X^{d}  \norm{ \Fi }_{L^2(\cu_m)}\right)\Ind{\cu_m \in \Pcu_*} .
\end{equation}

\textit{Step 4: Quantification.}
The remaining part is to quantify the estimate $\norm{\vb - w}_{L^2(\cltm)}\Ind{\cu_m \in \Pcu_*}$ and $\norm{ \Fi }_{L^2(\cu_m)}\Ind{\cu_m \in \Pcu_*}$.

The two random variables used in the estimation are defined as 
\begin{equation}\label{eq:defY}
\begin{split}
\Y_1 &:= \max_{1\leq i,j,k \leq d, \dist(x, \cu_m) \leq 1} \left|\phim_{\e_k}(x)\right| + \left|\Sm_{\e_k,ij}(x)\right|\\
\Y_2 &:= \max_{1\leq i,j,k \leq d,  \dist(x, \cu_m) \leq 1} \left|\Phi_{\ec^{-1}} \star \Dr{i}\coa{\phi_{\e_k}}^{\eta}_{\Pcu}(x)\right| + \left| \Phi_{\ec^{-1}} \star [\Dr{j}^{*}\S_{\e_k,ij}](x)\right|,
\end{split}
\end{equation}
and they have estimates following \cref{eq:OMax}, \cref{eq:SpatialAverage1}, \cref{eq:SpatialAverage3} and also \Cref{prop:LqFlux}, \Cref{prop:LpCluster} that there exists $0< s(d,\pp,\Lambda) < \infty$ and $0 < C(d, \pp, \Lambda) < \infty$ such that
\begin{equation}\label{eq:QuantY}
\begin{split}
\Y_1  \leq \GO_s\left(C(d,\pp,\Lambda, s)\ell(\ec)m^{\frac{1}{s}}\right) \qquad \Y_2  \leq \GO_s\left(C(d,\pp,\Lambda, s)\ec^{\frac{d}{2}}m^{\frac{1}{s}}\right).
\end{split}
\end{equation}
For $\norm{\vb - w}_{L^2(\cltm)}\Ind{\cu_m \in \Pcu_*}$, we have 
\begin{equation}\label{eq:QuantWV}
\begin{split}
\norm{w - \vb}^2_{L^2(\cu_m)}\Ind{\cu_m \in \Pcu_*} &= \norm{\sum_{j=1}^d (\Upsilon \Dr{j}\vb) \phim_{\e_j}}^2_{L^2(\cu_m)}\Ind{\cu_m \in \Pcu_*}\\
& \leq d \underbrace{\left(\max_{1 \leq j \leq d, x\in \cu_m} \phim_{\e_j}(x)\right)^2}_{\leq \Y_1^2} \sum_{j=1}^d \sum_{x \in \cu_m} \left((\Upsilon \Dr{j}\vb)(x)\right)^2 \\
& \leq d \Y_1^2 \norm{\D \vb}^2_{L^2(\cu_m)}.
\end{split}
\end{equation}

For $\norm{ \Fi }_{L^2(\cu_m)}\Ind{\cu_m \in \Pcu_*}$, we use the formula \cref{eq:defF}
\begin{equation}\label{eq:QuantFDecom}
\begin{split}
\norm{\Fi}^2_{L^2(\cu_m)} &\leq C(d)\left( \underbrace{\sum_{i=1}^d\norm{(1 - \Upsilon)\left(\ac - \ab \right)\left(\Dr{i}\vb\right)}^2_{L^2(\cu_m)}}_{\text{\cref{eq:QuantFDecom}-a}}  + \underbrace{\sum_{i,k=1}^d \norm{\phim_{\e_k}(\cdot + \e_k) \ac\Dr{i}\left(\Upsilon \Dr{k}\vb\right)}^2_{L^2(\cu_m)}}_{\text{\cref{eq:QuantFDecom}-b}} \right.\\
& \qquad + \underbrace{\sum_{i,j,k=1}^d \norm{\Sm_{\e_k,ij}(\cdot-\e_j)\Dr{j}^*\left(\Upsilon\Dr{k}\vb\right)}^2_{L^2(\cu_m)}}_{\text{\cref{eq:QuantFDecom}-c}}  + \underbrace{\sum_{i,j,k=1}^d \norm{\Dr{j}^*\left(  \coa{\S_{\e_k,ij}} \star \Phi_{\ec^{-1}}  \right) \left(\Upsilon \Dr{k}\vb\right)}^2_{L^2(\cu_m)}}_{\text{\cref{eq:QuantFDecom}-d}} \\
& \qquad \left. + \underbrace{\sum_{i,k=1}^d  \norm{\ac\Dr{i}\left( \coa{\phi_{\e_k}}_{\Pcu}^{\eta}  \star  \Phi_{\ec^{-1}}\right) \left(\Upsilon\Dr{k}\vb\right)}^2_{L^2(\cu_m)}}_{\text{\cref{eq:QuantFDecom}-e}}\right).
\end{split}
\end{equation}
We treat them term by term. For $\cref{eq:QuantFDecom}$-a, noticing that $(1 - \Upsilon) \leq \Ind{\dist(\cdot, \partial) \leq 2\ell(\ec)}$, we apply the trace formula \cref{eq:Trace}
\begin{equation}\label{eq:F1}
\begin{split}
\text{\cref{eq:QuantFDecom}-a} & = \sum_{i=1}^d\norm{(1 - \Upsilon)\left(\ac - \ab \right)\left(\Dr{i}\vb\right)}^2_{L^2(\cu_m)}\\
& \leq 2 \sum_{i=1}^d\norm{\left(\Dr{i}\vb\right)\Ind{\dist(\cdot, \partial) \leq 2\ell(\ec)}}^2_{L^2(\cu_m)}\\
& \leq C(d)\ell(\ec) \left(3^{-m} \norm{\D \vb}^2_{L^2(\cu_m)} + \norm{\D \vb}_{L^2(\cu_m)}\norm{\D^* \D \vb}_{L^2(\itr(\cu_m))} \right).
\end{split}
\end{equation}

For the term \cref{eq:QuantFDecom}-b, we notice that 
\begin{align*}
\Dr{i}\left(\Upsilon \Dr{k}\vb\right) = (\Dr{i}\Upsilon)(\Dr{k}\vb) + \Upsilon(\cdot + \e_i) (\Dr{i}\Dr{k}\vb),
\end{align*}
and the support of $\Dr{i}\Upsilon$ is contained in the region of distance between $\ell(\ec)$ and $2\ell(\ec)$ from $\partial \cu_m$ i.e.
$$
\Dr{i}\Upsilon \leq \frac{1}{\ell(\ec)}\Ind{\cdot \in \cu_m, \frac{1}{2}\ell(\ec) \leq \dist(\cdot, \partial) \leq 3 \ell(\ec)},
$$
then we apply these in \cref{eq:QuantFDecom}-b and also \cref{eq:Trace} and obtain that
\begin{equation}\label{eq:F2}
\begin{split}
\text{\cref{eq:QuantFDecom}-b} &= \sum_{i,k=1}^d \norm{\phim_{\e_k}(\cdot + \e_k) \ac\Dr{i}\left(\Upsilon \Dr{k}\vb\right)}^2_{L^2(\cu_m)} \\
&\leq \sum_{i,k=1}^d \norm{\underbrace{\phim_{\e_k}(\cdot + \e_k)}_{\leq \Y_1} (\Dr{i}\Upsilon)(\Dr{k}\vb)}^2_{L^2(\cu_m)} + \sum_{i,k=1}^d \norm{ \underbrace{\phim_{\e_k}(\cdot + \e_k)}_{\leq \Y_1}  \Upsilon(\cdot + \e_i) (\Dr{i}\Dr{k}\vb)}^2_{L^2(\cu_m)}\\
&\leq \sum_{i,k=1}^d \Y_1^2 \norm{  \frac{1}{\ell(\ec)}\Ind{\frac{1}{2}\ell(\ec) \leq \dist(\cdot, \partial) \leq 3 \ell(\ec)} (\Dr{k}\vb)}^2_{L^2(\cu_m)} + \Y_1^2 \sum_{i,k=1}^d \norm{  \Ind{\dist(\cdot, \partial \cu_m) \geq \ell(\ec)} (\Dr{i}\Dr{k}\vb)}^2_{L^2(\cu_m)}\\
&\leq C(d)\Y_1^2\left( \frac{1}{\ell(\ec)}\left(3^{-m} \norm{\D \vb}^2_{L^2(\cu_m)} + \norm{\D \vb}_{L^2(\cu_m)}\norm{\D^* \D \vb}_{L^2(\itr(\cu_m))} \right) + \norm{\D^* \D \vb}^2_{L^2(\itr(\cu_m))}  \right).
\end{split}
\end{equation}
In the last step, we apply \cref{eq:Trace} and use the interior $H^2$ norm of $\vb$ since the function $\Upsilon$ is supported just in the interior with distance $\ell(\ec)$ from $\partial \cu_m$. \cref{eq:QuantFDecom}-c follows the similar estimate.

For the term \cref{eq:QuantFDecom}-d, we use the quantity $\Y_2$  to estimate it 
\begin{equation}\label{eq:F4}
\begin{split}
\text{\cref{eq:QuantFDecom}-d} &= \sum_{i,j,k=1}^d \norm{ \underbrace{\Dr{j}^*\left(  \coa{\S_{\e_k,ij}} \star \Phi_{\ec^{-1}}  \right)}_{\leq \Y_2} \left(\Upsilon \Dr{k}\vb\right)}^2_{L^2(\cu_m)} \\
& \leq C(d) \Y_2^2 \norm{\D \vb}^2_{L^2(\cu_m)}.
\end{split}
\end{equation}
The term \cref{eq:QuantFDecom}-e follows the similar estimate.

We combine \cref{eq:F1},\cref{eq:F2} and \cref{eq:F4} together and obtain that 
\begin{equation}\label{eq:QuantF}
\begin{split}
\norm{\Fi}_{L^2(\cu_m)} &\leq C(d)\left(\norm{\D \vb}_{L^2(\cu_m)} \left(3^{-\frac{m}{2}}\ell^{\frac{1}{2}}(\ec) + 3^{-\frac{m}{2}} \ell^{-\frac{1}{2}}(\ec)\Y_1 + \Y_2 \right)\right.\\
& \qquad \left.\norm{\D \vb}^{\frac{1}{2}}_{L^2(\cu_m)} \norm{\D^* \D \vb}^{\frac{1}{2}}_{L^2(\itr(\cu_m))}\left(\ell^{\frac{1}{2}}(\ec) + \ell^{-\frac{1}{2}}(\ec)\Y_1 \right)  + \norm{\D^* \D \vb}_{L^2(\itr(\cu_m))} \Y_1 \right).
\end{split}
\end{equation}

We put these two estimates \cref{eq:QuantWV} and \cref{eq:QuantFDecom} into \cref{eq:TestInOut4} and get the desired result.
\end{proof}

\subsection{Construction of a vector field}
In this part we calculate the vector field $\Fi$ used in the last paragraph.
\begin{lemma}\label{lem:defF}
There exists a vector field $\Fi : \Zd \rightarrow \Rd$ such that 
$$
\D^* \cdot \left(\ac \D w - \ab \D \vb \right) = \D^* \cdot \Fi,
$$
with the formula

\begin{equation}\label{eq:defF}
\begin{split}
\Fi_i &= (1 - \Upsilon)\left(\ac - \ab \right)\left(\Dr{i}\vb\right) 
+ \sum_{k=1}^d \phim_{\e_k}(\cdot + \e_k) \ac\Dr{i}\left(\Upsilon \Dr{k}\vb\right) \\
& \qquad - \sum_{j,k=1}^d \Sm_{\e_k,ij}(\cdot-\e_j)\Dr{j}^*\left(\Upsilon\Dr{k}\vb\right) + \sum_{j,k=1}^d \Dr{j}^*\left(  \coa{\S_{\e_k,ij}} \star \Phi_{\ec^{-1}}  \right) \left(\Upsilon \Dr{k}\vb\right) \\
& \qquad - \sum_{k=1}^d  \ac\Dr{i}\left( \coa{\phi_{\e_k}}_{\Pcu}^{\eta}  \star  \Phi_{\ec^{-1}}\right) \left(\Upsilon\Dr{k}\vb\right).
\end{split}
\end{equation}
\end{lemma}
\begin{proof}
We write
\begin{equation}\label{eq:FStep1}
\begin{split}
\left[\ac\D w - \ab \D \vb\right]_i(x) &= \left[(\ac - \ab)\D \vb + \sum_{k=1}^d \ac \underbrace{\D\left( \left(\Upsilon \Dr{k}\vb \right) \phim_{\e_k}\right)}_{\text{Using } \cref{eq:DiffProduct}}\right](x) \\
& = \underbrace{\left[(1-\Upsilon)(\ac - \ab)\D \vb\right]_i(x)}_{\text{\cref{eq:FStep1}-a}} + \underbrace{\sum_{k=1}^d \phim_{\e_k}(x+\e_i)\ac(x,x+\e_i) \Dr{i}\left(\Upsilon \Dr{k}\vb \right)(x)}_{\text{\cref{eq:FStep1}-b}} \\
& \qquad + \sum_{k=1}^d \left[\left(\ac \D\phim_{\e_k} + (\ac - \ab)\D l_{\e_k} \right) \left(\Upsilon \Dr{k}\vb \right) \right]_i(x)
\end{split}
\end{equation}
The terms \cref{eq:FStep1}-a and \cref{eq:FStep1}-b appear in the \cref{eq:defF} as the first and second term on the right hand side, so it suffices to treat the remaining terms in \cref{eq:FStep1}, where we apply the definition of $\S_{\e_k}$ \cref{eq:defS}
\begin{equation}\label{eq:FStep2}
\begin{split}
& \sum_{k=1}^d \left[\left(\ac \D\phim_{\e_k} + (\ac - \ab)\D l_{\e_k} \right) \left(\Upsilon \Dr{k}\vb \right) \right]_i(x) \\
= & \sum_{k=1}^d \left[\underbrace{\left(\ac \left(\D\phi_{\e_k} + \D l_{\e_k}\right)  - \ab\D l_{\e_k} \right)}_{= \D^* \cdot \S_{\e_k}} \left(\Upsilon \Dr{k}\vb \right) \right]_i(x)  - \sum_{k=1}^d \left[\ac \left(\D \coa{\phi_{\e_k}}_{\Pcu}^{\eta} \star \Phi_{\ec^{-1}} \right)\left(\Upsilon \Dr{k}\vb \right)\right]_i(x) \\
= & \underbrace{\sum_{k=1}^d \left[\D^* \cdot \Sm_{\e_k} \left(\Upsilon \Dr{k}\vb \right) \right]_i(x)}_{\text{\cref{eq:FStep2}-a}}  - \underbrace{\sum_{k=1}^d \left[\D^* \cdot \left(\coa{\S_{\e_k}} \star \Phi_{\ec^{-1}} \right) \left(\Upsilon \Dr{k}\vb \right) \right]_i(x)}_{\text{\cref{eq:FStep2}-b}}\\
& \qquad - \underbrace{\sum_{k=1}^d \left[\ac \left(\D \coa{\phi_{\e_k}}_{\Pcu}^{\eta} \star \Phi_{\ec^{-1}} \right)\left(\Upsilon \Dr{k}\vb \right)\right]_i(x)}_{\text{\cref{eq:FStep2}-c}}.
\end{split}
\end{equation}
The terms \cref{eq:FStep2}-b and \cref{eq:FStep2}-c also appear in the definition of $\Fi_i$ \cref{eq:defF} as the forth and fifth term. We study the term \cref{eq:FStep2}-a and use the anti-symmetry that $\S_{\e_k,ij} = - \S_{\e_k,ji}$ 
\begin{align*}
\D^* \cdot \text{\cref{eq:FStep2}-a} &= \sum_{i,k=1}^{d} \Dr{i}^* \left[\D^* \cdot \Sm_{\e_k} \left(\Upsilon \Dr{k}\vb \right) \right]_i(x) \\
&= \sum_{i,j,k=1}^{d} \Dr{i}^*\left( \left(\Dr{j}^*\Sm_{\e_k,ij}\right) \left(\Upsilon \Dr{k}\vb \right)\right)(x) \\
&= \underbrace{\sum_{i,j,k=1}^{d} \Dr{i}^* \Dr{j}^* \left( \Sm_{\e_k,ij} \left(\Upsilon \Dr{k}\vb\right)\right)(x)}_{ = 0 \text{ by anti-symmetry}} -  \sum_{i,j,k=1}^{d} \Dr{i}^*\left(\Sm_{\e_k,ij}(\cdot - \e_j)  \Dr{j}^*\left(\Upsilon \Dr{k}\vb \right)\right)(x) \\
&=  \sum_{i=1}^d \Dr{i}^*\left(- \sum_{j,k=1}^{d} \left(\Sm_{\e_k,ij}(\cdot - \e_j)  \Dr{j}^*\left(\Upsilon \Dr{k}\vb \right)\right)\right)(x)
\end{align*}
This gives the formula in \cref{eq:defF}.
\end{proof}

\section{Analysis of the algorithm}\label{sec:Analysis}
In this part, we give a rigorous analysis on the performance of our algorithm. At first, we prove that we can implement the algorithm on the whole domain $\cu_m$ instead of the cluster with the help of the mask operation.
\begin{proposition}[Arbitary extension]\label{prop:ArbitaryExtension}
After an arbitrary extension of the function $u_0, u_1, u_2$ defined in \cref{eq:iterativeCluster} on $\itr(\cu_m) \backslash \cltm$, the functions $u_1, \ub, u_2$ also satisfy
\begin{align}
\label{eq:iterativeClusterEquivalent}
\left\{
	\begin{array}{lll}
	(\ec_{\clt,m}^2 - \nabla \cdot \acf{m} \nabla)u_1 &= f_{\clt,m} + \nabla \cdot \acf{m} \nabla u_0  &  \text{ in } \itr(\cu_m) , \\
	-\nabla \cdot \ab \nabla \ub &= \ec_{\clt,m}^2 u_1  & \text{ in } \itr(\cu_m), \\
	(\ec_{\clt,m}^2 - \nabla \cdot \acf{m} \nabla) u_2 &= (\ec_{\clt, m}^2 - \nabla \cdot \ab \nabla) \ub &  \text{ in } \itr(\cu_m). \\	
	\end{array}
\right.
\end{align}

\end{proposition}
\begin{proof}
In the first equation of \cref{eq:iterativeClusterEquivalent} the left hand side can be rewritten as
$$
(\ec_{\clt,m}^2 - \nabla \cdot \acf{m} \nabla)u_1(x) = \ec_{\clt,m}^2(x) u_1(x) + \sum_{y \sim x}(\acf{m}(x,y))(u_1(x) - u_1(y)),
$$
while the right hand side equals 
$$
f_{\clt,m}(x) + \nabla \cdot \acf{m} \nabla u_0(x) =  f_{\clt,m}(x) + \sum_{y \sim x}(\acf{m}(x,y))(u_0(y) - u_0(x)).
$$
If $x \in \cltm \backslash \partial \cu_m$ the left hand side and the right hand side both equal to the first equation in \cref{eq:iterativeCluster}, so the equation is established. If $x \in \itr(\cu_m) \backslash \cltm$, no matter what values $u_1,u_0$ takes on the extension, the factors and function $f_{\clt,m}(x) = \ec_{\clt,m}(x)=\a_{\clt,m}(x,y)=0$ make both left hand side and right hand side $0$.

In the second equation, on the right hand side $\ec_{\clt,m}^2 u_1$ coincides with that in \cref{eq:iterativeCluster} so the equation is also established.

The third equation is valid, if $x \in \cltm \backslash \partial \cu_m$ for the similar reason as described in the first equation. If $x \in \itr(\cu_m) \backslash \cltm$, the left hand side equals $0$ since all the factors and conductance are $0$. The right hand side is also $0$ thanks to a simple manipulation using the second equation
$$
(\ec_{\clt,m}^2 - \nabla \cdot \ab \nabla) \ub(x) = \ec_{\clt,m}^2(x) (\ub(x)+u_1(x)) = 0,
$$
and this finishes the proof.
\end{proof}

The same idea also works for $u$ defined in \cref{eq:main}, which can also be defined as the solution 
\begin{equation}
\label{eq:mainEquivalent}
\left\{
	\begin{array}{ll}
	-\nabla \cdot \acf{m} \nabla u = f  & \qquad \text{ in } \itr(\cu_m), \\
	u = g & \qquad \text{ on } \cltm \cap \partial \cu_m,
	\end{array}
\right.
\end{equation}
with an arbitrary extension outside $\clt_*(\cu_m)$. We are now ready to complete the proof of \Cref{thm:main}, and we start by analyzing our algorithm with standard $H^1$ and $H^2$ estimates for $\ub$ in \cref{eq:iterativeCluster}.

\begin{lemma}[$H^1$ and $H^2$ estimates]\label{lem:H2}
In the iteration \cref{eq:iterativeCluster} we have the following estimates 
\begin{align}
\norm{\nabla \ub}_{L^2(\cu_m)} &\leq |\ab|^{-1}(1+\Lambda)\norm{\nabla(u-u_0)\Ind{\a \neq 0}}_{L^2(\cltm)}, \label{eq:H2Step1} \\
\norm{\D^* \D \ub}_{L^2(\itr(\cu_m))} &\leq C(d,\Lambda) |\ab|^{-1}\ec \norm{\nabla(u-u_0)\Ind{\a \neq 0}}_{L^2(\cltm)}, \label{eq:H2Step2} \\
\norm{\nabla (\hat{u} - u) \Ind{\a \neq 0}}_{L^2(\cltm)} &\leq 2|\ab|^{-1}(1+\Lambda)^2\norm{\nabla(u-u_0)\Ind{\a \neq 0}}_{L^2(\cltm)}. \label{eq:H1Trivial}
\end{align}
\end{lemma}
\begin{proof}
We start by testing \cref{eq:iterativeClusterEquivalent} and \cref{eq:mainEquivalent} with the function $u_1$, and we also use the trick that $\ec_{\clt,m}$ and $\acf{m}$ restrict the problem on $(\cltm, \Eda(\cltm))$ 
\begin{align*}
\bracket{\ec_{\clt,m} u_1, \ec_{\clt,m} u_1}_{\cu_m} + & \bracket{\nabla u_1, \acf{m} \nabla u_1}_{\cu_m} 
= \bracket{\nabla u_1, \acf{m} (u - u_0)}_{\cu_m} ,\\
\Longrightarrow \lambda^2 \norm{u_1}^2_{L^2(\cltm)} + & \Lambda^{-1}\norm{\nabla u_1 \Ind{\a \neq 0}}^2_{L^2(\cltm)} \\
&\leq\norm{\nabla(u-u_0)\Ind{\a \neq 0}}_{L^2(\cltm)}\norm{\nabla u_1 \Ind{\a \neq 0}}_{L^2(\cltm)}.
\end{align*}
We obtain that 
\begin{align}
\ec \norm{u_1}_{L^2(\cltm)} &\leq \Lambda \norm{\nabla(u-u_0)\Ind{\a \neq 0}}_{L^2(\cltm)}, \label{eq:H2Pre1}\\ 
\norm{\nabla u_1 \Ind{\a \neq 0}}_{L^2(\cltm)} &\leq \Lambda \norm{\nabla(u-u_0)\Ind{\a \neq 0}}_{L^2(\cltm)}. \label{eq:H2Pre2} 
\end{align}
Combining the first equation and the second equation in \cref{eq:iterativeClusterEquivalent} and \cref{eq:mainEquivalent}, we obtain that 
$$
-\nabla \cdot \ab \nabla  \ub = - \nabla \cdot \acf{m} \nabla (u - u_0 -u_1) \qquad \text{ in } \itr(\cu_m),
$$ 
then we test it by the function $\ub$ and use Cauchy's inequality to obtain that 
\begin{align*}
\bracket{\nabla \ub, \ab \nabla \ub}_{\cu_m} &= \bracket{\nabla \ub, \acf{m} \nabla(u-u_0-u_1)}_{L^2(\cu_m)} \\
& \leq \norm{\nabla \ub}_{L^2(\cu_m)} \norm{\nabla (u-u_0-u_1) \Ind{\a \neq 0}}_{\cltm}\\
\Longrightarrow \norm{\nabla \ub}_{L^2(\cu_m)} &\leq |\ab|^{-1}\norm{\nabla (u-u_0-u_1) \Ind{\a \neq 0}}_{\cltm}.
\end{align*}
Using \cref{eq:H2Pre2} we obtain that 
\begin{align*}
\norm{\nabla \ub}_{L^2(\cu_m)} &\leq |\ab|^{-1} \left(\norm{\nabla (u-u_0-u_1) \Ind{\a \neq 0}}_{\cltm}\right) \\
& \leq |\ab|^{-1} \left( \norm{\nabla (u-u_0) \Ind{\a \neq 0}}_{\cltm} + \norm{\nabla u_1 \Ind{\a \neq 0}}_{\cltm}\right) \\
& \leq |\ab|^{-1}(1+\Lambda) \norm{\nabla (u-u_0) \Ind{\a \neq 0}}_{\cltm}.
\end{align*}
This proves the formula \cref{eq:H2Step1}. 

Concerning \cref{eq:H2Step2}, we use the estimation of $H^2$ regularity \cref{eq:H2DisInter} for $-\nabla \cdot \ab \nabla \ub = \ec^2_{\clt,m}u_1$ since $\ab$ is constant and obtain that 
\begin{align*}
\sum_{i,j=1}^d \norm{\Dr{i}^* \Dr{j} \ub}^2_{L^2(\itr(\cu_m))} &\leq C(d)|\ab|^{-2} \norm{\ec^2_{\clt, m}u_1}^2_{L^2(\cu_m)} .
\end{align*}
We put the result from \cref{eq:H2Pre1} and \cref{eq:H2Step1} and obtain that 
$$
\norm{\D^* \D \ub}_{L^2(\itr(\cu_m))} \leq C(d,\Lambda)|\ab|^{-1} \ec \norm{\nabla (u-u_0) \Ind{\a \neq 0}}_{L^2(\cltm)}.
$$

To prove \cref{eq:H1Trivial}, we put \cref{eq:mainEquivalent}, the first equation and the second equation of \cref{eq:iterativeClusterEquivalent} into the right hand side of the third equation and obtain that 
$$
(\ec_{\clt,m}^2 - \nabla \cdot \acf{m} \nabla) u_2 = \ec_{\clt,m}^2 \ub^2 -\nabla \cdot \acf{m} \nabla (u - u_0 - u_1)  \qquad \text{ in } \itr(\cu_m).
$$
We subtract $(\ec_{\clt,m}^2 - \nabla \cdot \acf{m} \nabla) \ub$ on the two sides to obtain 
$$
(\ec_{\clt,m}^2 - \nabla \cdot \acf{m} \nabla)(u_2 - \ub) =  -\nabla \cdot \acf{m} \nabla (u - u_0 - u_1 - \ub)  \qquad \text{ in } \itr(\cu_m),
$$
and then we test it by $(u_2 - \ub)$ to obtain that 
\begin{equation}\label{eq:u2ub}
\norm{\nabla (u_2 - \ub) \Ind{\a \neq 0}}_{L^2(\cltm)} \leq \Lambda \norm{\nabla (u - u_0 - u_1 - \ub) \Ind{\a \neq 0}}_{L^2(\cltm)}.
\end{equation}
Therefore, combining \cref{eq:u2ub}, \cref{eq:H2Pre1} and \cref{eq:H2Step1} we can obtain a trivial bound for our algorithm
\begin{align*}
\norm{\nabla (\hat{u} - u) \Ind{\a \neq 0}}_{L^2(\cltm)} &\leq \norm{\nabla (u - u_0 - u_1 - \ub) \Ind{\a \neq 0}}_{L^2(\cltm)} +  \norm{\nabla (u_2 - \ub) \Ind{\a \neq 0}}_{L^2(\cltm)} \\
& \leq 2|\ab|^{-1}(1+\Lambda)^2 \norm{\nabla (u - u_0) \Ind{\a \neq 0}}_{L^2(\cltm)}.
\end{align*}
\end{proof}

The trivial bound \cref{eq:H1Trivial} is not optimal. In the typical case $\cu_m \in \Pcu_*$ in large scale, we can use \Cref{thm:TwoScale} to help us get a better bound, and this help use conclude the performance of our algorithm.
\begin{proof}[Proof of \Cref{thm:main}]
We analyze the algorithm in two cases: $\cu_m \in \Pcu_*$ and $\cu_m \notin \Pcu_*$. In the case $\cu_m \notin \Pcu_*$, we use \cref{eq:H1Trivial} that 
$$
\norm{\nabla (\hat{u} - u) \Ind{\a \neq 0}}_{L^2(\cltm)} \Ind{\cu_m \notin \Pcu_*} \leq 2|\ab|^{-1}(1+\Lambda)^2 \norm{\nabla (u - u_0) \Ind{\a \neq 0}}_{L^2(\cltm)}  \Ind{\cu_m \notin \Pcu_*}.
$$

In the case $\cu_m \in \Pcu_*$, we combine the first equation and the second equation of \cref{eq:iterativeClusterEquivalent} and \cref{eq:mainEquivalent}, together with the third term they give
\begin{align*}
	- \nabla \cdot \acf{m} \nabla(u - u_0 - u_1) &= -\nabla \cdot \ab \nabla \ub &  \text{ in } \itr(\cu_m), \\
	(\ec_{\clt,m}^2 - \nabla \cdot \acf{m} \nabla) u_2 &= (\ec_{\clt, m}^2 - \nabla \cdot \ab \nabla) \ub &  \text{ in } \itr(\cu_m). \\	
\end{align*}
This gives us two equations of two-scale expansion. As in \cref{eq:Expansion}, we define $w := \ub + \sum_{k=1}^d (\Upsilon \Dr{k}\ub) \phim_{\e_k}$ and apply \Cref{thm:TwoScale}
\begin{equation}
\begin{split}
&\norm{\nabla(\hat{u} - u)\Ind{\a \neq 0}}_{L^2(\cltm)}  \Ind{\cu_m \in \Pcu_*} \\
& \leq \left(\norm{\nabla(w - (u - u_0 - u_1))\Ind{\a \neq 0}}_{L^2(\cltm)} + \norm{\nabla(u_2 - w)\Ind{\a \neq 0}}_{L^2(\cltm)}\right) \Ind{\cu_m \in \Pcu_*}.
\end{split}
\end{equation}
The last equation gives a bound of type \cref{prop:TwoScaleDetail}. Together with the \Cref{lem:H2} and the estimate for case $\cu_m \notin \Pcu_*$, we obtain that
$$
\norm{\nabla(\hat{u} - u)\Ind{\a \neq 0}}_{L^2(\cltm)} \leq \Z \norm{\nabla(u_0 - u)\Ind{\a \neq 0}}_{L^2(\cltm)}
$$
where $\Z$ is given by 
\begin{equation}\label{eq:defZ}
\begin{split}
\Z = C(d,\Lambda)\left( 3^{-\frac{m}{2}} \ell^{-\frac{1}{2}}(\ec)\Y_1\X^{d} + \Y_2\X^d   + \ec^{\frac{1}{2}} \ell^{-\frac{1}{2}}(\ec)\Y_1\X^{d} + \ec \Y_1 \X^d + \Ind{\cu_m \notin \Pcu_*}\right).
\end{split}
\end{equation}
This gives the exact expression of the quantity $\Z$. To conclude, we have to quantify $\Z$ and we use \cref{eq:QuantX}, \cref{eq:QuantY}, \cref{eq:PcuScale} and \cref{eq:OSum} that there exist two positive constants $s(d,\pp,\Lambda)$ and $C(d,\pp,\Lambda,s)$ such that
\begin{equation}\label{eq:QuantZ}
\Z \leq \GO_s\left(C(d,\pp,\Lambda,s)\left( \left(3^{-\frac{m}{2}}+\ec+\ec^{\frac{d}{2}} + \ec^{\frac{1}{2}}\right) m^{\frac{1}{s}+d}\ell^{\frac{1}{2}}(\ec) + \ec m^{\frac{1}{s}}\ell(\ec) + 3^{-m}\right)\right).
\end{equation}
Observing that $3^{-m} < \ec$, then the dominating order writes $\Z \leq \GO_s\left(C \ec^{\frac{1}{2}} \ell^{\frac{1}{2}}(\ec)  m^{\frac{1}{s}+d} \right)$ and this concludes the proof of \Cref{thm:main}.
\end{proof}

\section{Numerical experiments}
\label{sec:Numeric}
We report on numerical experiments corresponding to our algorithm. In a cube $\cu$ of size $L$, we try to solve a localized corrector problem, that is, we look for the function $\phi_{L,p} \in C_0(\cu)$ such that 
\begin{equation}
	\begin{array}{ll}
	-\nabla \cdot \a \nabla (\phi_{L,p} + l_p) = 0  & \text{ in } \clt_*(\cu). \\
	\end{array}
\end{equation}
The quantity $\phi_{L,p}$ is very similar to the corrector $\phi_p$ and has sublinear growth. This is a good example for illustrating the usefulness of our algorithm, since the homogenized approximation to this function is simply the null function, which is not very informative.

\begin{figure}[h!]
\centering
\includegraphics[scale=0.6]{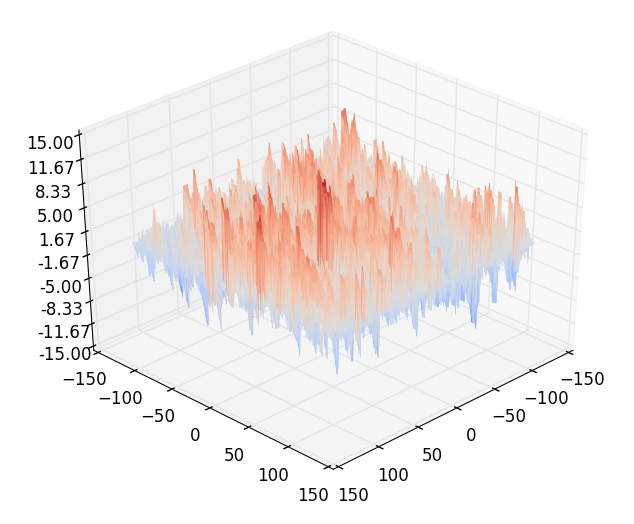}
\caption{A simulation for the corrector on the maximal cluster in a cube $243 \times 243$.}\label{Fig:Corrector}
\end{figure}

In our example, we take $d=2$, $p = \e_1$ and $L = 243$. We implement the algorithm to get a series of approximated solutions $\hat{u}_n$ where $\hat{u}_0 = 0$. Moreover, we use the residual error to see the convergence
$$
\text{res}(\hat{u}_n) := \frac{1}{|\cu|} \norm{-\nabla \cdot \a \nabla (\hat{u}_n + l_p)}_{L^2(\clt_*(\cu))} = \frac{1}{|\cu|} \norm{-\nabla \cdot \a \nabla (\hat{u}_n - \phi_{L,p})}_{L^2(\clt_*(\cu))}.
$$
See the \Cref{Fig:Corrector} for a simulation of the corrector $\phi_{L,p}$ with high resolution, and \Cref{Fig:Errors} for its residual errors.

\begin{figure}[h!]
\begin{adjustwidth}{}{}
\centering
\begin{tabular}{|c | c } 
 \hline
 round & errors  \\ [0.5ex] 
 \hline
	1 & 0.0282597982969 \\
	2 & 0.0126490361046 \\
	3 & 0.00707540548365 \\
	4 & 0.00435201077274 \\
	5 & 0.00282913420116 \\
	6 & 0.00190945842802 \\
	7 & 0.00132483912845 \\
	8 & 0.000939101476657 \\
 \hline
\end{tabular}
\end{adjustwidth}
\caption{A table of errors $\text{res}(\hat{u}_n)$.}
\label{Fig:Errors}
\end{figure}

\appendix
\section{Proof of some discrete functional inequality}\label{sec:Appendix}

\begin{lemma}[$H^2$ interior estimate for elliptic equation]\label{lem:H2DisInter}
Given two functions $v,f \in C_0(\cu_m)$ satisfying the discrete elliptic equation
\begin{equation}\label{eq:HarmonicAPP}
- \Delta v = f, \qquad \text{ in } \itr(\cu_m),
\end{equation}
we have an interior estimate
\begin{equation}\label{eq:H2DisInterAPP}
\norm{\D^* \D v}^2_{L^2(\itr(\cu_m))} := \sum_{i,j = 1}^d \norm{\Dr{i}^* \Dr{j} v}^2_{L^2(\itr(\cu_m))} \leq d \norm{f}^2_{L^2(\itr(\cu_m))}.
\end{equation}
\end{lemma}

\begin{proof}
We extend the elliptic equation to the whole space at first. The function $v,f$ have a natural null extension on $\Zd$ satisfying 
$$
\D^* \cdot \D v = f + (\D^* \cdot \D v) \Ind{\partial \cu_m} , \qquad \text{ in } \Zd.
$$
To simplify the notation, we denote by $\bar{f}$ the term on the right hand side. Then, by one step difference of direction $\e_j$, we have 
\begin{align*}
  \D^* \cdot \D (\Dr{j}v(x)) &= \Dr{j}\bar{f}(x). 
\end{align*}
We test this equation with a function $\phi$ of compact support, then by \cref{eq:IPP} we obtain
$$
\bracket{\D \phi, \D (\Dr{j} v) }_{\Zd} = \bracket{\Dr{j}^*\phi, \bar{f}}_{\Zd}.
$$
Putting $\phi = (\Dr{j} v)$ in this formula, we obtain that 
\begin{align*}
\bracket{\D (\Dr{j} v), \D (\Dr{j} v) }_{\Zd} &= \bracket{\Dr{j}^* \Dr{j} v, \bar{f}}_{\Zd} \\
&= \bracket{\Dr{j}^* \Dr{j} v, f}_{\Zd} + \bracket{\Dr{j}^* \Dr{j} v, (\D^* \cdot \D v) \Ind{\partial \cu_m}}_{\Zd}.  
\end{align*}
We do the sum over the $d$ canonical directions and get 
\begin{align*}
\sum_{i,j=1}^d \bracket{\Dr{i}\Dr{j} v, \Dr{i}\Dr{j} v }_{\Zd} &= \sum_{j=1}^d\bracket{\Dr{j}^* \Dr{j} v, f}_{\Zd} + \sum_{j=1}^d \bracket{\Dr{j}^* \Dr{j} v, (\D^* \cdot \D v) \Ind{\partial \cu_m}}_{\Zd} \\
&= \sum_{j=1}^d\bracket{\Dr{j}^* \Dr{j} v, f}_{\Zd} + \bracket{\D^* \cdot \D v , (\D^* \cdot \D v) \Ind{\partial \cu_m}}_{\Zd}.  
\end{align*}
Since $\Dr{j}^*v(x) = -\Dr{j}v(x-\e_j)$, we have 
\begin{align*}
\sum_{i,j=1}^d \bracket{\Dr{i}^*\Dr{j} v, \Dr{i}^*\Dr{j} v }_{\Zd} 
= \sum_{j=1}^d\bracket{\Dr{j}^* \Dr{j} v, f}_{\Zd} + \bracket{\D^* \cdot \D v , (\D^* \cdot \D v) \Ind{\partial \cu_m}}_{\Zd}.  
\end{align*}
There are three observations for this equation.
\begin{itemize}
\item $\supp(\Dr{i}^*\Dr{j} v) \subset \cu_m$.
\item For any $x \in \partial \cu_m$,
$$
(\D^* \cdot \D v )^2 = \left(\sum_{j=1}^d \Dr{j}^* \Dr{j} v(x)\right)^2 =  \sum_{j=1}^d \left(\Dr{j}^* \Dr{j} v(x)\right)^2,
$$
since $v=0$ on the boundary and only one term of $\left\{\Dr{j}^* \Dr{j} v(x)\right\}_{j=1 \cdots d}$ is not null.
\item On the boundary $\partial \cu_m, f = 0$ since $f\in C_0(\cu_m)$.
\end{itemize}
Combining the three observations, we get 
\begin{align*}
\sum_{i,j=1}^d \bracket{\Dr{i}^*\Dr{j} v, \Dr{i}^*\Dr{j} v }_{\cu_m} 
= \sum_{j=1}^d\bracket{\Dr{j}^* \Dr{j} v, f}_{\itr(\cu_m)} + \sum_{j=1}^d \bracket{\Dr{j}^* \Dr{j} v, \Dr{j}^* \Dr{j} v }_{\partial \cu_m}.  
\end{align*}
Thus, all the terms in the last sum on the right hand side can be found on the left hand side. We use Cauchy's inequality and Young's inequality 
\begin{align*}
\sum_{i,j=1}^d \bracket{\Dr{i}^*\Dr{j} v, \Dr{i}^*\Dr{j} v }_{\itr(\cu_m)} 
& \leq \sum_{j=1}^d\bracket{\Dr{j}^* \Dr{j} v, f}_{\itr(\cu_m)} \\
& \leq \sum_{j=1}^d \bracket{\Dr{j}^* \Dr{j} v, \Dr{j}^* \Dr{j} v}_{\itr(\cu_m)}^{\frac{1}{2}}  \bracket{f, f}_{\itr(\cu_m)}^{\frac{1}{2}} \\
& \leq \sum_{j=1}^d \left( \frac{1}{2}\bracket{\Dr{j}^* \Dr{j} v, \Dr{j}^* \Dr{j} v}_{\itr(\cu_m)} + \frac{1}{2} \bracket{f, f}_{\itr(\cu_m)} \right) \\
\Longrightarrow \sum_{i,j=1}^d \bracket{\Dr{i}^*\Dr{j} v, \Dr{i}^*\Dr{j} v }_{\itr(\cu_m)} & \leq d\bracket{f, f}_{\itr(\cu_m)},
\end{align*}
which concludes the proof.
\end{proof}

\begin{figure}[h]
  	\centering
    \includegraphics[scale = 0.5]{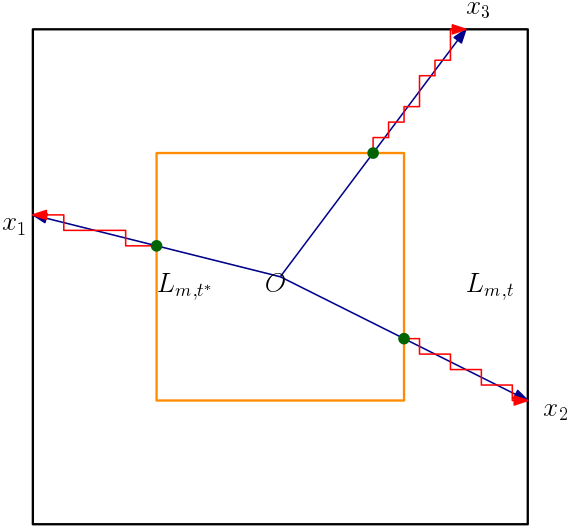}
    \caption{To construct the path $\gamma^{x,t^*}$ for every $x \in L_{m,t}$, one can find at first the pivot level $L_{m,t^*}$. Then we connect $O$ and $x$ and find one of its closest discrete path in $(\Zd, \Ed)$ and denote by $\gamma^{x, t^*}$ the segment from $L_{m,t^*}$ to $x$. By this consturction, every edge and vertex is passed by the paths $\{\gamma^{x,t^*}\}_{x\in L_{m, t}}$ at most $C(d)$ times. In this picture, the arrows in blud indicate the vectors $Ox_1, Ox_2, Ox_3$ and the segments with arrow in red are the paths $\gamma^{x_1, t^*}, \gamma^{x_2, t^*}, \gamma^{x_3, t^*}$.}\label{Fig:Path2.png}
\end{figure}

The same technique to do an integration along the path helps us to get an estimate of trace.
\begin{lemma}[Trace inequality]\label{lem:Trace}
For every $ u : \cu_m \rightarrow \mathbb{R}$ and $0 \leq K \leq \frac{3^m}{4}$, we have the following inequality 
\begin{equation}\label{eq:TraceAPP}
\norm{u \Ind{\dist(\cdot, \partial \cu_m) \leq K}}^2_{L^2(\cu_m)} 
\leq C(d) (K+1) \left(3^{-m}\norm{u}^2_{L^2(\cu_m)} + \norm{u}_{L^2(\cu_m)} \norm{\nabla u}_{L^2(\cu_m)}\right).
\end{equation}
\end{lemma}
\begin{proof}
We use the notation $L_{m,t}$ to define the level set in $\cu_m$ with distance $t$ to the boundary 
$$
L_{m,t} := \left\{x \in \cu_m : \dist(x, \partial \cu_m) = t \right\}.
$$
Then, we observe that $L_{m,0} = \partial \cu_m$ and we have the partition 
$$
\cu_m = \bigsqcup_{t = 0}^{\left\lfloor \frac{3^m}{2} \right\rfloor} L_{m,t}.
$$

Using the pigeonhole principle, it is easy to prove that there exists a $t^* \in \left[0, \left\lfloor \frac{3^m}{4} \right\rfloor -1\right]$ such that 
\begin{equation}\label{eq:pivot}
\norm{u}^2_{L^2(L_{m,t^*})} \leq \frac{4}{3^m}\norm{u}^2_{L^2(\cu_m)},
\end{equation}
and we define $t^* := \arg\min_{t \in \left[0, \left\lfloor \frac{3^m}{4} \right\rfloor \right]} \norm{u}^2_{L^2(L_{m,t})}$. We call $L_{m,t^*}$ the \textit{pivot level} and it plays the same role as the null boundary in the proof of Poincar\'e's inequality. In the following, we will apply the trick of integration along the path to prove the \cref{eq:Trace} for one lever $L_{m,t}$. For every $x \in L_{m,t}$, we denote by $\rt(x,t^*)$ a \textit{root} on the pivot level $L_{m,t^*}$, and choose a path $\gamma_{x, t^*} = \{\gamma^{x, t^*}_k\}_{0 \leq k \leq n}$ such that 
$$
\gamma^{x, t^*}_0 = \rt(x,t^*),\quad \gamma^{x, t^*}_k \sim \gamma^{x, t^*}_{k+1}, \quad \gamma^{x, t^*}_n=x.   
$$
Moreover, we use $|\gamma^{x, t^*}|$ to represent the number of steps of the path, for example here $|\gamma^{x, t^*}| = n$. We apply a discrete Newton-Leibniz formula to get
\begin{align*}
u^2(x) - u^2(\rt(x,t^*)) = \sum_{k=0}^{|\gamma^{x, t^*}|} \left(u^2(\gamma^{x, t^*}_{k+1}) - u^2(\gamma^{x, t^*}_k)\right) =\sum_{k=0}^{|\gamma^{x, t^*}|} \nabla u(\gamma^{x, t^*}_k, \gamma^{x, t^*}_{k+1}) \left(u(\gamma^{x, t^*}_{k+1}) + u(\gamma^{x, t^*}_k)\right).
\end{align*}
We put this formula into the norm of $\norm{u}_{L^2(L_{m,t})}, t  \in \left[0, \left\lfloor \frac{3^m}{4} \right\rfloor \right] $ and and apply Cauchy's inequality to obtain that 
\begin{align*}
\norm{u}^2_{L^2(L_{m,t})} &= \sum_{x \in L_{m,t}} \left(u^2(\rt(x,t^*)) +  \sum_{k=0}^{|\gamma^{x, t^*}|} \nabla u(\gamma^{x, t^*}_k, \gamma^{x, t^*}_{k+1}) \left(u(\gamma^{x, t^*}_{k+1}) + u(\gamma^{x, t^*}_k)\right)\right)\\
&\leq \sum_{x \in L_{m,t}} u^2(\rt(x,t^*)) \\
& \qquad + 2\left(\sum_{x \in L_{m,t}} \sum_{k=0}^{|\gamma^{x, t^*}|} \left(\nabla u(\gamma^{x, t^*}_k, \gamma^{x, t^*}_{k+1})\right)^2\right)^{\frac{1}{2}} \left( \sum_{x \in L_{m,t}} \sum_{k=0}^{|\gamma^{x, t^*}|} \left(u^2(\gamma^{x, t^*}_{k+1}) + u^2(\gamma^{x, t^*}_k)\right)\right)^{\frac{1}{2}}\\
&\leq \sum_{y \in L_{m,t^*}} u^2(y)\left(\sum_{x \in L_{m,t}} \Ind{y = \rt(x, t^*)}\right) \\
& \qquad + 4\left(\sum_{\{y_1,y_2\} \in \Ed(\cu_m)} \left(\nabla u(y_1, y_2)\right)^2 \left(\sum_{x \in L_{m,t}} \Ind{\{y_1,y_2\}\in \gamma^{x, t^*}}\right)\right)^{\frac{1}{2}} \\
& \qquad \qquad \times \left(\sum_{y \in \cu_m}u^2(y) \left(\sum_{x \in L_{m,t}} \Ind{y \in \gamma^{x, t^*}}\right) \right)^{\frac{1}{2}}.
\end{align*}
The next step is to decide how to choose the root $\rt(x,t^*)$ and the path. The main idea is to make every edge and every vertex as root is passed by $\{\gamma^{x, t^*}\}_{x \in L_{m, t}}$ a finite number of times bounded by a constant $C(d)$. One possible plan is to choose the root $\rt(x,t^*)$ and the path $\gamma^{x, t^*}$ a discrete path in $(\Zd, \Ed)$ which is the closest to the vector $\overrightarrow{Ox}$, then it is a simple exercise to see that it gives us a bound $C(d)$. See \Cref{Fig:Path2.png} as a visualization. Then we get that 
\begin{align*}
\norm{u}^2_{L^2(L_{m,t})} \leq C(d) \left(\norm{u}^2_{L^2(L_{m,t*})} +  \norm{u}_{L^2(\cu_m)} \norm{\nabla u}_{L^2(\cu_m)}\right).
\end{align*}
Then we put the \cref{eq:pivot} and get 
\begin{align*}
\norm{u}^2_{L^2(L_{m,t})} \leq 4C(d) \left(3^{-m}\norm{u}^2_{L^2(\cu_m)} +  \norm{u}_{L^2(\cu_m)} \norm{\nabla u}_{L^2(\cu_m)}\right).
\end{align*}
\cref{eq:Trace} is just a result by summing all the levels of distance less than $K$.
\end{proof}

\section{Small clusters}\label{sec:SmallCluster}
This part is devoted to studying the small clusters in the percolation. Many of the arguments presented here have appeared in the previous work \cite{armstrong2018elliptic}. We extract those results from \cite{armstrong2018elliptic} and expand upon certain points that are useful for our purposes. The motivation to state these results comes from the technique of partition of good cubes:
\begin{question}\label{Q:SmallCluster1}
In a cube $\cu \in \Tcu$ and its enlarged domain $\clp(\cu)$, besides the maximal cluster $\clt_*(\cu)$, what is the behavior of the other finite connected clusters ?
\end{question}
\begin{question}\label{Q:SmallCluster2}
When we apply \Cref{lem:FusionCluster}, since $\clt_*(\cu)$ and $\bigcup_{z \in \cu} \clt_*(\cu_{\Pcu}(z))$ are not necessarily equal, how can we describe the difference between the two ?
\end{question}
\begin{question}\label{Q:SmallCluster3}
What is the difference between $\clt_{\infty} \cap \cu$ and $\clt_*(\cu)$ ?
\end{question}

We start with a first very elementary lemma:
\begin{lemma}\label{lem:SmallClusterElementary}
For any $\cu \in \Tcu$ and $ z \in (\clt_{\infty} \cap \clp(\cu)) \backslash \clt_*(\cu)$, there exists a cluster $\clt'$ such that $z \in \clt'$ and $\clt' \xleftrightarrow{\a} \partial \clp(\cu)$.
\end{lemma}
\begin{proof}
For a cube $\cu \in \Tcu$ and its enlarged domain $\clp(\cu)$, there exist three types of clusters: 
\begin{enumerate}
\item One unique maximal cluster $\clt_*(\cu)$;
\item The isolated clusters which connect neither to $\clt_*(\cu)$ nor to the boundary $\partial \clp(\cu)$;
\item The clusters which do not connect to $\clt_*(\cu)$ but connect to the boundary.
\end{enumerate}
Then it is clear the cluster $\clt'$ containing $z \in (\clt_{\infty} \cap \clp(\cu)) \backslash \clt_*(\cu)$ can only be of the third type and this proves the lemma.
\end{proof} 

We define the third class above as small cluster. For any $z \in \Zd$, we denote by $\clt'(z)$ the clusters containing $z$.
\begin{definition}[Small clusters]\label{def:SmallCluster}
For any $\cu \in \Tcu$, we define \textit{small clusters} in $\cu$ as the union of clusters , restricted to $\clp(\cu)$, different from $\clt_*(\cu)$ but connecting to $\partial \clp(\cu)$, and we denote it by $\slt(\cu)$ i.e.
$$
\slt(\cu) := \bigcup_{z \in \partial \clp(\cu) \backslash \clt_*(\cu)} \clt'(z).
$$
\end{definition}

Intuitively, these small clusters should be of order $\size(\cu)^{d-1}$ when the cube $\cu$ is large. This is indeed true, as we prove the following lemma:
\begin{lemma}\label{lem:SmallCluster}
For any $\cu \in \Tcu$, the set $\slt(\cu)$ has the following decomposition
\begin{equation}\label{eq:SmallClusterDecomposition}
\slt(\cu) \subset \bigcup_{z \in \partial \clp(\cu)} \cu_{\Pcu}(z),
\end{equation}
and has the estimate
\begin{equation}\label{eq:SmallClusterSize}
|\slt(\cu)| \Ind{\cu \in \Pcu_*} \leq \GO_1(C\size(\cu)^{d-1}).
\end{equation}
\end{lemma}
\begin{proof}
We prove at first \cref{eq:SmallClusterDecomposition}. In the case that $\cu \notin \Pcu_*$, it is obvious since it has to enlarge to $\clp(\cu)$ which is a larger cube, and all the terms on the right hand side of \cref{eq:SmallClusterDecomposition} refer to $\clp(\cu)$.

In the case that $\cu \in \Pcu_*$, we consider one cluster $\clt'$ connecting to $x \in \partial \cu$. We suppose that it is not contained in the union of the elements of $\Pcu$ lying on $\partial \cu$, then it has to cross into the interior. As illustration in \Cref{Fig:SmallCases}, it has several situations: 
\begin{figure}[h]
\includegraphics[scale=0.5]{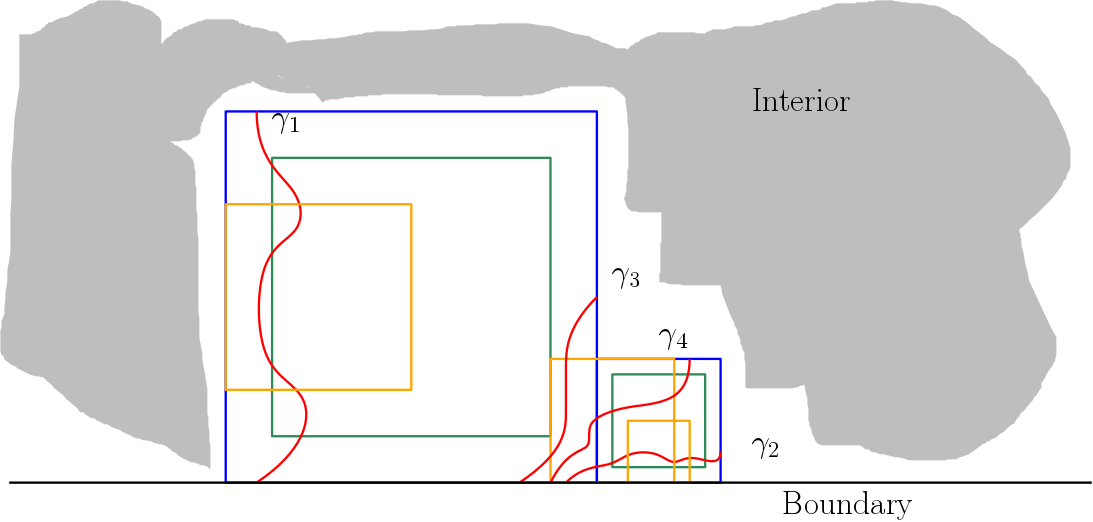}
\centering
\caption{The image explain why $\clt'$ is contained in the union of partition cubes lying at the $\partial \cu$. Without loss of generality, we suppose the big blue cube is $\cu_{\Pcu}(x)$ and the small one is its neighbor good cube. The cube in color of green is the part of size $\frac{3}{4}$ of the good cube. The paths in red represent different typical situations that if a finite cluster connects to the boundary $\partial \cu$ and wants to cross $\bigcup_{z \in \partial \clp(\cu)} \cu_{\Pcu}(z)$.}\label{Fig:SmallCases}
\end{figure}
\begin{enumerate}
\item The first case is that $\clt'$ cross at least one pair of $(d-1)$-dimensional opposite face of partition cube, as showed as $\gamma_1$ or  $\gamma_2$, so we have $|\diam(\clt')| > \size(\cu_{\Pcu}(x))$ for the case $\gamma_1$( or $|\diam(\clt')| > \frac{1}{3}\size(\cu_{\Pcu}(x))$ for the case $\gamma_2$). Then by the definition of partition cube, we can find a cube $\cu'$ of  $\frac{1}{2}\size(\cu_{\Pcu}(x))$ to contain parts of $\gamma_1$ and $\cu'$ intersects $\frac{3}{4}\cu_{\Pcu}(x)$, so by the definition of good cube we have necessarily $\clt' \xleftrightarrow{\a} \clt_*(\cu_{\Pcu}(x))$. Same discussion can be also applied to the case $\gamma_2.$ This gives a contradiction.

\item The second case is that $\clt'$ does not cross any pair of $(d-1)$-dimensional opposite face of partition cube, but also enter the interior of $\cu$ by $\partial \cu_{\Pcu}(x)$ or the boundary of its neighbor, so $|\diam(\clt')| > \frac{1}{3}\size(\cu_{\Pcu}(x))$. One can always find a cube $\cu'$ of size $\frac{1}{3}\size(\cu_{\Pcu}(x))$ crossed by $\clt'$. If it is the case in $\gamma_3$ that $\cu'$ intersects $\frac{3}{4}\cu_{\Pcu}(x)$, then we apply the definition of good cubes and $\clt' \xleftrightarrow{\a} \clt_*(\cu_{\Pcu}(x))$. Otherwise, in the case $\gamma_4$, $\clt'$ must cross a cube $\cu''$ of size $\frac{1}{6} \size(\cu_{\Pcu}(x)$ in its neighbor and we apply the same discussion, which also gives a contradiction.
\end{enumerate}

To estimate the upper bound \cref{eq:SmallClusterSize}, we use the decomposition above and calculate the volume of $\bigcup_{z \in \partial \cu} \cu_{\Pcu}(z) $ by doing a contour integration along $\partial \cu$ of height function $\size(\cu_{\Pcu}(z))$ and then applying \cref{eq:OSum},  
\begin{align*}
 |\slt(\cu)| \Ind{\cu \in \Pcu_*} \leq  \left| \bigcup_{z \in \partial \cu} \cu_{\Pcu}(z) 
\right|  \leq  \sum_{z \in \partial \cu} \size(\cu_{\Pcu}(z)) \leq \GO_1(C\size(\cu)^{d-1}).
\end{align*}

\begin{figure}[h]
\includegraphics[scale=0.7]{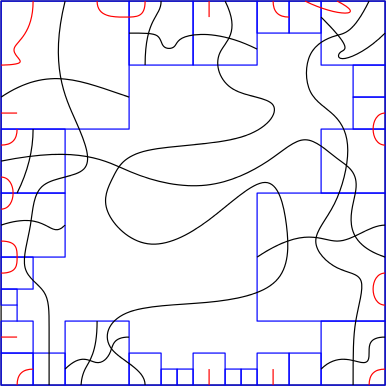}
\centering
\caption{The black cluster is $\clt_*(\cu)$. The cubes in blue are the good cubes at the boundary, which contains the small cluster $\slt(\cu)$ in color red. Its volume can be controlled by the integration along $\partial \cu$ of the size of the partition cubes.}
\end{figure}
\end{proof}

Thus, Lemmas \ref{lem:SmallClusterElementary} and \ref{lem:SmallCluster} answer \Cref{Q:SmallCluster1}, and the notation of $\slt(\cu)$ also helps us to solve \Cref{Q:SmallCluster2}:
\begin{lemma}\label{lem:SmallClusterFusion}
 For $\cu \in \Tcu$ such that $\size(\cu) > n > 0$, we have the estimate 
\begin{equation}\label{eq:SmallClusterFusion}
    \left\vert \clt_*(\cu) \backslash \left(\bigcup_{z \in 3^n \Zd \cap \cu} \clt_{*}(z + \cu_n) \right) \right\vert \Ind{\cu \in \Pcu_*}  \leq \GO_1(C \vert \cu \vert 3^{-n}).
\end{equation}
\end{lemma}
\begin{proof}
We decompose this difference in every cube of size $3^n$
\begin{align*}
& \left\vert \clt_*(\cu) \backslash \left(\bigcup_{z \in 3^n \Zd \cap \cu} \clt_{*}(z + \cu_n) \right) \right\vert \Ind{\cu \in \Pcu_*}\\
\leq & \sum_{z \in 3^n \Zd \cap \cu}  \left\vert (\clt_*(\cu) \cap (z + \cu_n)) \backslash \left(\bigcup_{z \in 3^n \Zd \cap \cu} \clt_{*}(z + \cu_n) \right) \right\vert \Ind{\cu \in \Pcu_*}\\
\leq & \sum_{z \in 3^n \Zd \cap \cu}  \left\vert (\clt_*(\cu) \cap (z + \cu_n)) \backslash  \clt_{*}(z + \cu_n) \right\vert \Ind{\cu \in \Pcu_*}\\
\leq & \sum_{z \in 3^n \Zd \cap \cu}  \left\vert (\clt_*(\cu) \cap (z + \cu_n)) \backslash  \clt_{*}(z + \cu_n) \right\vert \Ind{z + \cu_n \in \Pcu_*}\\
& \qquad  + \sum_{z \in 3^n \Zd \cap \cu}  \left\vert (\clt_*(\cu) \cap (z + \cu_n)) \backslash  \clt_{*}(z + \cu_n) \right\vert \Ind{z + \cu_n \notin \Pcu_*}\\
\end{align*}
The two terms can be treated separately. For the case $z + \cu_n \in \Pcu_*$, as we have mentioned, we have $\clp(z + \cu_n) = z + \cu_n$ and $\left\vert (\clt_*(\cu) \cap (z + \cu_n)) \backslash  \clt_{*}(z + \cu_n) \right\vert$ can be counted at the boundary $\partial( z + \cu_n)$.

 We turn this argument into the estimation using \cref{eq:sizeP} and \cref{eq:OSum}
\begin{align*}
& \sum_{z \in 3^n \Zd \cap \cu}  \left\vert (\clt_*(\cu) \cap (z + \cu_n)) \backslash  \clt_{*}(z + \cu_n) \right\vert \Ind{z + \cu_n \in \Pcu_*} \\
\leq &\sum_{z \in 3^n \Zd \cap \cu} |\slt(z + \cu_n)| \Ind{z + \cu_n \in \Pcu_*} \\
\leq &\GO_1(C|\cu|3^{-n}).
\end{align*}
For another part, we use \cref{eq:sizeP} and \cref{eq:OSum} directly that 
\begin{align*}
& \sum_{z \in 3^n \Zd \cap \cu}  \left\vert (\clt_*(\cu) \cap (z + \cu_n)) \backslash  \clt_{*}(z + \cu_n)  \Ind{z + \cu_n \notin \Pcu_*} \right\vert \\
\leq & \sum_{z \in 3^n \Zd \cap \cu}   |z + \cu_n|  \Ind{z + \cu_n \notin \Pcu_*}\\
\leq & \GO_1(C|\cu|3^{-n}).
\end{align*}
We combine all these estimates and conclude the result.
\end{proof}

Finally, we study \Cref{Q:SmallCluster3} on $(\clt_{\infty}(\cu) \cap \cu) \backslash \clt_*(\cu) $:
\begin{lemma}\label{lem:SmallClusterIntersect}
 Under the condition $\cu \in \Pcu_*$, and we use $\widetilde{\cu}$ to represent its predecessor, then we have $(\clt_{\infty} \cap \cu) = (\clt_*(\widetilde{\cu}) \cap \cu)$, and we have the estimate that 
$$
\left| (\clt_{\infty} \cap \cu) \backslash \clt_*(\cu) \right| \Ind{\cu \in \Pcu_*} \leq \GO_1(C|\cu|^{\frac{d-1}{d}}).
$$
\end{lemma}
\begin{proof}
 The lemma says when the cube $\cu_{m}$ is even better than a good cube, $\clt_*(\widetilde{\cu}) \cap \cu$ can contain all the part of $\clt_{\infty} \cap \cu$. One direction $(\clt_*(\widetilde{\cu}) \cap \cu) \subset (\clt_{\infty} \cap \cu)$ is obvious. We prove the other direction $ (\clt_{\infty} \cap \cu) \subset (\clt_*(\widetilde{\cu}) \cap \cu)$ by contradiction. We suppose that this direction is not correct so that there exists $z \in (\clt_{\infty} \cap \cu)$ but $z \notin (\clt_*(\widetilde{\cu}) \cap \cu)$. By \Cref{lem:SmallClusterElementary}, there exists a cluster $\clt'$ different from $\clt_*(\widetilde{\cu})$ and $z \in \clt'$ and $\clt'$ connects to $\partial \widetilde{\cu}$.($\cu \in \Pcu_* \Rightarrow \widetilde{\cu} \in \Pcu_*$.) Since $\clt_*(\cu)$ is part of $\clt_*(\widetilde{\cu})$, $\clt'$ cannot connect to $\clt_*(\cu)$. Thus, there exists an open path $\gamma$ such that $z \in \gamma \subset \clt'$ intersecting $\partial \cu$ and we  have $|\gamma| > \frac{1}{3}\size(\widetilde{\cu})$. This violate the second term in \Cref{prop:GeneralPartion} that a large path should belong to part of $\clt_*(\widetilde{\cu})$. We suppose that $\size(\cu) = 3^n$ and then apply \Cref{lem:SmallClusterFusion} to obtain that 
\begin{align*}
\left| (\clt_{\infty} \cap \cu) \backslash \clt_*(\cu) \right| \Ind{\cu \in \Pcu_*} &= \left| (\clt_*(\widetilde{\cu}) \cap \cu) \backslash \clt_*(\cu) \right| \Ind{\cu \in \Pcu_*} \\
& \leq  \left\vert \clt_*(\widetilde{\cu}) \backslash \left(\bigcup_{z \in 3^n \Zd \cap \widetilde{\cu}} \clt_{*}(z + \cu_n) \right) \right\vert \Ind{\cu \in \Pcu_*} \\
& \leq \GO_1(C |\cu| 3^{-n}). 
\end{align*}
\end{proof}
\begin{remark}
The same argument can prove even a stronger result that $\left(\clt_{\infty} \cap \frac{3}{4}\widetilde{\cu}\right) = \left(\clt_*(\widetilde{\cu}) \cap \frac{3}{4}\widetilde{\cu}\right)$.
\end{remark}

\section{Characterization of the effective conductance}\label{sec:Chara}
In the literature, there are several approaches to define the effective conductance $\ab$, and the object of this section is to give a proof of the equivalence of these definitions in the context of percolation. 

Let us at first recall the definition and some useful propositions in the previous work \cite[Definition 5.1]{armstrong2018elliptic}: we define the energy in the domain $U \subset \Zd$ with $l_p(x):= p \cdot x$ boundary condition
\begin{equation}\label{eq:defEnergy}
\begin{split}
\nu(U,p) := \inf_{v \in l_p + C_0(\clp(U))} \frac{1}{2\vert \clp(U) \vert} \bracket{\nabla v \cdot \a \nabla v}_{\Eda(\clt_*(U))},
\end{split}
\end{equation}
and we denote by $v(\cdot, U, p)$ its minimiser. The effective conductance $\ab$ is a deterministic positive scalar defined by
\begin{equation}\label{eq:defEffectiveConductance}
\frac{1}{2}p \cdot \ab p := \lim_{m \rightarrow \infty} \Expt[\nu(\cu_m,p)],
\end{equation}
with the rate of convergence \cite[Lemma 4.8]{armstrong2018elliptic}: there exists $s(d) > 0, \alpha(d, \pp, \Lambda) \in \left(0, \frac{1}{4} \right]$ and $C(d,p,\Lambda) < \infty$ such that for every $\cu \in \Tcu$
\begin{equation}\label{eq:ConvergenceEnergy}
\left\vert \frac{1}{2}p \cdot \ab p - \nu(\cu, p) \right\vert \leq \GO_s(C\vert p \vert^2 \size(\cu)^{-\alpha}).
\end{equation}
We will also use the following trivial bound several times in the proof
\begin{equation}\label{eq:TrivialBoundlp}
\nu(U,p) \leq  \frac{1}{2\vert \clp(U) \vert} \bracket{\nabla l_p \cdot \a \nabla l_p}_{\Eda(\clt_*(U))} \leq d \vert p \vert^2.
\end{equation}
The main theorem in this part is to prove the following characterization.
\begin{theorem}[Characterization of the effective conductance]\label{thm:CharaConductance}
In the context of homogenization in supercritical percolation, the effective conductance $\ab$ is a positive scalar constant and the following definitions are equivalent: 
\begin{align}
p \cdot \ab p &\stackrel{a.s}{=} \lim_{m \rightarrow \infty} \frac{1}{\vert \clp(\cu_m) \vert} \bracket{\nabla v(\cdot, \cu_m, p) \cdot \a \nabla v(\cdot, \cu_m, p)}_{\Eda(\clt_*(\cu_m))}. \label{eq:defab1} \\
p \cdot \ab p &\stackrel{a.s}{=} \lim_{m \rightarrow \infty} \frac{1}{\vert \cu_m \vert} \inf_{v \in l_p + C_0(\cu_m)} \bracket{\nabla v \cdot \a \nabla v}_{\Eda(\cu_m)}. \label{eq:defab2} \\
p \cdot \ab p &= \Expt[\D(\phi_p + l_p) \cdot \ac \D(\phi_p + l_p)].  \label{eq:defab3} \\
\ab p &= \Expt[\ac \D(\phi_p + l_p)].   \label{eq:defab4} 
\end{align}
\end{theorem}
Before starting the proof, we give some remarks on these definitions. Equation \eqref{eq:defab1} is just a variant of \cref{eq:defEffectiveConductance}. Equation \eqref{eq:defab2} differs from the first one in that just it does the minimization but does not enlarge the domain to $\clp(\cu_m)$ nor restricts the problem to $\clt_*(\cu_m)$. Equation \eqref{eq:defab3} uses the linear $\a$-harmonic function in the whole space instead of that in $\cu_m$, so it is stationary. The last one is a little different from the previous three ones, but we need it in \Cref{prop:SpatialAverage2}, thus we add it to the list of equivalent definitions.

\begin{proof}
Equation \eqref{eq:defab1} is a direct consequence of \cref{eq:defEffectiveConductance} and \cref{eq:ConvergenceEnergy}, Markov's inequality and the lemma of Borel-Cantelli to transform it to an ``almost sure" version.

Equation \eqref{eq:defab2} is a variant from the first one, especially when $\cu_m \in \Pcu_*$ they are very close. So we do the decomposition 
\begin{align*}
 & \left\vert \frac{1}{\vert \cu_m \vert} \inf_{v \in l_p + C_0(\cu_m)} \bracket{\nabla v \cdot \a \nabla v}_{\Eda(\cu_m)} - p\cdot \ab p \right\vert \\
\leq  &  \left\vert \frac{1}{\vert \cu_m \vert} \inf_{v \in l_p + C_0(\cu_m)} \bracket{\nabla v \cdot \a \nabla v}_{\Eda(\cu_m)} - p\cdot \ab p \right\vert \Ind{\cu_m \in \Pcu_*} \\
&\qquad +  \left\vert \frac{1}{\vert \cu_m \vert} \inf_{v \in l_p + C_0(\cu_m)} \bracket{\nabla v \cdot \a \nabla v}_{\Eda(\cu_m)} - p\cdot \ab p \right\vert \Ind{\cu_m \notin \Pcu_*} \\
\end{align*}
and the second one can be handled easily by a trivial bound by comparing with $l_p$ as in \cref{eq:TrivialBoundlp}
\begin{align*}
\left\vert \frac{1}{\vert \cu_m \vert} \inf_{v \in l_p + C_0(\cu_m)} \bracket{\nabla v \cdot \a \nabla v}_{\Eda(\cu_m)} - p\cdot \ab p \right\vert \Ind{\cu_m \notin \Pcu_*} \leq \GO_1(C(d,\pp,\Lambda)\vert p \vert^2 3^{-m}). 
\end{align*}
By an argument of Borel-Cantelli, we prove this term converges almost surely to $0$.
Then we focus on the case $\cu_m \in \Pcu_*$, In fact, in this case the minimiser on $\Eda(\cu_m)$ is the sum of the one on each clusters. Observing that the one on isolated cluster from $\partial \cu_m$ can be null since it has no boundary condition, so we have to deal with the one on $\clt_*(\cu_m)$ and the one on the small cluster $\slt(\cu_m)$. We apply \cref{eq:SmallClusterSize}, the estimate \cref{eq:sizeP} and a trivial bound \cref{eq:TrivialBoundlp} to get
\begin{align*}
& \left\vert \frac{1}{\vert \cu_m \vert} \inf_{v \in l_p + C_0(\cu_m)} \bracket{\nabla v \cdot \a \nabla v}_{\Eda(\cu_m)} - p\cdot \ab p \right\vert \Ind{\cu_m \in \Pcu_*} \\
\leq &  \left\vert \frac{1}{\vert \clp(\cu_m) \vert} \bracket{\nabla v(\cdot, \cu_m, p) \cdot \a \nabla v(\cdot, \cu_m, p)}_{\Eda(\clt_*(\cu_m))} - p\cdot \ab p \right\vert \\
& \qquad +  \left\vert \frac{1}{\vert \clp(\cu_m) \vert} \bracket{\nabla v(\cdot, \cu_m, p) \cdot \a \nabla v(\cdot, \cu_m, p)}_{\Eda(\clt_*(\cu_m))} - p\cdot \ab p \right\vert \Ind{\cu_m \notin \Pcu_*} \\
& \qquad + \left\vert \frac{1}{\vert \cu_m \vert} \inf_{v \in l_p + C_0(\cu_m)} \bracket{\nabla v \cdot \a \nabla v}_{\Eda(\slt(\cu_m))} \right\vert 
\Ind{\cu_m \in \Pcu_*} \\
\leq & \left\vert \frac{1}{\vert \clp(\cu_m) \vert} \bracket{\nabla v(\cdot, \cu_m, p) \cdot \a \nabla v(\cdot, \cu_m, p)}_{\Eda(\clt_*(\cu_m))} - p\cdot \ab p \right\vert  \\
& \qquad +  C(d)\vert p \vert^2  \Ind{\size(\cu_{\Pcu}(0)) > 3^m} +  C(d)\vert p \vert^2  \frac{\vert \slt(\cu_m) \vert \Ind{\cu_m \in \Pcu_*}}{\vert \cu_m \vert} \\
\leq & \left\vert \frac{1}{\vert \clp(\cu_m) \vert} \bracket{\nabla v(\cdot, \cu_m, p) \cdot \a \nabla v(\cdot, \cu_m, p)}_{\Eda(\clt_*(\cu_m))} - p\cdot \ab p \right\vert + \GO_1(C(d,\pp,\Lambda)\vert p \vert^2 3^{-m}). \\
\end{align*}
So its almost sure limit is the same as the first one when $m \rightarrow \infty$.

By a similar calculation, one can prove a variant of \cref{eq:defab2} that reads 
\begin{equation}
p \cdot \ab p \stackrel{a.s}{=} \lim_{m \rightarrow \infty} \frac{1}{\vert \cu_m \vert} \inf_{v \in l_p + C_0(\cu_m)} \bracket{\nabla v \cdot \ac \nabla v}_{\Eda(\cu_m)},
\end{equation}
and we recall, see \cite[Theorem 9.1]{jikov2012homogenization}, that this definition coincides with \cref{eq:defab3}. By a calculus of variation argument, we have 
$$
\forall p, q \in \Rd, \qquad q \cdot \ab p = \Expt[\D(\phi_q + l_q) \cdot \ac \D(\phi_p + l_p)].
$$ 
Moreover, observing that $\Ind{\a \neq 0} \D \phi_q + q$ and $\ac (\D\phi_p + p)$ are stationary, and the former is a the gradient and the latter is divergence free, we can use the Div-Curl and Birkhoff theorems 
$$
q \cdot \ab p = \Expt[\Ind{\a \neq 0} \D \phi_q + q]\Expt[\ac (\D\phi_p + p)] = q \cdot \Expt[\ac (\D \phi_p + p)].
$$
This concludes the equivalence with \cref{eq:defab4}.
\end{proof}

\subsection*{Acknowledgments}
I am grateful to Jean-Christophe Mourrat for his suggestion to study this topic and helpful discussions, and Paul Dario for inspiring discussions.


\bibliographystyle{abbrv}
\bibliography{AlgoPercoRef}

\end{document}